\documentclass[11pt]{article}
\usepackage[english]{babel}
\usepackage[latin1]{inputenc}
\usepackage{amsmath}
\usepackage{amssymb}
\usepackage{amsfonts}
\usepackage{amsthm}
\usepackage{color}

\newcommand{\p}{\mathfrak{p}}
\newcommand{\q}{\mathfrak{q}}
\newcommand{\s}{\mathfrak{s}}
\renewcommand{\u}{\mathfrak{u}}
\renewcommand{\v}{\mathfrak{v}}
\newcommand{\C}{\mathbb{C}}
\newcommand{\N}{\mathbb{N}}

\newcommand{\R}{\mathbb{R}}
\renewcommand{\S}{\mathbb{S}}
\newcommand{\T}{\mathbb{T}}

\newcommand{\boA}{\mathcal{A}}

\newcommand{\boE}{\mathcal{E}}
\newcommand{\boF}{\mathcal{F}}

\newcommand{\boH}{\mathcal{H}}

\newcommand{\boP}{\mathcal{P}}

\newcommand{\boS}{\mathcal{S}}

\newcommand{\boX}{\mathcal{X}}
\newcommand{\boZ}{\mathcal{Z}}
\newcommand{\eps}{\varepsilon}
\newcommand{\Card}{{\rm Card}}
\newcommand{\ch}{{\rm ch}}

\renewcommand{\div}{{\rm div}}

\renewcommand{\Im}{{\rm Im}}
\newcommand{\on}{\ {\rm on} \ }
\renewcommand{\Re}{{\rm Re}}

\newcommand{\sign}{{\rm sign}}
\renewcommand{\th}{{\rm th}}

\newtheorem{cor}{Corollary}
\newtheorem{lemma}{Lemma}
\newtheorem{prop}{Proposition}
\newtheorem{step}{Step}
\newtheorem{theorem}{Theorem}
\theoremstyle{definition}
\newtheorem{remark}{Remark}

\begin{document}

\title{Existence and properties of travelling waves for the Gross-Pitaevskii equation}
\author{\renewcommand{\thefootnote}{\arabic{footnote}}
Fabrice B\'ethuel \footnotemark[1], Philippe Gravejat \footnotemark[2], Jean-Claude Saut \footnotemark[3]}
\footnotetext[1]{Laboratoire Jacques-Louis Lions, Universit\'e Pierre et Marie Curie, Bo\^ite Courrier 187, 75252 Paris Cedex 05, France. E-mail: bethuel@ann.jussieu.fr}
\footnotetext[2]{Centre de Recherche en Math\'ematiques de la D\'ecision, Universit\'e Paris Dauphine, Place du Mar\'echal De Lattre De Tassigny, 75775 Paris Cedex 16, France. E-mail: gravejat@ceremade.dauphine.fr}
\footnotetext[3]{Laboratoire de Math\'ematiques, Universit\'e Paris Sud, B\^atiment 425, 91405 Orsay Cedex, France. E-mail: Jean-Claude.Saut@math.u-psud.fr}
\date{}
\maketitle

\begin{abstract}
This paper presents recent results concerning the existence and qualitative properties of travelling wave solutions to the Gross-Pitaevskii equation posed on the whole space $\R^N$. Unlike the defocusing nonlinear Schr\"odinger equations with null condition at infinity, the presence of non-zero conditions at infinity yields a rather rich and delicate dynamics. We focus on the case $N = 2$ and $N = 3$, and also briefly review some classical results on the one-dimensional case. The works we survey provide rigorous justifications to the impressive series of results which Jones, Putterman and Roberts \cite{JoneRob1,JonPuRo1} established by formal and numerical arguments.
\end{abstract}

\section{Introduction}

The Gross-Pitaevskii equation
\renewcommand{\theequation}{GP}
\begin{equation}
\label{GP}
i \partial_t \Psi = \Delta \Psi + \Psi (1 - |\Psi|^2) \on \R^N \times \R,
\end{equation} 
appears as a relevant model in various areas of physics: nonlinear optics, fluid mechanics, Bose-Einstein condensation... (see, for instance, \cite{Pitaevs1,Gross1,IordSmi1,JoneRob1,JonPuRo1,KivsLut1, BethSau2}). It corresponds to a version of the defocusing nonlinear Schr\"odinger equations. If one considers finite energy solutions, as suggested by the formal conservation of the energy (see \eqref{GLE} below), then $\Psi$ should not vanish at infinity, but instead $|\Psi(x, \cdot)|$ should in some sense tend to $1$ when $|x| \to + \infty$. As a matter of fact, this condition ensures that \eqref{GP} has a non-trivial dynamics, contrary to the case of null condition at infinity, where the dynamics is expected to be trivial (dispersion, scattering...).

In nonlinear optics, \eqref{GP} classically appears in the context of optical dark solitons, that is localized nonlinear waves (or "holes") which exist on a stable continuous wave background (see the review paper by Kivshar and Luther-Davies \cite{KivsLut1} or paper \cite{KivPeSt1}, where different scenarios for the transverse instability of one-dimensional black solitons are proposed in the context of \eqref{GP}). The boundary condition $|\Psi(x, \cdot)| \to 1$ is due to this non-zero background.

The Gross-Pitaevskii equation has also been intensively used as a model for superfluid Helium II and for Bose-Einstein condensation. We refer to \cite{Gross1,Coste1}, and to the survey by Berloff \cite{Berloff1}.

At least on a formal level, the Gross-Pitaevskii equation is hamiltonian. The conserved Hamiltonian is a Ginzburg-Landau energy, namely
\renewcommand{\theequation}{\arabic{equation}}
\numberwithin{equation}{section}
\setcounter{equation}{0}
\begin{equation}
\label{GLE}
E(\Psi) = \frac{1}{2} \int_{\R^N} |\nabla \Psi|^2 + \frac{1}{4} \int_{\R^N} (1 - |\Psi|^2)^2 \equiv \int_{\R^N} e(\Psi).
\end{equation}
Similarly, the momentum
\begin{equation}
\label{VectP}
\vec{P}(\Psi) = \frac{1}{2} \int_{\R^N} \langle i \nabla \Psi \ , \Psi - 1 \rangle
\end{equation}
is formally conserved. We will denote by $p$, the first component of $\vec{P}$, which is hence a scalar. Another quantity which is formally conserved by the flow is
$$m(\Psi) = \frac{1}{2} \int_{\R^N} \Big( |\Psi|^2 - 1 \Big).$$

If $\Psi$ does not vanish, one may write
$$\Psi = \sqrt{\rho} \exp i \varphi.$$
This leads to the hydrodynamic form of the equation
\begin{equation}
\label{HDGP}
\left\{ \begin{array}{ll} \partial_t \rho + \div (\rho v) = 0, \\ \rho (\partial_t v + v. \nabla v) + \nabla \rho^2 = \rho \nabla \left( \frac{\Delta \rho}{\rho} - \frac{|\nabla \rho|^2}{2 \rho^2} \right),\end{array} \right.
\end{equation}
where
$$v = - 2 \nabla \varphi.$$

If one neglects the right-hand side of the second equation, which is often referred to as the quantum pressure, system \eqref{HDGP} is similar to the Euler equation for a compressible fluid, with pressure law $p(\rho) = \rho^2$. In particular, the speed of sound waves near the constant solution $\Psi = 1$, that is $\rho = 1$ and $v = 0$, is given by
$$c_s = \sqrt{2}.$$

When $N = 1$, \eqref{GP} is integrable by the inverse scattering method, and it has been formally analyzed (with the non-trivial boundary condition at infinity) in \cite{ShabZak2}. The (local and global) Cauchy problem associated to \eqref{GP} with non-zero condition at infinity has been recently thoroughly investigated (see the survey by G\'erard \cite{Gerard2} in this volume). The possible existence of non-trivial travelling wave solutions highly suggests that the long time behaviour of global unsteady solutions should be non-trivial (as it would be in the presence of zero condition at infinity).

Travelling waves are solutions to \eqref{GP} of the form
$$\Psi(x,t) = u(x_1 - ct, x_\perp), \ x_\perp = (x_2, \ldots, x_N).$$
Here, the parameter $c \in \R$ corresponds to the speed of the travelling waves (we may restrict to the case $c \geq 0$ using complex conjugation). The equation for the profile $u$ is given by
\renewcommand{\theequation}{TWc}
\begin{equation}
\label{TWc}
i c \partial_1 u + \Delta u + u (1 - |u|^2) = 0.
\end{equation}

Travelling waves of finite energy play an important role in the dynamics of the Gross-Pitaevskii equation. Whereas the one-dimensional case can be integrated explicitly, solutions are not explicitly known in higher dimensions.

The first attempts to prove existence of travelling waves and to study their properties were based on formal expansions and numerics. They were mostly performed by physicists. The three-dimensional case was considered by Iordanskii and Smirnov \cite{IordSmi1}, whereas Jones, Putterman and Roberts \cite{JoneRob1,JonPuRo1} developed a thorough analysis in dimensions two and three: they found a branch of solutions with speeds $c$ covering the full subsonic range $(0, \sqrt{2})$, and they conjectured the non-existence of travelling waves for supersonic speeds. Moreover, they showed the connections with solitary waves for the Kadomtsev-Petviashvili equation \eqref{KP}.

More recently, the program of Jones, Putterman and Roberts \cite{JoneRob1,JonPuRo1} has been provided with rigorous mathematical proofs. One of the purposes of these notes is to review some of these mathematical progresses and to stress a number of open problems.

The next section provides a rather extensive survey of the one-dimensional case. We review classical results and also prove some new ones. Sections \ref{qualibat} and \ref{existence} are devoted to recent progresses on the rigorous mathematical proofs of the program of Jones, Putterman and Roberts \cite{JoneRob1,JonPuRo1}: qualitative properties and existence of travelling waves in higher dimensions. Section \ref{potaufeu} presents some related problems. In particular, we establish new results on the nonlinear Schr\"odinger flow around an obstacle, modelled by a potential.

\section{The one-dimensional case}
\label{formule1}

\subsection{The integration of equation \eqref{TWc}}
\label{testarossa}

In the one-dimensional case, equation \eqref{TWc} is entirely integrable using standard arguments from ordinary differential equation theory. Hence, it is possible to classify all the travelling waves with finite energy according to their speed $c$.

\begin{theorem}[\cite{Tsuzuki1,Maris3,Graveja4,Gallo2,DiMeGal1,Berloff1}]
\label{ferrari}
Assume $N = 1$ and $c \geq 0$, and let $v$ be a solution of finite energy to \eqref{TWc}.\\
i) If $c \geq \sqrt{2}$, $v$ is a constant of modulus one.\\
ii) If $0 \leq c < \sqrt{2}$, up to a multiplication by a constant of
modulus one and a translation, $v$ is either identically equal to $1$,
or to
\renewcommand{\theequation}{\arabic{equation}}
\numberwithin{equation}{section}
\setcounter{equation}{3}
\begin{equation}
\label{schumi}
v(x) = \v_c(x) \equiv \sqrt{\frac{2 - c^2}{2}} \th \Big( \frac{\sqrt{2 - c^2}}{2} x \Big) - i \frac{c}{\sqrt{2}}.
\end{equation}
Moreover, $\v_c$ may be written as
\begin{equation}
\label{hakkinen}
\v_c(x) = \sqrt{1 - \frac{2 - c^2}{2 {\rm ch}^2 \bigg( \frac{ \sqrt{2 - c^2}}{2} x \bigg)}} \exp i \Bigg( \arctan \bigg( \frac{\sqrt{2 - c^2}}{c} \th \Big( \frac{\sqrt{2 - c^2}}{2} x \Big) \bigg) - \frac{\pi}{2} \Bigg),
\end{equation}
if $0 < c < \sqrt{2}$.
\end{theorem}

For sake of completeness, we briefly recall the proof of Theorem \ref{ferrari}.

\begin{proof}
Denoting $v = v_1 + i v_2$, equation \eqref{TWc} becomes
\begin{align}
\label{massa}
& v_1'' - c v_2' + v_1 (1 - v_1^2 - v_2 ^2) = 0,\\
\label{raikonen}
& v_2'' + c v_1' + v_2 (1 - v_1^2 - v_2 ^2) = 0.
\end{align}
Multiplying \eqref{massa} by $v_2$ and \eqref{raikonen} by $v_1$ yields
\begin{equation}
\label{todt}
(v_1 v_2' - v_2 v_1')' = \frac{c}{2} \eta',
\end{equation}
where $\eta \equiv 1 - |v|^2$. Since $v$ is a finite energy solution to \eqref{TWc}, it is smooth and bounded on $\R$. Moreover, $v'(x)$ and $\eta(x)$ tend to $0$, as $x \to \pm \infty$. In particular, the integration of \eqref{todt} leads to
\begin{equation}
\label{alesi}
v_1 v_2' - v_2 v_1' = \frac{c}{2} \eta.
\end{equation}
Multiplying \eqref{massa} by $v_1'$, and \eqref{raikonen} by $v_2'$, and integrating the resulting equation, we also deduce that
\begin{equation}
\label{barrichelo}
|v'|^2 = \frac{\eta^2}{2}.
\end{equation}
Computing
$$\eta'' = - 2 |v'|^2 - 2 (v_1 v_1''+ v_2 v_2'') = - 2 |v'|^2 - 2c
(v_1 v_2' - v_2 v_1') + 2 \eta - 2 \eta^2,$$
where we invoke \eqref{TWc} for the second identity, we finally obtain from \eqref{alesi} and \eqref{barrichelo},
\begin{equation}
\label{scuderia}
\eta'' + (c^2 - 2) \eta + 3 \eta^2 = 0.
\end{equation}
Multiplying \eqref{scuderia} by $\eta'$ and integrating, it follows that
\begin{equation}
\label{imola}
{\eta'}^2 + (c^2 - 2) \eta^2 + 2 \eta^3 = 0.
\end{equation}
As previously mentioned, $\eta$ is a smooth, bounded function on $\R$, such that $\eta(x) \to 0$, as $x \to \pm \infty$. In case it is not identically equal to $0$ on $\R$, we can assume, up to a translation, that
\begin{equation}
\label{madmax1}
|\eta(0)| = \max \big\{ |\eta(x)|, x \in \R \big\} > 0,
\end{equation}
so that
\begin{equation}
\label{madmax2}
\eta'(0) = 0.
\end{equation}
Using \eqref{imola}, we are led to
\begin{equation}
\label{madmax3}
\eta(0) = \frac{2 - c^2}{2}.
\end{equation}
In particular, if $c = \sqrt{2}$, $\eta(0) = 0$, which gives a contradiction with \eqref{madmax1}. Hence, $\eta$ is identically equal to $0$ if $c = \sqrt{2}$. On the other hand, if $c > \sqrt{2}$, then, by \eqref{imola},
$$(c^2 - 2 + 2 \eta) \eta^2 = - {\eta'}^2 \leq 0.$$
so that, by \eqref{madmax3} and the continuity of $\eta$,
$$\eta(x) \leq \frac{2 - c^2}{2} < 0, \ \forall x \in \R.$$
This gives a contradiction with the fact that $\eta \in L^2(\R)$. Therefore, $\eta$ is identically equal to $0$ if $c > \sqrt{2}$. Finally, if $0 \leq c < \sqrt{2}$, it follows from \eqref{madmax2}, \eqref{madmax3}, and Cauchy-Lipschitz's theorem that there exists a unique local solution $\eta$ to \eqref{scuderia}-\eqref{madmax2}-\eqref{madmax3} on some neighbourhood of the point $0$. One can check that this solution is given by
\begin{equation}
\label{maranello}
\eta(x) = \frac{2 - c^2}{2 {\rm ch}^2 \Big( \frac{\sqrt{2 - c^2}}{2} x \Big)},
\end{equation}
so that it can be extended to be the unique, up to some translation, non-constant solution to \eqref{scuderia} in $L^2(\R)$. In conclusion, either $\eta$ is identically equal to $0$ (in particular, when $c \geq \sqrt{2}$), either it is given, up to a translation, by \eqref{maranello}.

When $c \neq 0$, it follows that $|v| = \sqrt{1 - \eta}$ does not vanish on $\R$, so that we can construct a lifting $\varphi$ of $v$, that is a smooth function $\varphi$ such that
$$v = |v| \exp i \varphi = \sqrt{1 - \eta} \exp i \varphi.$$
Invoking \eqref{alesi}, we are led to
\begin{equation}
\label{sanmarino}
\varphi' = \frac{c \eta}{2 - 2 \eta},
\end{equation}
so that $\varphi \equiv \varphi_0$ is identically constant when $\eta = 0$, whereas $\varphi$ is given by
$$\varphi(x) = \varphi_0 + \arctan \bigg( \frac{\sqrt{2 - c^2}}{c} \th \Big( \frac{\sqrt{2 - c^2}}{2} x \Big) \bigg),$$
when $\eta$ is given by \eqref{maranello}. This completes the proofs of assertion i), and of identity \eqref{hakkinen}. Identity \eqref{schumi} follows from standard trigonometric identities, in particular, from the fact that
$$\cos \big( \arctan(x) \big) = \frac{1}{\sqrt{1 + x^2}}, \ {\rm and} \ \sin \big( \arctan(x) \big) = \frac{x}{\sqrt{1 + x^2}}.$$

When $c = 0$, either $\eta$ is identically equal to $0$, so that, by \eqref{barrichelo}, $v$ is a constant of modulus one, or $\eta$ is given, up to some translation, by
$$\eta(x) = \frac{1}{\ch^2 \Big( \frac{x}{\sqrt{2}} \Big)}.$$
In this case, $v(0) = 0$, so that, by \eqref{barrichelo}, we may assume, up to the multiplication by a constant of modulus one, that
$$v'(0) = \frac{1}{\sqrt{2}}.$$
Hence, by Cauchy-Lipschitz's theorem, there exists a unique local solution to \eqref{TWc}, with $c = 0$, on some neighbourhood of the point $0$. One can check that this solution is given by
$$\v_0(x) = \th \Big( \frac{x}{\sqrt{2}} \Big),$$
so that, up to the invariances,
$$v(x) = \v_0(x) = \th \Big( \frac{x}{\sqrt{2}} \Big).$$
This completes the proof of Theorem \ref{ferrari}.
\end{proof}

Let us next draw some remarks concerning qualitative properties of non-constant finite energy travelling waves, which we deduce from formulae \eqref{schumi} and \eqref{hakkinen}.

The first observation is that the one-dimensional non-constant travelling waves form a smooth branch of subsonic solutions to \eqref{TWc}. As seen in Theorem \ref{ferrari}, there exist neither sonic, nor supersonic non-constant travelling waves. Moreover, formula \eqref{schumi} yields the spatial asymptotics of the non-constant solutions to \eqref{TWc}. Notice in particular that
$$\v_c(x) \to \v_c^{\pm \infty} \equiv \pm \sqrt{1 - \frac{c^2}{2}} - i \frac{c}{\sqrt{2}}, \ {\rm as} \ x \to \pm \infty.$$
Hence, $\v_c(x)$ converges to a constant $\v_c^{\pm \infty}$ of modulus one, as $x \to \pm \infty$. This constant is different for the limit $x \to - \infty$, and the limit $x \to + \infty$. In contrast, the solution we will find in the higher dimensional case will have a limit at infinity which is independent on the direction (see Theorem \ref{Decay}). Notice also that the function $\v_c - \v_c^{\pm \infty}$ has exponential decay, whereas the decay is algebraic in higher dimensions (see also Theorem \ref{Decay}). This fact is related to the properties of the kernel of the linear part of the equation. The Fourier transform of the kernel $K$ associated to \eqref{TWc} (more precisely, to \eqref{scuderia}) is given by
$$\widehat{K}(\xi) = \frac{1}{2 - c^2 + \xi^2},$$
hence has some exponential decay, when $0 \leq c < \sqrt{2}$.

A second observation is, in view of identity \eqref{hakkinen}, that $\v_c(x)$ does not vanish, unless $c = 0$. In the case $c = 0$, the non-constant solution of finite energy to \eqref{TWc} is given, up to the invariances, by
$$\v_0(x) = \th \Big( \frac{x}{\sqrt{2}} \Big),$$
which vanishes at $x = 0$. This stationary solution, which is real-valued, is often termed a kink solution. In contrast with the one-dimensional case, it can be shown in higher dimensions, using Pohozaev's identity, that the only finite energy solutions with $c = 0$ are constants (see \cite{BethSau1}). The energy of the kink solution is given by
$$E(\v_0) = \frac{2 \sqrt{2}}{3}.$$
The kink solution has the following remarkable minimization property.

\begin{lemma} 
\label{kinkmini}
We have
$$E(\v_0) =\inf \Big\{ E(v), v \in H^1_{\rm loc}(\R), \underset{x \in \R}{\inf} \big| v(x) \big| = 0 \Big\}.$$
In particular, if $E(v) < \frac{2 \sqrt{2}}{3}$, then,
$$\underset{x \in \R}{\inf} \big| v(x) \big| > 0.$$
\end{lemma}

\begin{proof}
We consider a minimizing sequence $(v_n)_{n \in \N}$ for the minimization problem
$$\boE_0 = \inf \bigg\{ \int_0^{+ \infty} e(v), v \in H^1_{\rm loc}([0, + \infty)), v(0) = 0 \bigg\},$$
which is well-defined by Sobolev's embedding theorem. We notice that the functions $v_n'$ are uniformly bounded in $L^2([0, + \infty))$, and that $v_n(0) = 0$. Hence, by Rellich's compactness theorem, there exists some function $u \in H^1_{\rm loc}([0, + \infty))$, with $u(0) = 0$, such that, up to a subsequence,
$$v_n' \rightharpoonup u' \ {\rm in} \ L^2([0, + \infty)), \ {\rm and} \ v_n \to u \ {\rm in} \ L^\infty_{\rm loc}([0, + \infty)), \ {\rm as} \ n \to + \infty.$$
By Fatou's lemma, we are led to
$$\int_0^{+ \infty} e(u) = \frac{1}{2} \int_0^{+ \infty} \big( u' \big)^2 + \frac{1}{4} \int_0^{+ \infty} \underset{n \to + \infty}{\liminf} \Big( \big( 1 - |v_n|^2 \big)^2 \Big) \leq \underset{n \to + \infty}{\liminf} \bigg( \int_0^{+ \infty} e(v_n) \bigg),$$
so that the infimum $\boE_0$ is achieved by the function $u$. In particular, the solution $u$ is critical for the Ginzburg-Landau energy, i.e. solves
$$u'' + u (1 - |u|^2) = 0.$$
Integrating this equation yields $u(x) = \v_0(x) = \th \Big( \frac{x}{\sqrt{2}} \Big)$ (see the proof of Theorem \ref{ferrari}), so that
$$\boE_0 = \int_0^{+ \infty} e(\v_0) = \frac{\sqrt{2}}{3}.$$
Next consider a map $v \in H^1_{\rm loc}(\R)$, with finite Ginzburg-Landau energy, and which vanishes at some point $x_0$. In view of the invariance by translation, we may assume that $x_0 = 0$, whereas the minimality of $\v_0$ yields
$$\int_0^{+ \infty} e(v) \geq \frac{\sqrt{2}}{3},$$
and the same inequality holds for the energy on $(- \infty, 0]$, so that
$$E(v) \geq \frac{2 \sqrt{2}}{3} = E(\v_0).$$
The proof of Lemma \ref{kinkmini} follows.
\end{proof}

\begin{remark}
\label{botox}
When $E(v) < \frac{2 \sqrt{2}}{3}$, $v$ has no zero, so that we may write
$$v = \varrho \exp i \varphi.$$
In particular, the energy density of $v$ can be expressed as
\begin{equation}
\label{thenergy}
e(v) = \frac{1}{2} \Big( (\varrho')^2 + \varrho^2 (\varphi')^2 \Big)+ \frac{1}{4} \Big( 1 - \varrho^2 \Big)^2,
\end{equation}
whereas, for the scalar momentum, we have
\begin{equation}
\label{scalare}
\langle i v', v \rangle = - \varrho^2 \varphi'.
\end{equation}
\end{remark}

\subsection{Energy and momentum of travelling wave solutions}
\label{energetic}

A general principle states that travelling wave solutions are critical points of the energy, keeping the momentum fixed. A natural question is therefore to determine whether this principle applies to the solutions $\v_c$. However, the precise definition of the momentum raises a serious difficulty.

In the context of nonlinear Schr\"odinger equations, the momentum of solutions should be defined as
$$P(\Psi) = \frac{1}{2} \int_\R \langle i \Psi', \Psi \rangle.$$
This quantity has a physical meaning and at least formally, describes the evolution center of mass in the \eqref{GP} dynamics. This fact is shown for instance by the relation
\begin{equation}
\label{interP}
\frac{d}{dt} \bigg( \frac{1}{2} \int_{\R} x \Big( |\Psi(x,t)|^2 - 1 \Big) dx \bigg) = 2 P(\Psi(t)).
\end{equation}
The momentum $P(\v_c)$ is well-defined in view of the exponential decay of the quantity $\langle i \v_c', \v_c \rangle$, and is equal to
$$P(\v_c) = - \frac{c \sqrt{2 - c^2}}{2}.$$
Notice that, if $c \neq 0$, then we may write $\v_c = \varrho_c \exp i \varphi_c$, so that we have 
$$P(\v_c) = - \frac{1}{2} \int_\R\varrho_c^2 \varphi_c'.$$
In the variational context we have in mind, we need to extend this definition to a larger class of maps $v$. S We restrict ourselves to the case the case $E(v) < \frac{2 \sqrt{2}}{3}$, since we wish to establish a minimality property for the solution $\v_c$, which is known to have an energy satisfying $E(\v_c) < \frac{2 \sqrt{2}}{3}$. Assuming $E(v) < \frac{2 \sqrt{2}}{3}$, it follows from Lemma \ref{kinkmini} that $v$ has no zero and the continuity of $v$ combined with the finiteness of its energy yield 
$$\delta = \underset{x \in \R}{\inf} |v(x)| > 0,$$
so that we may write $v = \varrho \exp i \varphi$. However, the quantity
$$P(v) = -\frac{1}{2} \int_\R \varrho^2 \varphi' = \frac{1}{2} \int_\R (1 - \varrho^2) \varphi' - \frac{1}{2} \int_\R \varphi',$$
is only well-defined if we are able to give a meaning to the last integral $\int_\R \varphi'$, since the first one may be bounded by Cauchy-Schwarz's inequality,
$$\bigg| \int_\R (1 - \varrho^2) \varphi' \bigg| \leq \frac{1}{2} \int_\R (1 - \varrho^2)^2 + \frac{1}{2 \delta^2} \int_\R \varrho^2 (\varphi')^2.$$
The integral $\int_\R \varphi'$ is well-defined in particular if the limits $\varphi(\pm \infty) = \underset{x \to \pm \infty}{\lim} \varphi(x)$ do exist, in which case we have $\int_\R \varphi' = \varphi(+ \infty) - \varphi(- \infty)$. This leads us to consider the set
$$\boZ_0(\R) = \Big\{ v \in H^1_{\rm loc}(\R), \ {\rm s.t.} \ E(v) < \frac{2 \sqrt{2}}{3}, \ {\rm and} \ v_{\pm \infty} = \lim_{x \to \pm \infty} v(x) \ {\rm exist} \Big\}.$$
In view of the above discussion, the momentum $P(v)$ is well-defined for any $v \in \boZ_0(\R)$. However, we have

\begin{lemma}
\label{alonso}
Let $\p > 0$. We have
$$\inf \Big\{ E(v), v \in \boZ_0(\R), \ {\rm s.t.} \ P(v) = \p \Big\} = 0.$$
\end{lemma}

\begin{proof}
Let $n \geq 1$, and consider the map $v_n = \exp i \psi_n \in H^1_{\rm loc}(\R)$ defined by
$$\psi_n(x) = \left\{ \begin{array}{lll} 0, \ {\rm if} \ x \leq - n,\\ - \p \frac{x + n}{n}, \ {\rm if} \ - n \leq x \leq n,\\
- 2 \p, \ {\rm if} \ x \geq n. \end{array} \right.$$
Since
\begin{equation}
\label{alonso2}
E(v_n) = \frac{\p^2}{n} \to 0, \ {\rm as} \ n \to + \infty,
\end{equation}
$v_n(x) \to 1$, as $x \to - \infty$, and $v_n(x) \to \exp (- 2 i \p)$, as $x \to + \infty$, the map $v_n$ belongs to $\boZ_0(\R)$, for any $n$ sufficiently large. Moreover, we compute
$$P(v_n) = \p,$$
so that Lemma \ref{alonso} follows from \eqref{alonso2} taking the limit $n \to + \infty$.
\end{proof}

We do not actually know if the space $\boZ_0(\R)$, and the limits at infinity for maps in that space are preserved by the \eqref{GP} flow. However, this is the case if one imposes more regularity on the initial condition. Consider, for instance, the space
$$\boZ_2(\R) = \Big\{ v \in \boZ_0(\R), \ {\rm s.t.} \ v_{|(0, \pm \infty)} \in \{ v_{\pm \infty} \} + H^2(0, \pm \infty) \Big\}.$$
We have

\begin{lemma}
\label{grouic}
Assume $\Psi_0 \in \boZ_2(\R)$, and let $\Psi$ be the solution to \eqref{GP} with initial datum $\Psi_0$. Then, $\Psi(t)$ belongs to $\boZ_2(\R)$ for any $t \in \R$, and $\Psi_{\pm \infty}(t) = \Psi_{0, \pm \infty}$.
\end{lemma}

We will not provide a proof of Lemma \ref{grouic} in this survey.

Lemmas \ref{alonso} and \ref{grouic} lead us to introduce a notion of renormalized momentum $p$, whose physical meaning seems less obvious, but whose variational properties are more relevant for our study. We set, for any function $v \in H^1_{\rm loc}(\R)$ such that $E(v) < \frac{2 \sqrt{2}}{3}$,
\begin{equation}
\label{prost}
p(v) \equiv \frac{1}{2} \int_{\R} (1 - \varrho^2) \varphi',
\end{equation}
where we use the polar form $v = \varrho \exp i \varphi$ of $v$. Notice that, by construction,
$$P(v) = p(v) - \frac{1}{2} \Big( \varphi(+ \infty) - \varphi(- \infty) \Big),$$
when the limits $\varphi(\pm \infty)$ exist. We next compute the values of the energy and renormalized momentum.

\begin{prop}
\label{sauber}
Assume $N = 1$ and $0 \leq c < \sqrt{2}$, and let $\v_c$ be the non-constant solution of equation \eqref{TWc} given by \eqref{schumi}. Then, the energy of $\v_c$ is equal to
\begin{equation}
\label{piquet}
E(\v_c) = \frac{(2 - c^2)^\frac{3}{2}}{3},
\end{equation}
whereas its renormalized momentum is given by
\begin{equation}
\label{fangio}
p(\v_c) = \frac{\pi}{2} - \arctan \Big( \frac{c}{\sqrt{2 - c^2}} \Big) - \frac{c}{2} \sqrt{2 - c^2}.
\end{equation}
\end{prop}

\begin{remark}
It should also be noticed that
$$|\v_c|^2 - 1 = \frac{c^2 - 2}{{2 \ch^2 \big( \frac{\sqrt{2 - c^2}}{2} x \big) }} < 0,$$
and that this quantity is integrable, so that the mass of $\v_c$ is equal to
$$m(\v_c) = - \sqrt{2 - c^2},$$
whereas its mass center is given by
$$\frac{1}{2} \int_{\R} x \big( |\v_c(x)|^2 - 1 \big) dx = 0.$$
\end{remark}

\begin{remark}
In view of identity \eqref{fangio}, it is tempting to pass to the limit $c \to 0$, and state that
\begin{equation}
\label{prost2}
p(\v_0) = \frac{\pi}{2}.
\end{equation}
As a matter of fact, this identity might be recovered using the singular lifting
$$\v_0(x) = \Big| \th \Big( \frac{x}{\sqrt{2}} \Big) \Big| \exp \Big( i \big( \sign(x) - 1 \big) \frac{\pi}{2} \Big),$$
so that the derivative of the phase $\varphi_0$ is equal to $\pi \delta_0$, and formula \eqref{prost}, in the sense of measures, yields \eqref{prost2}. 
\end{remark}

\begin{proof}
It follows from \eqref{schumi} and \eqref{barrichelo} that
$$2 |\v_c'(x)|^2 = (1 - |\v_c(x)|^2)^2 = \Big( 1 - \frac{c^2}{2} \Big)^2 \bigg( 1 - \th \Big( \frac{\sqrt{2 - c^2}}{2} x \Big)^2 \bigg)^2,$$
so that
$$E(\v_c) = \frac{(2 - c^2)^\frac{3}{2}}{4} \int_{- \infty}^{+ \infty} \Big( 1 - \th(u)^2 \Big)^2 du = \frac{(2 - c^2)^\frac{3}{2}}{3}.$$ 
Using definition \eqref{prost} and identity \eqref{sanmarino}, we next compute
$$p(\v_c) = \frac{c}{4} \int_{\R} \frac{\eta_c^2}{1 - \eta_c} = \frac{c}{4} \bigg( - \int_{\R} \eta_c + \int_{\R} \frac{\eta_c}{1 - \eta_c} \bigg),$$
where $\eta_c \equiv 1 - |\v_c|^2$. By \eqref{schumi}, we are led to
\begin{align*}
p(\v_c) = & \frac{c \sqrt{2 - c^2}}{4} \bigg( - \int_{- \infty}^{+ \infty} \big( 1 - \th(u)^2 \big) du + 2 \int_{- \infty}^{+ \infty} \frac{1 - \th(u)^2}{c^2 + (2 - c^2) \th(u)^2} du \bigg)\\
= & - \frac{c \sqrt{2 - c^2}}{2} + \frac{\sqrt{2 - c^2}}{c} \int_0^1 \frac{dv}{1 + \frac{(2 - c^2) v^2}{c^2}} dv = - \frac{c \sqrt{2 - c^2}}{2} + \arctan \bigg( \frac{\sqrt{2 - c^2}}{c} \bigg),
\end{align*}
which gives formula \eqref{fangio}.
\end{proof}

Invoking formula \eqref{fangio}, we notice that the function $c \mapsto p(\v_c)$ is smooth, decreasing, and satisfies 
\begin{equation}
\label{monaco}
\frac{\rm d}{\rm dc} \Big( p(\v_c) \Big) = - \sqrt{2 - c^2}.
\end{equation}
Hence, it performs a diffeomorphism from $(0, \sqrt{2})$ on $(0, \frac{\pi}{2})$, so that we can express $E(\v_c) \equiv \boE(p(\v_c))$ as a function of $p(\v_c)$ to obtain the following graph.

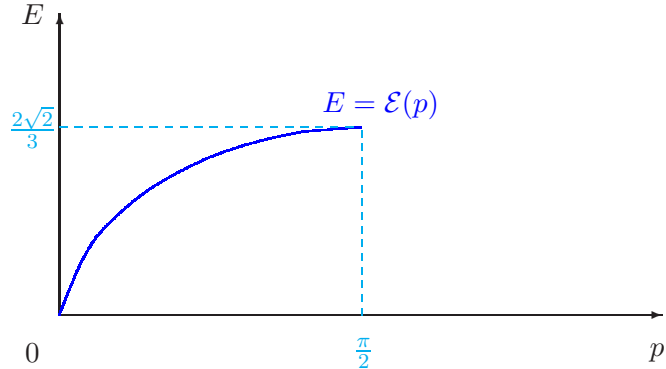
\begin{center}
\begin{picture}(90,50)(0,0)
\linethickness{0.2mm}
\put(10,10){\line(1,0){80}}
\put(90,10){\vector(1,0){0.12}}
\linethickness{0.2mm}
\put(10,10){\line(0,1){40}}
\put(10,50){\vector(0,1){0.12}}
\linethickness{0.2mm}
{\color{blue}
\qbezier(8.7,10)(8.77,10.35)(10.25,13.92)
\qbezier(10.25,13.92)(11.72,17.78)(13.7,20.5)
\qbezier(13.7,20.5)(16.8,24)(20.61,26.7)
\qbezier(20.61,26.7)(24.42,29.16)(28.7,31)
\qbezier(28.7,31)(34.21,33.2)(41.14,34.4)
\qbezier(41.14,34.4)(48.06,34.98)(48.7,35)
}
\linethickness{0.1mm}
{\color{cyan}
\multiput(46.3,10)(0,1.82){14}{\line(0,1){0.91}}
}
\linethickness{0.1mm}
{\color{cyan}
\multiput(3.6,35)(2,0){20}{\line(1,0){1}}
}
\put(0,5){\makebox(0,0)[cc]{$0$}}
\put(0,50){\makebox(0,0)[cc]{$E$}}
\put(83,5){\makebox(0,0)[cc]{$p$}}
{\color{cyan}
\put(44,5){\makebox(0,0)[cc]{$\frac{\pi}{2}$}}
\put(0,35){\makebox(0,0)[cc]{$\frac{2 \sqrt{2}}{3}$}}
}
{\color{blue}
\put(45,38){\makebox(0,0)[cc]{$E = \boE(p)$}}
}
\end{picture}
\end{center}

The curve $E = \boE(p)$ is smooth, increasing and strictly concave, and lies below the line $E = \sqrt{2} p$. Each point of the curve represents a non-constant solution $\v_c$ to \eqref{TWc} of energy $E(\v_c)$ and scalar momentum $p(\v_c)$. The speed of the solution $\v_c$ (and as a result, its position on the curve) is given by the slope of the curve. Indeed, it follows from \eqref{piquet} and \eqref{monaco} that
$$\frac{{\rm d} \boE}{{\rm d} p} \big( p(\v_c) \big) = \frac{\rm d}{{\rm d} c} \Big( E(\v_c) \Big) \bigg( \frac{\rm d}{\rm dc} \Big( p(\v_c) \Big) \bigg)^{-1} = c.$$

\subsection{The variational interpretation of equation \eqref{TWc}}
\label{carrera}

It is a rather general principle that travelling waves could be obtained minimizing the Ginzburg-Landau energy keeping the momentum fixed. In view of Lemma \ref{alonso}, the appropriate notion of momentum is given by formula \eqref{prost}. In this context, we consider the set
$$X_\p = \Big\{ v \in H^1_{\rm loc}(\R), {\rm s.t.} \ E(v) < \frac{2 \sqrt{2}}{3}, \ {\rm and } \ p(v) = \p \Big\},$$
for any $\p \geq 0$, and the minimization problem
\renewcommand{\theequation}{$\boP_\p$}
\begin{equation}
\label{pepe}
E_{\min}(\p) = \inf \Big\{ E(v), v \in X_\p \Big\}.
\end{equation}
For any $0 < \p < \frac{\pi}{2}$, let $c = c(\p)$ be the only speed $c$ such that
$$\p = \frac{\pi}{2} - \arctan \Big( \frac{c}{\sqrt{2 - c^2}} \Big) - \frac{c}{2} \sqrt{2 - c^2}.$$
We have

\begin{theorem}
\label{grouic-grouic}
Let $0 < \p < \frac{\pi}{2}$. Then, $X_\p$ is not empty, and we have
$$E(\v_{c(\p)}) = \inf \Big\{ E(v), v \in X_\p \Big\},$$
that is \eqref{pepe} is achieved, and the only minimizer, up to invariances, is $\v_{c(\p)}$.
\end{theorem}

\begin{remark}
For $\p = 0$, the solutions to the minimization problem \eqref{pepe} are the constants of modulus one, for which
$$E_{\min}(0) = 0.$$
\end{remark}

\begin{remark}
Notice that the set $X_\p$ is empty when $\p \geq \frac{\pi}{2}$. It is formally tempting to assert that 
$$E_{\min}(\p) = \frac{2\sqrt{2}}{3},$$
for any $\p \geq \frac{\pi}{2}$. Indeed, any map which vanishes at some point has an energy larger than $\frac{2\sqrt{2}}{3}$, and we can construct minimizing sequences with the appropriate limiting energy. For instance, we may consider the functions $v_n$ defined by
$$v_n(x) = \Big| \th \big( \frac{x}{\sqrt{2}} \big) \Big| \exp i \psi_n(x),$$
where
$$\psi_n(x) = \left\{ \begin{array}{lll} 0, \ {\rm if} \ x \geq \frac{1}{n},\\ q_n (n x - 1), \ {\rm if} \ - \frac{1}{n} \leq x \leq \frac{1}{n},\\ - 2 q_n, \ {\rm if} \ x \leq - \frac{1}{n}, \end{array} \right.$$
with $n \geq 1$ and $q_n = \frac{\p}{\sqrt{2} n \th \big( \frac{1}{\sqrt{2} n} \big)}$, and compute
$$E(v_n) = \frac{2 \sqrt{2}}{3} + 2 \sqrt{2} q_n^2 n^2 \bigg( \frac{1}{\sqrt{2} n} - \th \Big( \frac{1}{\sqrt{2} n} \Big) \bigg) \to \frac{2 \sqrt{2}}{3}, \ {\rm as} \ n \to + \infty,$$
whereas
$$p(v_n) = \frac{1}{2} \int_{\R} \Big( 1 - \th \Big( \frac{x}{\sqrt{2}} \Big)^2 \Big) \psi_n'(x) dx = \sqrt{2} q_n n \th \Big( \frac{1}{\sqrt{2} n} \Big) = \p.$$
\end{remark}

In order to prove Theorem \ref{grouic-grouic}, we need to analyse minimizing sequences for \eqref{pepe}. More generally, for given $0 < \p < \frac{\pi}{2}$, we consider a sequence $(u_n)_{n \in \N}$ verifying
\renewcommand{\theequation}{\arabic{equation}}
\numberwithin{equation}{section}
\setcounter{equation}{26}
\begin{equation}
\label{aminimality}
p_n \equiv p(u_n) \to \p, \ {\rm and} \ E(u_n) \to E_{\min}(\p), \ {\rm as} \ n \to + \infty.
\end{equation}
Minimizing sequences for $E_{\min}(\p)$ are a special example of sequences satisfying \eqref{aminimality}. We have

\begin{theorem}
\label{secondaire}
Let $0 < \p < \frac{\pi}{2}$ be given, and let $(u_n)_{n \in \N}$ be a sequence of maps satisfying \eqref{aminimality}. Then, there exist a subsequence $(u_{\sigma(n)})_{n \in \N}$, a sequence of points $(a_n)_{n \in \N}$, and a real number $\theta$ such that
$$u_{\sigma(n)} \big( \cdot + a_{\sigma(n)} \big) \to \exp i \theta \ \v_{c(\p)}(\cdot), \ {\rm as} \ n \to + \infty,$$
uniformly on any compact subset of $\R$. Moreover,
$$1-\vert u_{\sigma(n)} \big( \cdot + a_{\sigma(n)} \big) \vert^2 \to 1 - \vert \v_{c(\p)}(\cdot) \vert^2 \ {\rm in} \ L^2(\R), \ {\rm as} \ n \to + \infty,$$
and
$$u_{\sigma(n)}' \big( \cdot + a_{\sigma(n)} \big) \to \exp i \theta \ \v_{c(\p)}'(\cdot) \ {\rm in} \ L^2(\R), \ {\rm as} \ n \to + \infty.$$
\end{theorem}

The proof of Theorem \ref{secondaire} relies on several mostly uncorrelated observations, which we next present as separate lemmas. The first elementary observation emphasizes the role of the sonic speed $\sqrt{2}$.

\begin{lemma}
\label{colisee}
Let $\varrho$ and $\varphi$ be real-valued, smooth functions on some interval of $\R$, such that $\varrho$ is positive. Set $v = \varrho \exp i \varphi$. Then, we have the pointwise bound
$$\Big| (1 - \varrho^2) \varphi' \Big| \leq \frac{\sqrt{2}}{\varrho} e(v).$$
\end{lemma}

\begin{proof}
Notice that we have by \eqref{thenergy}, 
$$e(v) = \frac{1}{2} \Big( (\varrho')^2 + \varrho^2 (\varphi')^2 \Big) + \frac{1}{4} \Big( 1 - \varrho^2 \Big)^2 \geq \frac{1}{2} \Big( \varrho^2 | \varphi'|^2 + \frac{1}{2} (1 - \varrho^2)^2 \Big).$$
The conclusion follows from the inequality $|a b| \leq \frac{1}{2} (a^2 + b^2)$ applied to $a = \frac{1}{\sqrt{2}} (1 - \varrho^2)$ and $b = \varrho \varphi'$.
\end{proof}

As a consequence, we have

\begin{cor}
\label{coro}
Let $0 < \p < \frac{\pi}{2}$, and $v \in X_\p$. Then,
$$\underset{x \in \R}{\inf} |v(x)| \leq \frac{E(v)}{\sqrt{2} \p}.$$
In particular, if $\delta(v) \equiv 1 - \frac{E(v)}{\sqrt{2} \p} > 0$, then, given any $0 < \delta < \delta(v)$, there exists some point $x_\delta \in \R$ such that
$$1 - |v(x_\delta)| \geq \delta.$$
\end{cor}

\begin{proof}
Set $\varrho_0 = \underset{x \in \R}{\inf} |v(x)|$. It follows from Lemma \ref{colisee} that we have the pointwise bound
$$\Big| (1 - \varrho^2) \varphi' \Big| \leq \frac{\sqrt{2}}{\varrho_0} e(v),$$
and the conclusion follows by integration.
\end{proof}

The next observation shows that the size of regions where $|v|$ is away from $1$ can be bounded in terms of the energy.

\begin{lemma}
\label{changi}
Let $E > 0$ and $0 < \delta_0 < 1$ be given. There exists an integer $\ell_0 = \ell_0(E, \delta_0)$, depending only on $E$ and $\delta_0$, such that the following property holds: given any map $v \in H^1_{\rm loc}(\R)$ satisfying $E(v) \leq E$, either
$$\big| 1 - |v(x)| \big| < \delta_0, \ \forall x \in \R,$$
or there exists $\ell$ points $x_1$, $x_2$, $\ldots$, and $x_\ell$ satisfying $\ell \leq \ell_0$,
$$\big| 1 - |v(x_i)| \big| \geq \delta_0, \ \forall 1 \leq i \leq \ell,$$
and
$$\big| 1 - |v(x)| \big| \leq \delta_0, \ \forall x \in \R \setminus \underset{i = 1}{\overset{\ell}{\cup}} \big[ x_i - 1, x_i + 1 \big].$$
\end{lemma}

\begin{proof}
Set
$$\boA = \big\{ z \in \R, \ {\rm s.t.} \ \big| 1 - |v(z)| \big| \geq \delta_0 \big\},$$
and assume that $\boA$ is not empty. Considering the covering $\R = \underset{i \in \N}{\cup} I_i$, where $I_i = [i - \frac{1}{2}, i + \frac{1}{2}]$, we claim that, if $I_i \cap \boA \neq \emptyset$, then
\begin{equation}
\label{claim0}
\int_{\tilde{I_i}} e(v) \geq \mu_0,
\end{equation}
where $\tilde{I_i} = [i - 1, i + 1]$, and $\mu_0$ is some positive constant. To prove the claim, we first notice that
\begin{equation}
\label{tranquille}
|v(x) - v(y)| \leq \| v' \|_{L^2(\R)} |x - y|^\frac{1}{2} \leq \sqrt{2} E^\frac{1}{2} |x - y|^\frac{1}{2},
\end{equation}
for any $(x, y) \in \R^2$. Therefore, if $z \in \boA$, then,
$$\big| 1 - |v(y)| \big| \geq \frac{\delta_0}{2}, \ \forall y \in [z - r, z + r],$$
where $r = \frac{\delta_0^2}{8 E} $. Choosing $r_0 = \min \{ r, \frac{1}{2} \}$, we are led to
$$\int_{z - r_0}^{z + r_0} e(v) \geq\frac{1}{4} \int_{z - r_0}^{z + r_0} (1 - |v|)^2 \geq \mu_0 \equiv \frac{r_0 \delta_0^2}{8}.$$
In particular, if $z \in I_i \cap \boA$ for some $i \in \N$, then $[z - r_0, z + r_0] \subset \tilde{I_i}$, and claim \eqref{claim0} follows. To conclude the proof, we notice that
$$\underset{i \in \N}{\sum} \int_{\tilde{I_i}} e(v) = 2 E(v) \leq 2 E,$$
so that, in view of \eqref{claim0},
$$\ell \mu_0 \leq 2 E,$$
where $\ell = \Card \{ i \in \N, \ {\rm s.t.} \ I_i \cap \boA \neq \emptyset \}$. The conclusion follows setting $\ell_0 = \frac{2 E}{\mu_0} = \frac{128 E^2}{\delta_0^4}$, and choosing some point $x_i \in I_i \cap \boA$, for any $i \in \N$ such that $I_i \cap \boA \neq \emptyset$ (relabelling if necessary, the points $x_i$).
\end{proof}

We will also need the following construction.

\begin{lemma}
\label{construct}
Let $0 < |\q| \leq \frac{1}{32}$ and $0 \leq \mu \leq \frac{1}{4}$. There exists some number $\ell > 1$, and a map $w = |w| \exp i \psi \in H^1([0, \ell])$, such that
\begin{equation}
\label{quoidonc0}
w(0) = w(\ell), \ 0 < \big|1 - |w(0)| \big| \leq \mu,
\end{equation}
\begin{equation}
\label{quoidonc1}
\q = \frac{1}{2} \int_0^\ell (1 - |w|^2) \psi',
\end{equation}
and
\begin{equation}
\label{quoidonc3}
E(w) \leq 14 |\q|.
\end{equation}
\end{lemma}

\begin{proof}
Consider the functions $f_1$ and $\psi_1$ defined on the interval $[0, 2]$ by
$$f_1(s) = s \ {\rm on } \ \Big[ 0, \frac{1}{2} \Big], \ f_1(s) = 1 - s \ {\rm on } \ \Big[ \frac{1}{2}, 1 \Big], \ {\rm and} \ f_1(s) = 0 \ {\rm on} \ \Big[ 1, 2 \Big],$$
and
$$\psi_1(s) = s \ {\rm on} \ [0, 1], \ {\rm and} \ \psi_1(s) = 2 - s \ {\rm on} \ [0, 1].$$
For a given positive number $\lambda > 0$, we consider the functions defined on $[0, 2 \lambda]$ by
$$f_\lambda(s) = \frac{1}{\lambda} f \Big( \frac{s}{\lambda} \Big), \ {\rm and} \ \psi_\lambda(s) = \psi \Big( \frac{s}{\lambda} \Big),$$
so that $|f_\lambda| \leq \frac{1}{2 \lambda}$, $|\psi_\lambda'| = \frac{1}{\lambda}$, $f_\lambda(0) = f_\lambda(2 \lambda) = 0$, 
 $\psi_\lambda(0) = \psi_\lambda(2 \lambda) = 0$, and
\begin{equation}
\label{vroum-vroum}
\int_0^{2 \lambda} f_\lambda \psi_\lambda' = \frac{1}{4 \lambda}, \ \int_0^{2 \lambda} f_\lambda = \frac{1}{4}, \ \int_0^{2 \lambda} f_\lambda^2 = \frac{1}{12 \lambda}, \ \int_0^{2 \lambda} (f_\lambda')^2 = \frac{1}{\lambda^3}, \ {\rm and} \ \int_0^{2 \lambda} (\psi_\lambda')^2 = \frac{2}{\lambda}.
\end{equation}
We then choose $\lambda = \frac{1}{8 |\q|}$, so that, $\frac{1}{\lambda} \leq \frac{1}{4}$, introduce a new parameter $\delta > 0$ to be determined later, and consider the function
$$\rho_{\lambda, \delta} = \sqrt{1 - \delta - f_\lambda},$$
so that $1 - \rho_{\lambda, \delta}^2 = f_\lambda + \delta$. It follows from our choice of parameter $\lambda$ that
\begin{equation}
\label{revroum}
|\q| = \frac{1}{2} \int_0^{2 \lambda} f_\lambda \psi_\lambda' = \frac{1}{2} \int_0^{2 \lambda} (f_\lambda + \delta) \psi_\lambda' =
\frac{1}{2} \int_0^{2 \lambda} (1 - \rho_{\lambda, \delta}^2) \psi_\lambda'.
\end{equation}
We finally choose $\ell = 2\lambda$ and
$$w = \left\{ \begin{array}{ll} \rho_{\lambda, \delta} \exp i \psi_\lambda, \ {\rm if} \ \q > 0,\\ \rho_{\lambda, \delta} \exp (- i \psi_\lambda), \ {\rm if} \ \q < 0. \end{array} \right.$$
Condition \eqref{quoidonc1} is fullfilled with this choice of $w$ in view of \eqref{revroum}. Moreover, by construction, $w(0) = w(\ell) = \rho_{\lambda, \delta} = \sqrt{1 - \delta}$, so that conditions \eqref{quoidonc0} are satisfied for any $\delta \leq \mu^2$. We finally compute
$$E(w) = \int_0^{2 \lambda} \bigg( \frac{(f_\lambda')^2}{8 (1 - \delta - f_\lambda)} + \Big( 1 - \delta - f_\lambda \Big) \frac{(\psi_\lambda')^2}{2} + \frac{f_\lambda^2}{4} + \frac{\delta f_\lambda}{2} + \frac{\delta^2}{4} \bigg),$$
so that, since
$$0 \leq f_\lambda + \delta \leq \frac{1}{2 \lambda} + \delta \leq \mu^2 + \frac{1}{8} \leq \frac{1}{2},$$
it follows from \eqref{vroum-vroum} that
$$E(w) \leq \int_0^{2 \lambda} \bigg( \frac{(f_\lambda')^2}{4} + \frac{(\psi_\lambda')^2}{2} + \frac{f_\lambda^2}{4} + \frac{\delta f_\lambda}{2} + \frac{\delta^2}{4} \bigg) \leq \frac{1}{4 \lambda^3} + \frac{1}{\lambda} + \frac{1}{48 \lambda} + \frac{\delta}{8} + \frac{\delta^2 \lambda}{2}.$$
Inequality \eqref{quoidonc3} follows choosing $\delta = \min \{ \mu^2, \frac{1}{\lambda} \}$. 
\end{proof}

We are now in position to undertake the study of sequences satisfying \eqref{aminimality}. We set
$$\delta_n = \delta(u_n) = 1 - \frac{E(u_n)}{\sqrt{2} p(u_n)},$$
so that, by \eqref{aminimality}, $\delta_n \to \delta_\p \equiv 1 - \frac{E_{\min}(\p)}{\sqrt{2} \p}$, as $n \to + \infty$. Moreover, it follows from the properties of the curve $p \mapsto \boE(\p)$ (see Subsection \ref{energetic}) that
\begin{equation}
\label{deltan}
\delta_\p > 0,
\end{equation}
for any $0 < \p < \frac{\pi}{2}$, so that
\begin{equation}
\label{deltaplane}
\delta_n > \frac{\delta_\p}{2},
\end{equation}
for $n$ sufficiently large. We first have

\begin{lemma}
\label{onlyweak}
Let $0 < \p < \frac{\pi}{2}$, and let $(u_n)_{n \in \N}$ be a sequence of maps satisfying \eqref{aminimality}. Then, there exists a subsequence $(u_{\sigma(n)})_{n \in \N}$, some numbers $\theta$ and $\tilde{x}$, and a solution $\v_c$ to \eqref{TWc} such that
$$u_{\sigma(n)} \rightharpoonup \exp i \theta \ \v_c( \cdot + \tilde{x}) \ {\rm weakly} \ {\rm in} \ H^1([- A, A]), \ {\rm as} \ n \to + \infty,$$
for any $A > 0$. In particular, we have
\begin{equation}
\label{semis}
\int_{- A}^A e(\v_c) \leq \underset {n \to + \infty}{\liminf} \bigg( \int_{- A}^A e(u_n) \bigg), \ {\rm and} \ \int_{- A}^A (1 - \varrho_c^2) \varphi_c' = \underset{n \to + \infty}{\lim} \bigg( \int_{- A}^A (1 - \varrho_n^2) \varphi_n' \bigg),
\end{equation}
where we have written $u_n = \varrho_n \exp i \varphi_n$, and $\v_c = \varrho_c \exp i \varphi_c$.
\end{lemma}

\begin{proof}
Since $(E(u_n))_{n \in \N}$ is bounded by assumption \eqref{aminimality}, it follows from standard compactness results that there exists a subsequence $(u_{\sigma(n)})_{n \in \N}$, and a map $u \in H^1_{\rm loc}(\R)$ such that
$$u_{\sigma(n)} \rightharpoonup u \ {\rm in} \ H^1([- A, A]), \ {\rm as} \ n \to + \infty,$$
for any $A > 0$. Assertions \eqref{semis} follow using Rellich's compactness theorem, and the lower semi-continuity of $E$. Hence, it remains to prove that the limiting map $u$ solves \eqref {TWc} on $(- A, A)$. For that purpose, we consider a smooth map $\xi$, with compact support in $(-A, A)$, such that
\begin{equation}
\label{ortho}
\int_{\R} \langle i u, \xi' \rangle = 0.
\end{equation}
We claim that, for any $t$ sufficiently small,
\begin{equation}
\label{claim}
\int_{- A}^A e(u_{\sigma(n)} + t \xi) \geq \int_{- A}^A e(u_{\sigma(n)}) + O(t^2) + \underset{n \to \infty}{o(1)}.
\end{equation}
To establish the claim, we first observe that, in view of \eqref{aminimality}, we may assume without loss of generality that $E(u_n) < \frac{2 \sqrt{2}}{3}$, so that $u_n$ has no zero. Expandind $p(u_n + t \xi)$, we obtain
\begin{align*}
p(u_n + t\xi) & = p(u_n) + t \int_{\R} \langle i u_n, \xi' \rangle + O(t^2)\\
& = p_n + O(t^2) = \p + O(t^2) + \underset{n \to \infty}{o(1)},
\end{align*}
so that, setting $\q_{n, t} = \p - \p(u_n + t \xi)$, we are led to
$$\q_n = O(t^2) + \underset{n \to \infty}{o(1)}.$$
We next construct a comparison map $v_{n, t}$ for $E_{\min}(\p)$ applying several modifications to the map $u_n + t \xi$. For that purpose, we invoke Lemma \ref{construct} with $\q = \q_{n, t}$, and $\mu = \mu_{n, t} = \inf \{ \frac{1}{4}, \nu_{n, t} \}$, where $\nu_{n, t} = \sup \{ |1 - |u_n(x)||, x \not \in [-A, A] \}$. This yields a positive number $\ell_{n, t} > 1$, and a map $w_{n, t} = |w_n| \exp i \psi_n$, defined on $[0, \ell_n(t)]$ such that
$$w_n(0) = w_n(\ell_{n, t}), \ {\rm and} \ \big| 1 - |w_{n, t}(0)| \big| \leq \mu_{n, t},$$
and such that
$$\q_{n, t} = \frac{1}{2} \int_0^{\ell_{n, t}} (1 - |w_{n, t}|)^2 \psi_n',$$
and
\begin{equation}
\label{diplodocus}
E(w_{n, t}) \leq 14 |\q_{n ,t}| = O(t^2) + \underset{n \to \infty}{o(1)}.
\end{equation}
In view of the mean value theorem, there exists some point $x_n$ in $[A, + \infty)$ such that
$|u_n(x_n)| = |w_{n, t}(0)|$. Multiplying possibly $w_{n, t}$ by a constant complex number of modulus one, we may therefore assume, without loss of generality, that
$u_n(x_n) = w_n(0)$. We define the comparison map $v_{n, t}$ as follows
\begin{equation}
\label{defv}
\begin{split}
v_{n, t}(x) & = u_n(x) + t \xi(x), \ \forall x < x_n,\\
v_{n, t}(x) & = w_n(x - x_n), \ \forall x_n \leq x \leq x_n + \ell_{n, t},\\
v_{n, t}(x) & = u_n(x - \ell_{n, t}) + t \xi(x - \ell_{n,t}), \ \forall x \geq x_n + \ell_{n, t}.
\end{split}
\end{equation}
We verify that
\begin{equation}
\label{oggi}
E(v_{n, t}) = E(u_n + t \xi) + E(w_{n, t}), \ {\rm and} \ p(v_{n, t}) = p(u_n + t \xi) + \q_{n, t} = \p,
\end{equation}
so that $v_{n, t}$ is a comparison map for $E_{\min}$, and therefore
\begin{equation}
\label{compa}
E(v_{n, t}) \geq E_{\min}(\p).
\end{equation}
On the other hand, we have in view of assumption \eqref{aminimality},
\begin{equation}
\label{brachiosaurus}
E(u_n) = E_{\min}(\p) + o(1),
\end{equation}
whereas, since $\xi$ has compact support in $(-A, A)$,
\begin{equation}
\label{albertosaurus}
E(u_n + t \xi) - E(u_n) = \int_{- A}^A \Big( e(u_n + t \xi) - e(u_n)\Big).
\end{equation}
Combining \eqref{albertosaurus} with \eqref{brachiosaurus}, \eqref{compa} and \eqref{diplodocus}, we establish claim \eqref{claim}.

To complete the proof of Lemma \ref{onlyweak}, we expand the integral in \eqref{claim} so that
$$t \int_{-A}^A \Big( u_n' \xi' + \xi u_n (1 - |u_n|^2) \Big) \geq O(t^2) + \underset{n \to \infty}{o(1)}.$$
We then let $n$ tend to $+ \infty$. This yields, in view of the compact embedding of $H^1([-A, A])$ in $C^0([-A, A])$,
$$t \int_{- A}^A \Big( u' \xi'+ \xi u(1 - |u|^2) \Big) \geq O(t^2).$$
Letting $t$ tend to $0^+$ and $0^-$, we deduce
$$\int_{- A}^A u' \xi' + \xi u (1 - |u|^2) = 0.$$
Since $\xi$ is any arbitrary function with compact support verifying \eqref{ortho}, this shows that there exists some constant $c$ such that $u$ solves \eqref{TWc}, and hence is of the form $u = \exp i \theta \ \v_c(\cdot + \tilde{x})$, for some $c > 0$.
\end{proof}

\begin{remark}
At this stage, we might only consider maps $\v_c$ with positive speed $c$, using the convention that, for any $c < 0$, $\v_c(x) = \v_{- c}(- x)$.
\end{remark}

It might happen that the limit map provided by Lemma \ref{onlyweak} is a constant one. To capture the possible losses at infinity, we need to implement a concentration-compactness argument. The first step in this argument is to invoke Lemma \ref{changi}, and Corollary \ref{coro} to assert

\begin{prop}
\label{carbonifere}
Let $0 < \p < \frac{\pi}{2}$ be given, and let $(u_n)_{n \in \N}$ be a sequence of maps satisfying \eqref{aminimality}. There exists an integer $\ell_\p$, depending only on $\p$, such that there exists $\ell_n$ points $x_1^n$, $x_2^n$, $\ldots$, and $x_{\ell_n}^n$ satisfying $\ell_n \leq \ell_\p$, and 
\begin{equation}
\label{megalodon}
\big| 1 - |u_n(x_i^n)| \big| \geq \frac{\delta_\p}{4}, \ \forall 1 \leq i \leq \ell_n,
\end{equation}
and
\begin{equation}
\label{megalodon1}
\big| 1 - |u_n(x)| \big| \leq \frac{\delta_\p}{4}, \ \forall x \in \R \setminus \underset{i = 1}{\overset{\ell_n}{\cup}} \big[ x_i^n - 1, x_i ^n + 1 \big].
\end{equation}
\end{prop}

Passing possibly to a subsequence, we may assume that the number $\ell_n$ does not depend on $n$, and set $\ell = \ell_n$. A standard compactness argument shows, that passing again possibly to a further subsequence, and relabelling possibly the points $x_i^n$, we may find some integer
$1 \leq \tilde{\ell} \leq \ell$, and some number $R > 0$ such that
\begin{equation}
\label{trias}
|x_i^n - x_j^n| \to + \infty, \ {\rm as} \ n \to + \infty, \ \forall 1 \leq i \neq j \leq \tilde{\ell},
\end{equation}
and
\begin{equation}
\label{trias1}
x_i^n \in \underset{i = 1}{\overset{\tilde{\ell}}{\cup}} B(x_i^n, R), \ \forall \tilde{\ell} < i \leq \ell.
\end{equation}
Going back to Proposition \ref{carbonifere}, we deduce
\begin{equation}
\label{megalodon2}
\big| 1 - |u_n(x)| \big| \leq \frac{\delta_\p}{4}, \ \forall x \in \R \setminus \underset{i = 1}{\overset{\tilde{\ell}}{\cup}} B(x_i^n, R+1),
\end{equation}
so that, invoking Lemma \ref{colisee}, we have on $\R \setminus \underset{i = 1}{\overset{\tilde{\ell}}{\cup}} B(x_i^n, R + 1)$,
\begin{equation}
\label{megalodon3}
\frac{1}{2} \Big| (|u_n| ^2 - 1) \varphi_n' \Big| \leq \frac{e(u_n)}{\sqrt{2} \big( 1 - \frac{\delta_\p}{4} \big)}.
\end{equation}
We are now in position to provide the proof to Theorem \ref{secondaire}.

\begin{proof}[Proof of Theorem \ref{secondaire}]
We divide the proof into several steps.

\begin{step}
\label{stepun}
For any $1 \leq i \leq \tilde{\ell}$, there exists numbers $c_i \neq 0$, $\tilde{x_i}$ and $\theta_i$ such that
$$u_n(\cdot + x_i^n) \rightharpoonup \exp(i \theta_i) \v_{c_i} (\cdot + \tilde{x_i}).$$
\end{step}

Applying Lemma \ref{onlyweak} to the sequence $u_n(\cdot + x_i^n)_{n \in \N}$ yields the existence of the limiting solution to \eqref{TWc}, $\exp(i \theta_i) \v_{c_i}(\cdot + \tilde{x_i})$. It remains to prove that the function $\v_{c_i}$ is not a constant function. This is a consequence of the fact that
$$|u_n(x_i^n)| \leq 1 - \frac{\delta_\p}{4},$$
so that, since by compact embedding, we have uniform convergence on compact sets, we obtain
$$|\v_{c_i}(\tilde{x_i})| \leq 1 - \frac{\delta_\p}{4}.$$
Hence, $\v_{c_i}$ is not a constant map.

\begin{step}
\label{stepdeux}
Given any number $\mu > 0$, there exists a number $A_\mu > 0$, and $n_\mu \in \N$, such that, if $n \geq n_0$, then
$$\int_{\underset{i = 1}{\overset{\tilde{\ell}}{\cup}} B(x_i^n, A_\mu)} e(u_n) \geq \underset{i = 1}{\overset{\tilde{\ell}}{\sum}} E(\v_{c_i}) - \mu,$$
and
$$\bigg| \frac{1}{2} \int_{\underset{i = 1}{\overset{\tilde{\ell}}{\cup}} B(x_i^n, A_\mu)} (\varrho_n^2 - 1) \varphi_n' - \underset{i = 1}{\overset{\tilde{\ell}}{\sum}} \p_i \bigg| \leq \mu,$$
where $\p_i = p(\v_{c_i})$.
\end{step}

To prove Step \ref{stepdeux}, we choose $A > R + 1$ so that, for any $1 \leq i \leq \tilde{\ell}$, we have
$$\int_{- A}^A e(v_{c_i}) \geq E(v_{c_i}) - \frac{\mu}{2 \tilde{\ell}},$$
and
$$\frac{1}{2} \bigg| \int_{- A}^A (|v_{c_i}|^2 - 1) \varphi_c' - \p_i \bigg| \leq \frac{\mu}{2 \tilde{\ell}}.$$
The conclusion follows from the convergences stated in \eqref{semis}.

\begin{step}
\label{steptrois}
We have
$$\bigg| \frac{1}{2} \int_{\R \setminus \underset{i = 1}{\overset{\tilde{\ell}}{\cup}} B(x_i^n, A)} (\varrho_n^2 - 1) \varphi_n' \bigg| \leq \frac{1}{\sqrt{2} \big( 1 - \frac{\delta_\p}{4} \big)} \int_{\R \setminus \underset{i = 1}{\overset{\tilde{\ell}}{\cup}} B(x_i^n, A)} e(u_n).$$
\end{step}

To establish this inequality, it is sufficient to integrate \eqref{megalodon3}.

Passing possibly to a further subsequence, we may assume that there exist some numbers $\p_\mu$ and $E_\mu$ such that
$$\frac{1}{2} \int_{\R \setminus \underset{i = 1}{\overset{\tilde{\ell}}{\cup}} B(x_i^n, A)} (\varrho_n^2 - 1) \varphi_n' \to \p_{\mu}, \ {\rm and} \ \int_{\R \setminus \underset{i = 1}{\overset{\tilde{\ell}}{\cup}} B(x_i^n, A)} e(u_n) \to E_\mu, \ {\rm as} \ n \to + \infty.$$
Going back to Step \ref{stepdeux}, and letting $n \to + \infty$, we are led to the estimates
$$\bigg| \p - \underset{i = 1}{\overset{\tilde{\ell}}{\sum}} \p_i - \p_\mu \bigg| \leq \mu, \ {\rm and } \ \boE(\p) \geq \underset{i = 1}{\overset{\tilde{\ell}}{\sum}} \boE(\p_i) + E_\mu - \mu,$$
whereas Step \ref{steptrois} yields
$$\sqrt{2} \Big( 1 - \frac{\delta_\p}{4} \Big) \p_\mu \leq E_\mu.$$
Letting $\mu \to 0$, we may assume that for some subsequence $(\mu_m)_{m \in \N}$ tending to $0$, we have
$$\p_{\mu_m} \to \tilde{\p}, \ {\rm and} \ E_{\mu_m} \to \tilde{E}, \ {\rm as} \ {m \to + \infty}.$$
Our previous inequalities then yield
\begin{equation}
\label{ouida}
\begin{split}
\p & = \underset{i = 1}{\overset{\tilde{\ell}}{\sum}} \p_i + \tilde{p},\\
\boE(\p) & \geq \underset{i = 1}{\overset{\tilde{\ell}}{\sum}} \boE(\p_i) + \tilde{E},
\end{split}
\end{equation}
with
\begin{equation}
\label{nenni}
\sqrt{2} \Big( 1-\frac{\delta_\p}{4} \Big) \tilde{\p} \leq \tilde{E}.
\end{equation}

\begin{step}
\label{stepquatre}
We have $\tilde{E} = \tilde{\p} = 0$, and $\tilde{\ell} = 1$.
\end{step}

This statement is a consequence of properties of the curve $\p \mapsto \boE(\p)$, and the definition of $\delta_\p$. Assume first that $\tilde{p} \leq 0$. Then, it follows from \eqref{ouida} that
\begin{equation}
\label{argentinosaurus}
\p\geq \underset{i = 1}{\overset{\tilde{\ell}}{\sum}} |\p_i|,
\end{equation}
and
\begin{equation}
\label{argentinosaurus2}
\boE(\p) \geq \underset{i = 1}{\overset{\tilde{\ell}}{\sum}} \boE(|\p_i|).
\end{equation}
On the other hand, as observed before, the function $\boE$ is concave on $\R^+$, non-decreasing, and $\boE(0) = 0$, so that $\boE$ is subadditive, and hence in view of \eqref{argentinosaurus},
$$\boE(\p) \leq \underset{i = 1}{\overset{\tilde{\ell}}{\sum}} \boE(|\p_i|),$$
so that \eqref{argentinosaurus2} becomes an equality. However, since $\boE$ is strictly concave, this is possible if and only if $\ell = 1$, and \eqref{argentinosaurus} is an equality. Step \ref{stepquatre} therefore follows in the case considered.

It remains to handle the case $\tilde{p} \geq 0$. In view of the definition on $\delta_\p$, we have
\begin{equation}
\label{brachiosaure}
\frac{\boE(\p)}{\p} = \sqrt{2} \big( 1 - \delta_\p \big) < \sqrt{2} \Big( 1 - \frac{\delta_\p}{4} \Big),
\end{equation}
so that by concavity, we have for any $0 < \s < \p$,
\begin{equation}
\label{corde}
\boE(\s) > \s \frac{\boE(\p)}{\p} = \s \sqrt{2} \big( 1 - \delta_\p \big).
\end{equation}
Setting $\check{\p} = \p - \tilde{\p} =\underset{i = 1}{\overset{\tilde{\ell}}{\sum}} \boE(\p_i)$, we obtain using \eqref{ouida}, \eqref{nenni}, \eqref{brachiosaure} and \eqref{corde},
\begin{equation}
\label{therope}
\begin{split}
\boE(\check{\p}) & \geq \check{\p} \frac{\boE(\p)}{\p} = \boE(\p) - \tilde{p} \frac{\boE(\p)}{\p} = \boE(\p) - \sqrt{2} \tilde{p} \Big( 1 - \delta_\p \Big)\\
& > \boE(\p) - \sqrt{2} \tilde{p} \big( 1 - \frac{\delta_\p}{4} \big) \geq \boE(\p) - \tilde{E}\\
& > \underset{i = 1}{\overset{\tilde{\ell}}{\sum}} \boE(|\p_i|).
\end{split}
\end{equation}
We then conclude following the first case, replacing $\p$ by $\check{\p}$, and using the subadditivity of the function $\boE$.

\begin{step}
\label{stepcinq}
Proof of Theorem \ref{carbonifere} completed.
\end{step}

At this stage, we have shown that there exist some $\theta = \theta_1$, and some solution $\v_c = \v_{c_1}$ such that
$$u_n(\cdot + a_n) \rightharpoonup \exp i \theta \ \v_{c_1},$$
weakly in $H^1(\R)$, as $n \to + \infty$, where $a_n = x_1^n - \tilde{x_1}$, and that we have the convergence of the energy $E(u_n) \to E(\v_{c_1})$. In particular,
$$u'_n(\cdot + a_n) \to \exp i \theta \ \v_{c_1}',$$
strongly in $L^2(\R)$ as $n \to + \infty$, and the conclusion follows by compact embedding.
\end{proof}

\subsection{Orbital stability}
\label{Orbit1}

We recall first the classical notion of orbital stability (see, for instance, \cite{Benjami1}). A solution $\v_c$ to \eqref{TWc} is said to be orbitally stable in a metric space $X(\R)$, if and only if given any $\varepsilon > 0$, there exists some $\delta > 0$ such that for any solution $\Psi$ to \eqref{GP} in $X(\R)$, if
$$d_X \big( \Psi(\cdot , 0), \v_c \big) \leq \delta,$$
then
$$\sup_{t\in \R} \bigg( \inf_{(a, \theta) \in \R^2} d_X \big( \Psi(\cdot, t), \exp i \theta \ \v_c(\cdot - a) \big) \bigg) \leq \varepsilon,$$
As a preliminary step, this definition requires of course to prove first that the Cauchy problem for \eqref{GP} is globally well-posed in $X(\R)$ (see \cite{Gerard2}). A natural choice is the energy space
\begin{equation}
\label{equation}
\boX^1(\R) = \{ w \in L^\infty(\R), \ {\rm s.t.} \ w' \in L^2(\R), 1 - |w|^2 \in L^2(\R) \}.
\end{equation}
Given any $v_0 \in \boX^1(\R)$, Zhidkov \cite{Zhidkov1} (see also \cite{Gerard2}) established that \eqref{GP} has a global solution $\Psi$ with initial datum $v_0$. Moreover, the Ginzburg-Landau energy is conserved,
\begin{equation}
\label{helios}
E(\Psi(\cdot,t)) = E(v_0), \ \forall t \in \R.
\end{equation}
Notice that $E$ is a continuous function on $\boX^1(\R)$. If we assume moreover that
$$E(v_0) < \frac{2 \sqrt{2}}{3},$$
using \eqref{helios}, we are led to
$$E(\Psi(\cdot,t)) < \frac{2 \sqrt{2}}{3}, \ \forall t \in \R,$$
so that, by Remark \ref{botox}, we may write $\Psi(\cdot, t) = \varrho(\cdot, t) \exp i \varphi(\cdot, t)$ for any $t \in \R$. We may therefore define the scalar momentum $p(\Psi(\cdot, t))$ using \eqref{prost}. It is shown in \cite{LinZhiw1} that the scalar momentum $p$ is then a conserved quantity, i.e.
\begin{equation}
\label{eole}
p(\Psi(\cdot, t)) = p(v_0), \ \forall t \in \R.
\end{equation}

Given any $A > 0$, we consider on $\boX^1(\R)$ the distance $ d_{A, \boX^1}$ defined by 
\begin{equation}
\label{distxi}
d_{A, \boX^1}(v_1, v_2) \equiv \| v_1 - v_2 \|_{L^\infty([-A, A])} + \| v_1' - v_2' \|_{L^2(\R)} + \| |v_1| - |\v_2| \|_{L^2(\R)}.
\end{equation}
Following ideas from Grillakis-Shatah-Strauss \cite{GriShSt1}, Zhiwu Lin \cite{LinZhiw1} proved the orbital stability of the solutions $\v_c$ to \eqref{TWc}, for any $0 < c < \sqrt{2}$, when the perturbations $w$ are taken in the space $\boX^1(\R)$.

Using our study of the minimization problem \eqref{pepe}, we would like to give a short proof of the result of Zhiwu Lin \cite{LinZhiw1} which follows the classical compactness argument of \cite{CazeLio1}. We first recall

\begin{theorem}[\cite{LinZhiw1}]
\label{voyager}
For $v_0 \in \boX^1(\R)$, consider the global in time solution $\Psi$ having initial datum $v_0$. Let $0 < c < \sqrt{2}$ be given. For any numbers $\varepsilon > 0$ and $A > 0$, there exists some positive number $\delta$, such that, if 
\begin{equation}
\label{pioneer}
d_{A, \boX^1}(v_0, \v_c) \leq \delta,
\end{equation}
then, we have 
\begin{equation}
\label{solaris}
\sup_{t \in \R} \bigg( \inf_{(a, \theta) \in \R^2} d_{A, \boX^1} \Psi( \cdot, t) , \exp i \theta \ \v_c(\cdot - a)) \bigg) < \varepsilon.
\end{equation}
\end{theorem}

\begin{proof}
We assume by contradiction that $\v_c$ is not orbitally stable for the distance $d_{A, \boX^1}$. In this case, we may find a positive number $\varepsilon_0$, sequences of numbers $(\delta_n)_{ n \in \N}$ and $(t_n)_{n \in \N}$, and a sequence of functions $(v_0^n)_{n \in \N}$ such that $\delta_n \to 0$, as $n \to + \infty$,
\begin{equation}
\label{pioneer1}
d_{A, \boX^1} (v_0^n , \v_c) < \delta_n,
\end{equation}
and
\begin{equation}
\label{solaris1}
\inf_{(a, \theta) \in \R^2} d_{A, \boX^1} (\Psi^n( \cdot, t_n) , \exp i \theta \ \v_c(\cdot - a_n)) \geq \varepsilon_0,
\end{equation}
where $\Psi^n$ denotes the solution to \eqref{GP} with initial datum $v_0^n$. Denoting $w_n \equiv \Psi^n( \cdot, t_n)$, the conservation of the Ginzburg-Landau energy and of the scalar momentum implies that
\begin{equation}
\label{mars}
E(w_n) = E(v_0^n), \ {\rm and} \ p(w_n) = p(v_0^n),
\end{equation}
for any $n \in \N$. Moreover, the energy is continuous on $\boX^1(\R)$, so that, invoking \eqref{pioneer1},
$$E(v_0^n) \to E(\v_c) = E_{\min}(\p), \ {\rm as} \ n \to + \infty,$$
where $\p = p(\v_c)$ is given by formula \eqref{fangio}. Therefore, by \eqref{mars},
\begin{equation}
\label{voyager1}
E(w_n) \to E_{\min}(\p), \ {\rm as} \ n \to + \infty.
\end{equation}
In view of \eqref{prost}, the scalar momentum may be written as
\begin{equation}
\label{hubble}
p(v) = \frac{1}{2} \int_{\R} \frac{|v|^2 - 1}{|v|^2} \langle i v, v' \rangle,
\end{equation}
for any function $v \in \boX^1(\R)$. By \eqref{pioneer1}, we have
\begin{equation}
\label{explorer}
\langle i v_0^n, {v_0^n}' \rangle \to \langle i \v_c, \v_c' \rangle \ {\rm in} \ L^2(\R), \ {\rm as} \ n \to + \infty.
\end{equation}
On the other hand, in view of \eqref{schumi}, the modulus of $\v_c$ has a minimum value equal to $\frac{c}{\sqrt{2}}$ on $\R$ which is achieved at the origin. Since $(|v_0^n|)_{n \in \N}$ uniformly converges to $|\v_c|$ on $\R$ as $n \to + \infty$ by \eqref{pioneer1}, we may assume that
$$\underset{x \in \R}{\inf} |v_0^n(x)| \geq \frac{c}{2 \sqrt{2}},$$
for $n$ sufficiently large, so that
$$\bigg| \frac{|v_0^n|^2 - 1}{|v_0^n|^2} - \frac{|\v_c|^2 - 1}{|\v_c^2|} \bigg| \leq \frac{4}{c^4} \Big| |v_0^n|^2 - |\v_c|^2 \Big|.$$
By \eqref{pioneer1}, it follows that
$$\frac{|v_0^n|^2 - 1}{|v_0^n|^2} \to \frac{|\v_c|^2 - 1}{|\v_c^2|} \ {\rm in} \ L^2(\R), \ {\rm as} \ n \to + \infty,$$
so that by \eqref{hubble}, and \eqref{explorer},
$$p(v_0^n) \to p(\v_c) = \p, \ {\rm as} \ n \to + \infty.$$
By \eqref{mars}, we are led to
\begin{equation}
\label{voyager2}
p(w_n) \to \p, \ {\rm as} \ n \to + \infty,
\end{equation}
so that, by \eqref{voyager1}, the sequence $(w_n)_{n \in \N}$ verifies assumptions \eqref{aminimality}. Hence, by Theorem \ref{secondaire}, there exist some points $(a_n)_{n \in \N}$, and some real number $\theta$ such that, up to some subsequence,
\begin{equation}
\label{leverrier}
w_n \big( \cdot + a_n \big) \to \exp i \theta \ \v_{c(\p)}(\cdot), \ {\rm as} \ n \to + \infty,
\end{equation}
uniformly on any compact set of $\R$, and
\begin{equation}
\label{neptune}
w_n' \big( \cdot + a_n \big) \to \exp i \theta \ \v_{c(\p)}'(\cdot) \ {\rm in} \ L^2(\R), \ {\rm as} \ n \to + \infty.
\end{equation}
In view of \eqref{voyager1} and \eqref{neptune}, we have
$$\int_{\R} \big( 1 - |w_n|^2 \big)^2 \to \int_{\R} \big( 1 - |\v_c|^2 \big)^2, \ {\rm as} \ n \to + \infty,$$
whereas, by \eqref{voyager1} and \eqref{leverrier}, up to some subsequence, we have
$$1 - \big| w_n \big( \cdot + a_n \big) \big|^2 \rightharpoonup 1 - \big| \v_{c(\p)}(\cdot) \big|^2 \ {\rm in} \ L^2(\R), \ {\rm as} \ n \to + \infty.$$
This yields
$$1 - \big| w_n \big( \cdot + a_n \big) \big|^2 \to 1 - \big| \v_{c(\p)}(\cdot) \big|^2 \ {\rm in} \ L^2(\R), \ {\rm as} \ n \to + \infty.$$
Hence, by \eqref{leverrier} and \eqref{neptune},
$$w_n \big( \cdot + a_n \big) \to \exp i \theta \ \v_{c(\p)}(\cdot) \ {\rm in} \ \boX^1(\R), \ {\rm as} \ n \to + \infty,$$
which gives a contradiction with \eqref{solaris1}, and completes the proof of Theorem \ref{voyager}.
\end{proof}

\begin{remark}
Di Menza and Gallo \cite{DiMeGal1} proved the linear stability of $\v_0$ submitted to small perturbations in $\{ \v_0 \} + H^1(\R)$ (see also \cite{Gerard2}).
\end{remark}

\subsection{Relating \eqref{TWc} to the Korteweg-de Vries equation}

Travelling wave solutions to \eqref{TWc} are related to the soliton of the Korteweg-de Vries equation as follows.
Set $\varepsilon = \sqrt{2 - c^2}$, and consider the scaled function
\begin{equation}
\label{senna1}
N_\varepsilon(x) = \frac{1}{\varepsilon^2} \eta_c \Big( \frac{x}{\varepsilon} \Big),
\end{equation}
where $\eta_c \equiv 1 - |\v_c|^2$. Invoking \eqref{maranello}, we are led to
$$N_\varepsilon(x) = N(x) \equiv \frac{1}{2 \ch^2 \big( \frac{x}{2} \big)}.$$
A remarkable property of $N$ is that it represents the classical soliton to the Korteweg-de-Vries equation
\renewcommand{\theequation}{KdV}
\begin{equation}
\label{KdV}
\partial_t w + \partial_x^3 w + 6 w \partial_x w = 0.
\end{equation}
Concerning the phase $\varphi_c$ of $\v_c$, we consider the scale change
\renewcommand{\theequation}{\arabic{equation}}
\numberwithin{equation}{section}
\setcounter{equation}{74}
\begin{equation}
\label{senna2}
\Theta_\varepsilon(x) = \frac{\sqrt{2}}{\varepsilon} \varphi_c \Big( \frac{x}{\varepsilon} \Big),
\end{equation}
so that we similarly obtain from \eqref{sanmarino},
$$ \Theta_\varepsilon(x)' = \sqrt {1 - \frac{\varepsilon^2}{2}} \frac{N(x)}{1 - \varepsilon^2 N(x)} \longrightarrow N(x), \ {\rm as} \ \varepsilon \to 0.$$
It is also of interest to compute the corresponding energies. Using scale changes \eqref{senna1} and \eqref{senna2}, we are led to
$$E(\v_c) = \frac{\varepsilon^5}{8} \int_\R \frac{(N')^2}{1 - \varepsilon^2 N} + \frac{\varepsilon^3}{4} \Big( 1 - \frac{\varepsilon^2}{2} \Big) \int_\R \frac{N^2}{1 - \varepsilon^2 N}+ \frac{\varepsilon^3}{4} \int_\R N^2,$$
and
$$p(\v_c) = \frac{\varepsilon^3}{2 \sqrt{2}} \sqrt{1 - \frac{\varepsilon^2}{2}} \int_\R \frac{N^2}{1 - \varepsilon^2 N},$$
so that the quantity $\sqrt{2} p(\v_c) - E(\v_c)$ is governed, at the limit $\varepsilon \to 0$, by
$$\sqrt{2} p(\v_c) - E(\v_c) \sim - \frac{\varepsilon^5}{4} E_{KdV}(N),$$
where $E_{KdV}(N)$ is the energy of the soliton to \eqref{KdV}, namely
$$E_{KdV}(N) \equiv \frac{1}{2} \int_\R (N')^2 - \int_\R N^3.$$

\section{Qualitative properties of travelling waves in higher dimensions}
\label{qualibat}

We next turn to finite energy travelling waves in higher dimensions, and describe a number of qualitative properties, which have been rigorously established so far. Many of the results in this section were already guessed in the quoted seminal papers of Jones, Putterman and Roberts \cite{JoneRob1,JonPuRo1}. Some mathematical proofs turn out to be quite different from the physical intuition.

\subsection{Range of speeds and energies}

The following results are proved in \cite{Graveja2,Graveja4}. They describe the possible spectrum of speeds $c$, excluding in particular the possibility of non-constant supersonic travelling waves.

\begin{theorem}[\cite{Graveja2,Graveja4}]
\label{supersonic}
i) Let $N \geq 2$. Any travelling wave of finite energy and of supersonic speed $c > \sqrt{2}$ is constant.\\
ii) Let $N = 2$. Any travelling wave of finite energy and of sonic speed $c = \sqrt{2}$ is constant.
\end{theorem}

We briefly sketch the proof of Theorem \ref{supersonic} in the supersonic case. The proof is similar in the sonic case.

\begin{proof}
Let $v$ be a finite energy solution to \eqref{TWc} of speed $c > \sqrt{2}$. In order to prove that $v$ is a constant function, we compute some integral identities (relating the energy and the scalar momentum), so that the energy of $v$ is necessarily equal to $0$, i.e. $v$ is a constant function of modulus one. The first identities are the so-called Pohozaev's identities \cite{Pohozae1}, obtained by multiplying \eqref{TWc} by the test function $\langle x, \nabla v(x) \rangle$, and integrating by parts.

\begin{lemma}
\label{pohozaev}
Let $N \geq 1$ and $c \geq 0$, and consider a finite energy solution $v$ to \eqref{TWc}. Then,
\begin{equation}
\label{poho1}
E(v) = \int_{\R^N} |\partial_1 v|^2,
\end{equation}
and for any $2 \leq j \leq N$,
\begin{equation}
\label{pohoperp}
E(v) = \int_{\R^N} |\partial_j v|^2 + c p(v).
\end{equation}
\end{lemma}

Identities \eqref{poho1} and \eqref{pohoperp} are not sufficient to ensure that the energy of $v$ is equal to $0$
\footnote{Notice that \eqref{poho1} and \eqref{pohoperp} are sufficient to ensure that any travelling wave of finite energy and of speed $c = 0$ is constant (see \cite{JoneRob1,BethSau1}).}
, so that another integral relation is required.

\begin{lemma}
\label{newidentity}
Let $N \geq 2$ and $c > \sqrt{2}$, and consider a finite energy solution $v$ to \eqref{TWc}. Then,
\begin{equation}
\label{newone}
\int_{\R^N} \Big( |\nabla v|^2 + \big( 1 - |v|^2 \big)^2 \Big) = 2 c \Big( 1 - \frac{2}{c^2} \Big) p(v).
\end{equation}
\end{lemma}

Identity \eqref{newone} is obtained using the singularities of the kernels associated to \eqref{TWc} in case $c > \sqrt{2}$. Indeed, assuming that $v$ does not vanish on $\R^N$ (in order to simplify the proof), $v$ may be written as
$$v = \varrho \exp i \varphi,$$
where $\varrho = |v|$, and $\varphi$ are smooth, real-valued functions. By \eqref{TWc}, the function $\eta = 1 - \varrho^2$ is solution to
\begin{equation}
\label{equeta}
\Delta^2 \eta - 2 \Delta \eta + c^2 \partial^2_{1,1} \eta = - \Delta F + 2c \partial_1 \rm{div} (G),
\end{equation}
where
\begin{equation}
\label{FG}
F = 2|\nabla v|^2 + 2\eta^2 - 2 c \eta \partial_1 \varphi, \ {\rm and} \ G = - \eta \nabla \varphi.
\end{equation}
Notice in particular that the nonlinearities $F$ and $G$ are related to the density of energy and momentum (equal to $- G$ in view of \eqref{prost}). Taking the Fourier transform of \eqref{equeta}, we are led to
\begin{equation}
\label{etafou}
\widehat{\eta}(\xi) = \widehat{K_0}(\xi) \widehat{F}(\xi) - 2 c \sum_{j = 1}^N \widehat{K_j}(\xi) \widehat{G}(\xi),
\end{equation}
where $\widehat{K_0}(\xi) = \frac{|\xi|^2}{|\xi|^4 + 2 |\xi|^2 - c^2 \xi_1^2}$, and $\widehat{K_j}(\xi) = \frac{\xi_1 \xi_j}{|\xi|^4 + 2 |\xi|^2 - c^2 \xi_1^2}$, for any $1 \leq j \leq N$. The singular nature of the kernels $K_0$ and $K_j$ now gives some relation between $\widehat{F}(0)$ and $\widehat{G}(0)$. This in turn yields a relation between the energy and the scalar momentum. Indeed, since $v$ is of finite energy, the function $\widehat{\eta}$ in the left-hand side of \eqref{etafou} belongs to $L^2(\R^N)$, whereas the functions $\widehat{F}$ and $\widehat{G}$ in the right-hand side are continuous. When $0 \leq c < \sqrt{2}$, $\widehat{K_0}$ and $\widehat{K_j}$ are sufficiently smooth, so that \eqref{etafou} can hold without additional assumption. In contrast, when $c > \sqrt{2}$, $\widehat{K_0}$ and $\widehat{K_j}$ have singularities on the set
$$\Gamma = \big\{ \xi \in \R^N, |\xi|^4 + 2 |\xi|^2 - c^2 \xi_1^2 = 0 \big\}.$$
Since \eqref{etafou} holds, this leads to a relation between $\widehat{F}(\xi)$ and $\widehat{G}(\xi)$ on the set $\Gamma$. Taking some limit $\xi \to 0$, this yields identity \eqref{newone}. Combining \eqref{poho1}, \eqref{pohoperp} and \eqref{newone}, it may be shown that the energy of $v$ is equal to $0$, so that $v$ is a constant of modulus one. This completes the proof of Theorem \ref{supersonic}.
\end{proof}

In the three-dimensional subsonic case, small energy solutions are also excluded in view of the following theorem (see \cite{BetGrSa1}).

\begin{theorem}[\cite{BetGrSa1}]
\label{grenouille}
Let $N = 3$ and $0 < c < \sqrt{2}$. There exists some positive universal constant $\boE_0$ such that any non-constant finite energy solution to \eqref{TWc} on $\R^3$ satisfies
$$E(v) \geq \boE_0.$$
\end{theorem}

Theorem \ref{grenouille} improves an earlier result by Tarquini (see \cite{Tarquin1}), which states that any solutions of sufficiently small energy, with respect to their speed, are excluded in any dimension $N \geq 2$.

\begin{theorem}[\cite{Tarquin1}]
\label{tarquinlancien}
Let $N \geq 2$ and $0 < c < \sqrt{2}$. There exists some positive constant $\boE(c, N)$, depending only on $c$ and $N$, such that any non-constant finite energy solution to \eqref{TWc} on $\R^3$ satisfies
$$E(v) \geq \boE(c, N).$$
Moreover,
$$\boE(c, N) \to 0, \ {\rm as} \ c \to \sqrt{2}.$$
\end{theorem}

Theorem \ref{grenouille} leaves some hope for a complete scattering theory for solutions with small energy in the three-dimensional case. Such a theory has been established in any dimension $N \geq 4$ by Gustafson, Nakanishi and Tsai (see \cite{GusNaTs1}, or \cite{Nakanis1} in the present volume). In dimension three, they were able to establish the existence of global dispersive solutions to \eqref{GP} (see \cite{GusNaTs2} or \cite{Nakanis1} in the present volume). On the other hand, two-dimensional travelling waves of small energy are known to exist (see Theorem \ref{dim2}), excluding scattering in the energy space.

For sake of completeness, we briefly sketch the proof of Theorem \ref{grenouille}.

\begin{proof}
The proof relies on equation \eqref{etafou}, using the next elementary observation. If some quantity $E$ satisfies
\begin{equation}
\label{jagger}
E \leq K E^2,
\end{equation}
where $K$ is some positive constant, then $E$ is either equal to $0$, or
$$E \geq \frac{1}{K}.$$
Taking the $L^2$-norm of \eqref{etafou}, we are led to an inequality of the form \eqref{jagger}, where the energy $E(v)$ plays the role of the quantity $E$, whereas $K$ is equal to some linear combination of the $L^2$-integrals of the kernels $K_0$ and $K_j$.

Indeed, first notice that the nonlinearities $F$ and $G$ in the right-hand side of \eqref{etafou} are (almost) quadratic quantities with respect to $\eta$ and $\nabla \varphi$, related to the densities of energy and of momentum. In particular, their $L^1$-norms are bounded by the energy up to some positive universal constant $K$, so that
$$|\widehat{F}| + |\widehat{G}| \leq K E(v).$$
Taking the $L^2$-norm of \eqref{etafou}, we are led to
$$\| \eta \|_{L^2(\R^N)} = \| \widehat{\eta} \|_{L^2(\R^N)} \leq K E(v) \| \widehat{K_0} \|_{L^2(\R^N)}.$$
On the other hand, it is proved in \cite{BetGrSa1} that any finite energy solution $v$ to \eqref{TWc}, with the additional assumption $|v| \geq \frac{1}{2}$ (which holds in case $E(v)$ is sufficiently small (see \cite{BetGrSa1})), satisfies
$$E(v) \leq 7 c^2 \| \eta \|_{L^2(\R^3)}^2,$$
so that, by the above elementary observation,
\begin{equation}
\label{clapton}
E(v) \geq \frac{1}{7 c^2 K^2 \| \widehat{K_0} \|_{L^2(\R^N)}^2}.
\end{equation}
A direct computation now gives
$$\| \widehat{K_0} \|_{L^2(\R^3)}^2 = \frac{\pi^2}{c} \arcsin \Big( \frac{c}{\sqrt{2}} \Big),$$
so that
$$E(v) \geq \frac{1}{7 \pi^2 K^2 c \arcsin \big( \frac{c}{\sqrt{2}} \big)} \geq \boE_0 = \frac{\sqrt{2}}{7 \pi^3 K^2}.$$
This completes the proof of Theorem \ref{grenouille}.
\end{proof}

\begin{remark}
The proof above cannot be performed in the two-dimensional case. Indeed, a direct computation gives
$$\| \widehat{K_0} \|_{L^2(\R^2)}^2 = \frac{\pi}{\sqrt{2 (2 - c^2)}},$$
so that \eqref{clapton} now becomes
\begin{equation}
\label{goldman}
E(v) \geq \boE(c) = \frac{\sqrt{2 (2 - c^2)}}{7 \pi K^2 c^2},
\end{equation}
where $\boE(c) \to 0$, as $c \to \sqrt{2}$. In particular, \eqref{goldman} does not prevent the existence of solutions with energy as small as possible (see Theorem \ref{dim2}).
\end{remark}

\subsection{Regularity and decay at infinity}

The following results were proved in \cite{Farina1} and \cite{BetGrSa1} (see also \cite{BethSau1,Graveja3,Tarquin1,Tarquin2}). They describe the regularity of subsonic travelling waves (in particular their real-analyticity).

\begin{theorem}[\cite{Farina1,BetGrSa1}]
\label{tarquinlejeune}
Let $N \geq 2$ and $0 \leq c < \sqrt{2}$, and consider a finite energy solution $v$ to \eqref{TWc}. Then, $v$ is a real-analytic, bounded function on $\R^N$ such that
\begin{equation}
\label {elinfini}
|v| \leq \sqrt{1 + \frac{c^2}{4}}.
\end{equation}
\end{theorem}

For sake of completeness, we briefly sketch the proof of Theorem \ref{tarquinlejeune}.

\begin{proof}
The smoothness of $v$ follows from a standard bootstrap argument using the finiteness of the energy, and the elliptic nature of \eqref{TWc}. Bound \eqref{elinfini} essentially results from the maximum principle for the function $|v|^2$, which verifies
$$\Delta |v|^2 + 2 |v|^2 \Big( 1 + \frac{c^2}{4} - |v|^2 \Big) = 2 |\nabla v|^2 - 2 c \langle i \partial_1 v \ , v \rangle + \frac{c^2}{2} |v|^2 \geq 0,$$
by \eqref{TWc} (see \cite{Farina1,Tarquin1,Tarquin2} for more details).

Real-analyticity is established following an argument of Bona and Li \cite{BonaLi1,BonaLi2} (see also \cite{KatoPip1,Maris1,Maris2}). The idea is to prove the uniform convergence of Taylor's series of $v$,
$$T_{v,x}(z) = \sum_{\alpha \in \N^N} \frac{\partial^\alpha v(x)}{\alpha!} (z - x)^\alpha,$$
on a complex neighbourhood of any arbitrary point $x \in \R^N$. The required estimates for the derivatives are provided by \eqref{TWc}, using standard $L^q$-multiplier theory, Sobolev's embedding theorem, Gagliardo-Nirenberg's inequality, and the superlinear nature of the nonlinearities.

Indeed, denoting $v_1 = \Re(v) - 1$ and $v_2 = \Im(v)$, equation \eqref{TWc} may be recast as
\begin{align}
\label{eq-fou-v1}
\partial^2_{jk} v_1 = H_{j,k}*F_1(v_1, v_2) - i c H_{1,j,k} * F_2(v_1, v_2),\\
\label{eq-fou-v2}
\partial^2_{jk} v_2 = i c H_{1,j,k}*F_1(v_1, v_2) + K_{j,k} * F_2(v_1, v_2),
\end{align}
where the nonlinearities $F_1$ and $F_2$ are defined from $\C^2$ to $\C$ by
$$F_1(z_1, z_2) = 3 z_1^2 + z_2^2 + z_1^3 + z_1 z_2^2, \ {\rm and} \ F_2(z_1, z_2) = 2 z_1 z_2 + z_1^2 z_2 + z_2^3,$$
and the kernels $H_{j,k}$, $H_{1,j,k}$ and $K_{j,k}$ are given by
$$\widehat{H_{j,k}}(\xi) = \frac{\xi_j \xi_k |\xi|^2}{|\xi|^4 + 2 |\xi|^2 - c^2 \xi_1^2}, \ \widehat{H_{1,j,k}}(\xi) = \frac{\xi_1 \xi_j \xi_k}{|\xi|^4 + 2 |\xi|^2 - c^2 \xi_1^2}, \ \widehat{K_{j,k}}(\xi) = \frac{\xi_j \xi_k (2 + |\xi|^2)}{|\xi|^4 + 2 |\xi|^2 - c^2 \xi_1^2},$$
for any $1 \leq j, k \leq N$. By a result from Lizorkin \cite{Lizorki1}, the kernels $H_{j,k}$, $H_{1,j,k}$ and $K_{j,k}$ are $L^q$-multipliers for any $1 < q < + \infty$. Coupled to \eqref{eq-fou-v1} and \eqref{eq-fou-v2}, this provides $L^q$-estimates of any derivatives of $v_1$ and $v_2$.

\begin{lemma}
\label{Lq-estim}
Let $1 \leq j, k \leq N$, $\alpha \in \N^N$ and $1 < q < + \infty$. There exists some positive number $K_1(q)$, possibly depending on $q$, but not on $\alpha$, such that
\begin{equation}
\label{estim-Lq}
\| \partial^\alpha \partial^2_{jk} v_1 \|_{L^q(\R^N)} + \| \partial^\alpha \partial^2_{jk} v_2 \|_{L^q(\R^N)} \leq K_1(q) \Big( \| \partial^\alpha F_1(v_1, v_2) \|_{L^q(\R^N)} + \| \partial^\alpha F_2(v_1, v_2) \|_{L^q(\R^N)} \Big).
\end{equation}
\end{lemma}

The second step is to transform $L^q$-estimates \eqref{estim-Lq} into uniform ones. This follows from Sobolev's embedding theorem and Gagliardo-Nirenberg's inequality.

\begin{lemma}
\label{unif-estim}
Let $1 \leq j \leq N$, $\alpha \in \N^N$ and $\frac{N}{2} < q < + \infty$. There exist some positive numbers $K_2(q)$ and $K_3(q)$, possibly depending on $q$, but not on $\alpha$, such that
\begin{equation}
\label{estim-unif}
\begin{split}
\| \partial^\alpha v_1 \|_{L^\infty(\R^N)} + \| \partial^\alpha v_2 \|_{L^\infty(\R^N)} \leq & K_2(q) F_q(\alpha),\\
\| \partial^\alpha \partial_j v_1 \|_{L^q(\R^N)} + \| \partial^\alpha \partial_j v_2 \|_{L^q(\R^N)} \leq & K_3(q) F_q(\alpha),
\end{split}
\end{equation}
where we have set
$$F_q(\alpha) = \underset{0 \leq \beta \leq \alpha}{\max} \Big( \| \partial^\beta F_1(v_1, v_2) \|_{L^q(\R^N)} + \| \partial^\beta F_2(v_1, v_2) \|_{L^q(\R^N)} \Big).$$
\end{lemma}

In view of \eqref{estim-unif}, the convergence of Taylor's series $T_{v_1, x}(z)$ and $T_{v_2, x}(z)$ reduces to the convergence of the series
$$S_{q,x_0}(z) = \underset{\alpha \in \N^N}{\sum} \frac{F_q(\alpha)}{\alpha !}| z - x_0|^{|\alpha|},$$
for $z$ sufficiently close to $x_0$, and $q$ suitably chosen. Combining the superlinear nature of $F_1$ and $F_2$ with estimates \eqref{estim-Lq} and \eqref{estim-unif}, an inductive argument based on Abel's identity as in \cite{Maris1}, yields

\begin{lemma}
\label{F-alpha-estim}
Let $\alpha \in \N^N$ and $\frac{N}{2} < q < + \infty$. There exists some positive number $K_4(q)$, possibly depending on $q$, but not on $\alpha$, such that
\begin{equation}
\label{estim-F-alpha}
F_q(\alpha) \leq K_4(q)^{|\alpha|} \alpha^{\tilde{\alpha}},
\end{equation}
where we have set $\tilde{\alpha} = (\max \{ \alpha_1 - 1, 0 \}, \ldots, \max \{ \alpha_N - 1, 0 \})$.
\end{lemma}

In view of \eqref{estim-F-alpha}, choosing $q = N$, the series $S_{q,x_0}(z)$ is convergent for any $z$ such that $|z - x_0| < \frac{e}{K_4(N)}$. Taylor's series $T_{v_1, x_0}(z)$ and $T_{v_2, x_0}(z)$ converge the same way, so that $v$ is a real-analytic function on $\R^N$.
\end{proof}

Jones, Putterman and Roberts \cite{JoneRob1,JonPuRo1} investigated the decay properties of subsonic travelling waves in dimensions two and three. They computed some formal asymptotics of axisymmetric solutions to \eqref{TWc} in \cite{JoneRob1,JonPuRo1}. Their formal derivation was confirmed, and somewhat extended, in \cite{Graveja1,Graveja3,Graveja5,Graveja6}.

\begin{theorem}[\cite{Graveja1,Graveja3,Graveja5,Graveja6}]
\label{Decay}
Let $N \geq 2$ and $0 \leq c < \sqrt{2}$, and consider a finite energy solution $v$ to \eqref{TWc}. There exist a complex number $\lambda_\infty$, such that $|\lambda_\infty| = 1$,
and a smooth, real-valued function $v_\infty : \S^{N - 1} \to \R$, such that
\begin{equation}
\label{asymptotics}
|x|^{N-1} \Big( v(x) - \lambda_\infty \Big) - i \lambda_\infty v_\infty \Big( \frac{x}{|x|} \Big) \to 0 \ {\rm in} \ L^\infty(\S^{N - 1}), \ {\rm as} \ |x| \to + \infty.
\end{equation}
Moreover, there exist some real constants $\alpha$, $\beta_2$, $\ldots$, and $\beta_N$ such that the function $v_\infty$ is equal to
\begin{equation}
\label{v-infini}
v_\infty(\sigma) = \alpha \frac{\sigma_1}{\Big( 1 - \frac{c^2}{2} + \frac{c^2\sigma_1^2}{2} \Big)^\frac{N}{2}} + \sum_{j = 2}^N \beta_j \frac{\sigma_j}{\Big( 1 - \frac{c^2}{2} + \frac{c^2\sigma_1^2}{2} \Big)^\frac{N}{2}}, \forall \sigma \in \S^{N - 1}.
\end{equation}
The constants $\alpha$ and $\beta_j$ are given by
\begin{align}
\label{alpha}
\alpha & = \frac{\Gamma(\frac{N}{2})}{2 \pi^\frac{N}{2}} \Big(1 - \frac{c^2}{2} \Big)^\frac{N-3}{2} \bigg( \frac{4-N}{2} c E(v) + \Big( 2 + \frac{N-3}{2} c^2 \Big) p(v) \bigg),\\
\label{beta}
\beta_j & = \frac{\Gamma(\frac{N}{2})}{\pi^\frac{N}{2}} \Big(1 - \frac{c^2}{2} \Big)^\frac{N-1}{2} P_j(v),
\end{align}
where $P_j(v) = \frac{1}{2} \int_{\R^n} \langle i \partial_j v, v - 1 \rangle$. The constants $\beta_j$ are equal to $0$, when $v$ is axisymmetric around axis $x_1$.
\end{theorem}

Notice that a finite energy solution $v$ tends to some constant $\lambda_\infty$ at infinity. Since multiplication by a constant of modulus $1$ keeps \eqref{TWc} invariant, we may, without loss of generality, assume that $\lambda_\infty = 1$. Notice also that the decay of the function $v - 1$ is algebraic. More precisely, $v(x) - 1$ decays as $\frac{1}{|x|^{N - 1}}$, as $|x| \to + \infty$. This yields the following corollary of Theorem \ref{Decay}.

\begin{cor}[\cite{Graveja3,Graveja5}]
\label{Lp-decay}
Let $N \geq 2$ and $0 \leq c < \sqrt{2}$, and consider a finite energy solution $v$ to \eqref{TWc}.\\
i) The function $v - 1$ belongs to $L^q(\R^N)$ for any $q > \frac{N}{N - 1}$, and its gradient $\nabla v$ is in $L^q(\R^N)$ for any $q > 1$. Moreover, any higher order derivative $\partial^\alpha v$ belongs to $L^q(\R^N)$ for any $q \geq 1$.\\
ii) Assume $N = 2$, and $v$ is non-constant and axisymmetric around axis $x_1$. Then, the function $v - 1$ does not belong to $L^2(\R^2)$. In particular, it does not belong to $H^1(\R^2)$.
\end{cor}

Corollary \ref{Lp-decay} has some significant consequences. It first leads to a rigorous definition of the momentum $\vec{P}(v)$ of a finite energy solution $v$ to \eqref{TWc}. Indeed, formula \eqref{VectP} makes sense, in any dimension, for $\Psi = v$, since $v - 1$ and $\nabla v$ respectively belong to $L^4(\R^N)$ and $L^\frac{4}{3}(\R^N)$. Actually, it makes sense for any $w \in W(\R^N) = \{ 1 \} + V(\R^N)$, where $V(\R^N)$ is defined by
\begin{equation}
\label{def-V}
V(\R^N) = \{ v: \R^N \mapsto \C, \ {\rm s.t.} \ (\nabla v, \Re(v)) \in L^2(\R^N)^2, \Im(v) \in L^4(\R^N) \ {\rm and} \ \nabla \Re(v) \in L^\frac{4}{3}(\R^N) \}.
\end{equation}
Notice that
$$W(\R^N) \subset \boE(\R^N),$$
where
$$\boE(\R^N) = \{ v \in H_{\rm loc}^1(\R^N), \ {\rm s.t.} \ E(v) < + \infty \},$$
denotes the energy space. This last observation plays an important role in the variational argument of \cite{BetGrSa1} for the construction of finite energy solutions to \eqref{TWc} (see Section \ref{existence} below). Our choice of the variational space is indeed $W(\R^N) = \{ 1 \} + V(\R^N)$.

Statement ii) of Corollary \ref{Lp-decay} has to be considered in connection with the stability problem of two-dimensional travelling waves $v$. This analysis requires to find some functional space, on which the Gross-Pitaevskii equation is known to be globally well-posed, which preserves both the energy and the momentum, and which contains the travelling waves. As a matter of fact, there are several functional spaces where the two-dimensional Gross-Pitaevskii equation is known to be globally well-posed: for instance in $\{ 1 \} + H^1(\R^2)$ (see \cite{BethSau1,Gerard2}), in the energy space $\boE(\R^N)$ (see \cite{Gerard1,Gerard2}), and also in the space $\{ w \} + H^1(\R^N)$ (see \cite{Gallo1,Gerard2}), for any $w$ in the energy space $\boE(\R^N)$. Taking $w = v$, the space $\{ v \} + H^1(\R^2)$ seems appropriate to address the question of stability near a travelling wave $v$. Two important advantages of the space $\{ v \} + H^1(\R^2)$ are that $v$ belongs to this space (in contrast with $\{ 1 \} + H^1(\R^2)$, in view of Corollary \ref{Lp-decay}), and that the momentum is well-defined (in contrast with the energy space).

For sake of completeness, we briefly sketch the proof of Theorem \ref{Decay}.

\begin{proof}
The proof is reminiscent of a series of articles by Bona and Li \cite{BonaLi1}, De Bouard and Saut \cite{deBoSau2}, and Maris \cite{Maris1,Maris2}. It relies on the use of convolution equations, in particular on a careful analysis of the kernels they involve. Indeed, assuming that $v$ does not vanish on $\R^N$ (in order to simplify the proof), $v$ may be written as $v = \varrho \exp i \varphi$, where $\varrho$ and $\varphi$ are smooth, real-valued functions. By \eqref{TWc}, the functions $\eta = 1 - \varrho^2$ and $\nabla \varphi$ are solutions to \eqref{equeta}, and
\begin{equation}
\label{eqophi}
\Delta \varphi = \frac{c}{2} \partial_1 \eta + \div(G),
\end{equation}
where $F$ and $G$ are the nonlinearities given by \eqref{FG}. Taking the Fourier transform of \eqref{equeta} and \eqref{eqophi}, equation \eqref{TWc} reduces to the convolution equations
\begin{align*}
\eta = & K_0*F + 2 c \sum_{j = 1}^N K_j*G_j,\\
\partial_j \varphi = & \frac{c}{2} K_j*F + c^2 \sum_{k = 1}^N L_{j,k}*G_k + \sum_{k = 1}^N R_{j, k}*G_k,
\end{align*}
where $K_0$, $K_j$, $L_{j,k}$ and $R_{j,k}$ are the kernels given by
\begin{equation}
\label{kernels}
\begin{split}
\widehat{K_0}(\xi) = & \frac{|\xi|^2}{|\xi|^4 + 2 |\xi|^2 - c^2\xi_1^2}, \widehat{K_j}(\xi) = \widehat{R_{1, j}}(\xi) \widehat{K_0}(\xi),\\
\widehat{R_{j,k}}(\xi) = & \frac{\xi_j\xi_k}{|\xi|^2}, \ {\rm et} \ \widehat{L_{j, k}}(\xi) = \widehat{R_{1, j}}(\xi) \widehat{R_{1, k}}(\xi) \widehat{K_0}(\xi).
\end{split}
\end{equation}

Both the decay properties and the asymptotics of $v$ follow from the decay properties and the asymptotics of the kernels given by \eqref{kernels}. The first step is to derive them. In view of \eqref{kernels}, this reduces to the analysis of $K_0$ and $R_{j, k}$.

The kernels $R_{j, k}$ are directly related to Riesz's operators $R_j$, given by
$$\widehat{R_j}(\xi) = - i \frac{\xi_j}{|\xi|},$$
since they may be written as
$$R_{j, k} = - R_j * R_k.$$
The properties of Riesz's operators are well-known (see \cite{Stein1}), partly because of their explicit expression
$$R_j(x) = \frac{\Gamma \Big( \frac{N + 1}{2} \Big)}{\pi^\frac{N + 1}{2}} {\rm p.v.} \bigg( \frac{x_j}{|x|^{N + 1}} \bigg),$$
where p.v. stands for the principal value. There is a similar expression for the kernels $R_{j, k}$,
\begin{equation}
\label{explicit}
R_{j, k}(x) = \frac{\Gamma \Big( \frac{N}{2} \Big)}{2 \pi^\frac{N}{2}} {\rm p.v.} \bigg( \frac{\delta_{j, k} |x|^2 - N x_j x_k}{|x|^{N + 2}} \bigg),
\end{equation}
so that their asymptotics are straightforward.

In view of its non-homogeneity, the analysis of $K_0$ is more involved. However, we have

\begin{lemma}[\cite{Graveja3}]
\label{KO}
Let $N \geq 2$ and $N - 2 < \alpha \leq N$. There exists some positive constant $K_\alpha$ such that
$$|K_0(x)| \leq \frac{K_\alpha}{|x|^\alpha}, \ \forall x \in \R^N.$$
In particular, $K_0$ belongs to $L^q(\R^N)$ for any $\frac{N}{N - 1} < q < \frac{N}{N - 2}$.
\end{lemma}

Lemma \ref{KO} results from some integrability estimates of the derivatives of $\widehat{K_0}$. They provide the algebraic decay of $K_0$ by the inverse Fourier transform formula. Notice that this formula, coupled with some stationary phase estimates, states that the asymptotic properties of $K_0$ are mainly given by the behaviour of $\widehat{K_0}$ close to the origin. In view of \eqref{kernels}, this behaviour is given by
$$\widehat{K_0}(\xi) \sim \frac{|\xi|^2}{(2 - c^2) \xi_1^2 + 2 |\xi_\perp|^2} = \frac{1}{2 - c^2} \widehat{R_{1, 1}} \Big( \sqrt{2 - c^2} \xi_1, \sqrt{2} \xi_\perp \Big) + \frac{1}{2} \sum_{j = 1}^N \widehat{R_{j, j}} \Big( \sqrt{2 - c^2} \xi_1, \sqrt{2} \xi_\perp \Big).$$
Up to some scale changes, the asymptotics of $K_0$ are therefore the same as the asymptotics of $R_{j, j}$. This gives for instance,

\begin{lemma}[\cite{Graveja5}]
\label{KOalalimite}
Assume $N \geq 2$ and consider a smooth function $f$ such that
\begin{equation}
\label{decay}
|f(x)| \leq \frac{K}{1 + |x|^{2 N}}, \ {\rm and} \ |\nabla f(x)| \leq \frac{K}{1 + |x|^{2 N + 1}}, \ \forall x \in \R^N,
\end{equation}
where $K$ is some positive constant. Then,
\begin{equation}
\label{koalalimite}
R^N K_0*f(R \sigma) \to \frac{\Gamma \big( \frac{N}{2} \big) c^2 \big(2 - c^2 \big)^{\frac{N - 3}{2}} \Big( 2 - c^2 + (c^2 -2 N) \sigma_1^2 \Big)}{2 \sqrt{2} \pi^\frac{N}{2} \Big( 2 - c^2 + c^2 \sigma_1^2 \Big)^\frac{N + 2}{2}} \bigg( \int_{\R^N} f \bigg), \ {\rm as} \ R \to + \infty.
\end{equation}
\end{lemma}

The asymptotics of a travelling wave $v$ given by formulae \eqref{asymptotics}, \eqref{v-infini}, \eqref{alpha} and \eqref{beta} are derived using formulae like \eqref{koalalimite} for the kernels $K_0$, $K_j$, $L_{j, k}$ and $R_{j, k}$, and the nonlinearities $F$ and $G$. Notice in particular that the presence of the energy and the scalar momentum in formulae \eqref{alpha} and \eqref{beta} is due to the presence of the integrals of $F$ and $G$ in formula \eqref{koalalimite} above (see the proof of Theorem \ref{supersonic} above).

Notice also that the proof of \eqref{koalalimite} first requires to establish \eqref{decay}, that is some algebraic decay for the function $f$ convoluted to $K_0$, and its gradient $\nabla f$. This is the second step of the proof. As in earlier papers \cite{BonaLi1,deBoSau2,Maris1,Maris2}, this relies on some inductive argument using the superlinear nature of the nonlinearities $F$ and $G$, which are almost quadratic functions of the variables $\eta$ and $\nabla \varphi$ in view of \eqref{FG}. When the nonlinearities are superlinear, the algebraic decay of a solution is identical to the algebraic decay of the kernels of the convolution equations it satisfies.

To get a feeling for this claim, let us consider the simplified model
\begin{equation}
\label{toy}
f = K*f^r,
\end{equation}
for which we have

\begin{lemma}[\cite{Graveja3}]
\label{Induction}
Assume $r > 1$, and consider smooth solutions $f$ and $K$ to \eqref{toy}, such that $f$ belongs to $L^r(\R^N)$, and $K$ is in $L^1(\R^N)$, and satisfies
$$|K(x)| \leq \frac{A}{1 + |x|^{\alpha_K}}.$$
If there exists some positive constant $0 < \alpha_0 < \alpha_K$ such that
$$|f(x)| \leq \frac{A}{|x|^{\alpha_0}},$$
then,
$$|f(x)| \leq \frac{A}{|x|^{\alpha_K}},$$
where $A$ denotes some positive, possibly different, constants.
\end{lemma}

Lemma \ref{Induction} states that the algebraic decay of $f$ is exactly the same as the algebraic decay of $K$, if some small decay may be first established for $f$. Its assumptions are quite restrictive, but may be extended to more involved equations with only additional technicalities. The proof is by induction. Indeed, by \eqref{toy},
\begin{equation}
\label{E43}
\begin{split}
|x|^\alpha |f(x)| \leq & A \bigg( \int_{\R^N} |x-y|^\alpha |K(x-y)| |f(y)|^r dy + \int_{\R^N} |K(x-y)| |y|^\alpha |f(y)|^r dy \bigg) \\ \leq & A \bigg( \big\| |\cdot|^\alpha K \big\|_{L^\infty(\R^N)} \big\| f \big\|^r_{L^r(\R^N)} + \big\| K \big\|_{L^1(\R^N)} \big\| |\cdot|^\frac{\alpha}{r} f \big\|^r_{L^\infty(\R^N)} \bigg).
\end{split}
\end{equation}
Using the assumptions of Lemma \ref{Induction}, equation \eqref{E43} reduces to
\begin{equation}
\label{E44}
\big\| |\cdot|^\alpha f \big\|_{L^\infty(\R^N)} \leq A + A \| |\cdot|^\frac{\alpha}{r} f \|^r_{L^\infty(\R^N)},
\end{equation}
if $0 \leq \alpha \leq \alpha_K$. Equation \eqref{E44} now links the algebraic decay with exponent $\alpha$ of $f$ to its algebraic decay with exponent $\frac{\alpha}{q}$. In particular, if we know some algebraic decay with a small exponent $\alpha_0 > 0$, a bootstrap argument yields the algebraic decay of $f$ for $\alpha = q \alpha_0$, $\alpha = q^2 \alpha_0$, $\ldots$, that is for any $\alpha \in [\alpha_0, \alpha_K]$. This provides the optimal decay of $f$, which is the decay of the kernel $K$.

However, in order to perform the inductive argument of Lemma \ref{Induction}, and to get the decay of a travelling wave $v$, we must establish some small decay for $v$. Multiplying equations \eqref{equeta} and \eqref{eqophi} by $\eta$ and $\varphi$, and integrating by parts on the complementary of a ball, we may establish the following integral decay.

\begin{lemma}[\cite{BetOrSm1,Graveja3}]
\label{Nabot}
Let $N \geq 2$ and $0 \leq c < \sqrt{2}$, and consider a finite energy solution $v$ to \eqref{TWc}. There exists some positive constant $\alpha_0$ such that
$$\int_{\R^N} |x|^{\alpha_0} e(v)(x) dx < + \infty.$$
\end{lemma}

Lemma \ref{Nabot} gives the small decay required to use the inductive argument of Lemma \ref{Induction}. Coupled with the decay and integrability properties of the kernels (see Lemma \ref{KO} and formula \eqref{explicit}), this first gives the decay of a travelling wave $v$, then its asymptotics, using formulae like formula \eqref{koalalimite} of Lemma \ref{KOalalimite}. This completes the proof of Theorem \ref{Decay}.
\end{proof}

\section{ The existence problem in higher dimensions}
\label{existence}

\subsection{The variational approach}

As mentioned in the introduction, the existence problem in dimensions two and three has been widely investigated in the physical literature: rigorous mathematical proofs have been provided more recently, using a variational approach.

It is a classical observation that one may obtain travelling waves by minimizing the energy $E$ keeping the momentum $p$ fixed. For a given $\p \geq 0$, we therefore consider the minimization problem
$$E_{\min}(\p) = \inf \{ E(v), v \in W(\R^N), p(v) = \p \},$$
where the definition of the space $W(\R^N)$ is given by \eqref{def-V}. In this approach, equation \eqref{TWc} is the Euler-Lagrange equation to this constrained minimization problem. The speed $c$ appears as a Lagrange multiplier, and is therefore not fixed a priori. Instead of using minimization under constraint, an alternate approach is to introduce, for given $c$, the Lagrangian
$$F_c(v) = E(v) - c p(v),$$
whose critical points are solutions to \eqref{TWc}. Solutions may then be found applying a mountain-pass argument. As a matter of fact, both techniques have been used so far. The first existence results were based on asymptotic Ginzburg-Landau's theory: the solutions obtained in that context had small speeds $c > 0$ and possessed vorticity, i.e zeroes with non-trivial topological degree. For instance, in dimension two, one may obtain, using the mountain-pass lemma, a branch of solutions parametrized by the speed $c > 0$.

\begin{theorem}[\cite{BethSau1}]
\label{mpdim2}
Assume $N = 2$. There exists $c_0 > 0$ such that there exists a non-constant finite energy solution $v$ to \eqref{TWc} for any $0 < c < c_0$. Moreover, $v$ exactly has two vortices with degree $\pm1$ located at a distance $\sim \frac{2}{c}, $ as $c \to 0$, whereas
$$ E(v) \sim 2 \pi \ln(c), \ {\rm and} \ p(v) \sim \frac{2 \pi}{c}, \ {\rm as} \ c \to 0.$$
\end{theorem}

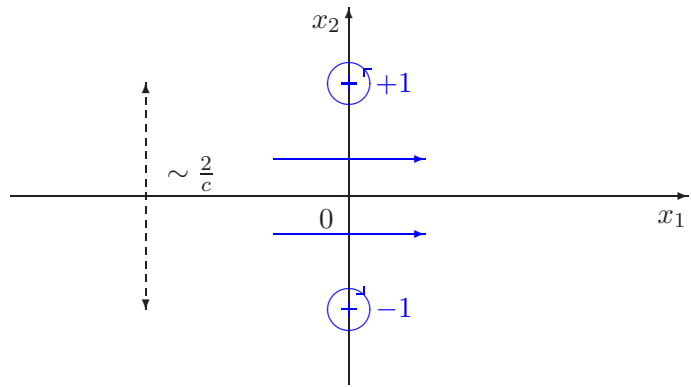
\begin{center}
\begin{picture}(90,50)(0,0)
\linethickness{0.1mm}
\put(0,25){\line(1,0){90}}
\put(90,25){\vector(1,0){0.12}}
\linethickness{0.1mm}
\put(45,0){\line(0,1){50}}
\put(45,50){\vector(0,1){0.12}}
{\color{blue}
\linethickness{0.2mm}
\put(44,40){\line(1,0){2}}
\put(45,39){\line(0,1){2}}
\put(44,10){\line(1,0){2}}
\put(45,9){\line(0,1){2}}
\linethickness{0.1mm}
\put(45,40){\circle{6}}
\put(45,10){\circle{6}}
\put(47,41){\line(0,1){1}}
\put(47,42){\line(1,0){1}}
\put(47,12){\line(0,1){1}}
\put(46,12){\line(1,0){1}}
\put(35,30){\line(1,0){20}}
\put(55,30){\vector(1,0){0.12}}
\put(35,20){\line(1,0){20}}
\put(55,20){\vector(1,0){0.12}}
}
\linethickness{0.1mm}
\multiput(18,10)(0,1.82){17}{\line(0,1){0.91}}
\put(18,40){\vector(0,1){0.12}}
\put(18,10){\vector(0,-1){0.12}}
\put(42,22){\makebox(0,0)[cc]{$0$}}
\put(42,48){\makebox(0,0)[cc]{$x_2$}}
\put(88,22){\makebox(0,0)[cc]{$x_1$}}
\put(24,28){\makebox(0,0)[cc]{$\sim \frac{2}{c}$}}
{\color{blue}
\put(51,10){\makebox(0,0)[cc]{$- 1$}}
\put(51,40){\makebox(0,0)[cc]{$+ 1$}}
}
\end{picture}
\end{center}

In dimension three, the constrained minimization approach provides the following theorem for large momentums.

\begin{theorem}[\cite{BetOrSm1}]
\label{min3d}
Assume $N = 3$. There exists $\p_\infty > 0$ such that, for any $\p \geq \p_\infty$, there exists a solution $\u_\p$ to \eqref{TWc}, with $c = c(\u_\p)$, verifying $p(\u_\p) = \p$,
$$E(\u_\p) \sim \pi \sqrt{\p} \ln(\p), \ {\rm and} \ c(\u_\p) \sim \frac{\pi \ln(\p)}{2 \sqrt{\p}}, \ {\rm as} \ \p \to + \infty.$$
Moreover, $v$ presents a vortex ring whose diameter is $\sim \frac{2 \sqrt{\p}}{\pi}, $ as $\p \to + \infty$.
\end{theorem}

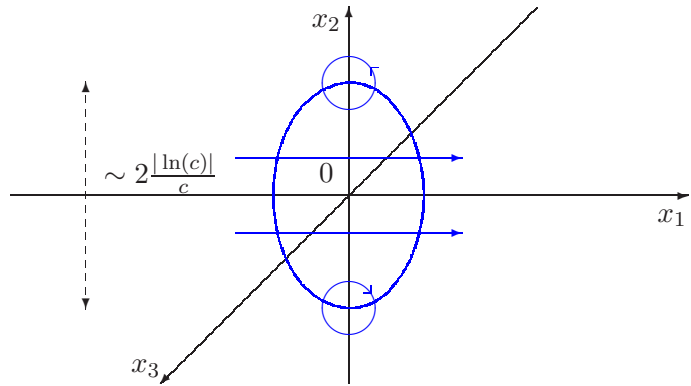
\begin{center}
\begin{picture}(90,50)(0,0)
\linethickness{0.1mm}
\put(0,25){\line(1,0){90}}
\put(90,25){\vector(1,0){0.12}}
\put(45,0){\line(0,1){50}}
\put(45,50){\vector(0,1){0.12}}
\put(20,0){\vector(-1,-1){0.12}}
\multiput(20,0)(0.12,0.12){417}{\line(1,0){0.12}}

{\color{blue}
\linethickness{0.2mm}
\put(55,24.75){\line(0,1){0.5}}
\multiput(54.99,25.75)(0.01,-0.5){1}{\line(0,-1){0.5}}
\multiput(54.97,26.25)(0.02,-0.5){1}{\line(0,-1){0.5}}
\multiput(54.93,26.75)(0.03,-0.5){1}{\line(0,-1){0.5}}
\multiput(54.89,27.25)(0.04,-0.5){1}{\line(0,-1){0.5}}
\multiput(54.83,27.74)(0.06,-0.49){1}{\line(0,-1){0.49}}
\multiput(54.76,28.23)(0.07,-0.49){1}{\line(0,-1){0.49}}
\multiput(54.69,28.72)(0.08,-0.49){1}{\line(0,-1){0.49}}
\multiput(54.6,29.2)(0.09,-0.48){1}{\line(0,-1){0.48}}
\multiput(54.5,29.68)(0.1,-0.48){1}{\line(0,-1){0.48}}
\multiput(54.39,30.16)(0.11,-0.47){1}{\line(0,-1){0.47}}
\multiput(54.27,30.62)(0.12,-0.47){1}{\line(0,-1){0.47}}
\multiput(54.14,31.09)(0.13,-0.46){1}{\line(0,-1){0.46}}
\multiput(54,31.54)(0.14,-0.45){1}{\line(0,-1){0.45}}
\multiput(53.85,31.99)(0.15,-0.45){1}{\line(0,-1){0.45}}
\multiput(53.69,32.43)(0.16,-0.44){1}{\line(0,-1){0.44}}
\multiput(53.52,32.86)(0.17,-0.43){1}{\line(0,-1){0.43}}
\multiput(53.34,33.28)(0.18,-0.42){1}{\line(0,-1){0.42}}
\multiput(53.15,33.69)(0.09,-0.21){2}{\line(0,-1){0.21}}
\multiput(52.95,34.1)(0.1,-0.2){2}{\line(0,-1){0.2}}
\multiput(52.74,34.49)(0.1,-0.2){2}{\line(0,-1){0.2}}
\multiput(52.53,34.87)(0.11,-0.19){2}{\line(0,-1){0.19}}
\multiput(52.3,35.25)(0.11,-0.19){2}{\line(0,-1){0.19}}
\multiput(52.07,35.61)(0.12,-0.18){2}{\line(0,-1){0.18}}
\multiput(51.83,35.96)(0.12,-0.17){2}{\line(0,-1){0.17}}
\multiput(51.58,36.29)(0.12,-0.17){2}{\line(0,-1){0.17}}
\multiput(51.33,36.62)(0.13,-0.16){2}{\line(0,-1){0.16}}
\multiput(51.07,36.93)(0.13,-0.16){2}{\line(0,-1){0.16}}
\multiput(50.8,37.22)(0.13,-0.15){2}{\line(0,-1){0.15}}
\multiput(50.52,37.51)(0.14,-0.14){2}{\line(0,-1){0.14}}
\multiput(50.24,37.78)(0.14,-0.13){2}{\line(1,0){0.14}}
\multiput(49.95,38.03)(0.14,-0.13){2}{\line(1,0){0.14}}
\multiput(49.66,38.27)(0.15,-0.12){2}{\line(1,0){0.15}}
\multiput(49.36,38.5)(0.15,-0.11){2}{\line(1,0){0.15}}
\multiput(49.06,38.71)(0.15,-0.11){2}{\line(1,0){0.15}}
\multiput(48.75,38.91)(0.15,-0.1){2}{\line(1,0){0.15}}
\multiput(48.44,39.09)(0.16,-0.09){2}{\line(1,0){0.16}}
\multiput(48.12,39.25)(0.32,-0.16){1}{\line(1,0){0.32}}
\multiput(47.8,39.4)(0.32,-0.15){1}{\line(1,0){0.32}}
\multiput(47.48,39.53)(0.32,-0.13){1}{\line(1,0){0.32}}
\multiput(47.16,39.65)(0.33,-0.12){1}{\line(1,0){0.33}}
\multiput(46.83,39.75)(0.33,-0.1){1}{\line(1,0){0.33}}
\multiput(46.5,39.83)(0.33,-0.08){1}{\line(1,0){0.33}}
\multiput(46.17,39.9)(0.33,-0.07){1}{\line(1,0){0.33}}
\multiput(45.83,39.95)(0.33,-0.05){1}{\line(1,0){0.33}}
\multiput(45.5,39.98)(0.33,-0.03){1}{\line(1,0){0.33}}
\multiput(45.17,40)(0.33,-0.02){1}{\line(1,0){0.33}}
\put(44.83,40){\line(1,0){0.33}}
\multiput(44.5,39.98)(0.33,0.02){1}{\line(1,0){0.33}}
\multiput(44.17,39.95)(0.33,0.03){1}{\line(1,0){0.33}}
\multiput(43.83,39.9)(0.33,0.05){1}{\line(1,0){0.33}}
\multiput(43.5,39.83)(0.33,0.07){1}{\line(1,0){0.33}}
\multiput(43.17,39.75)(0.33,0.08){1}{\line(1,0){0.33}}
\multiput(42.84,39.65)(0.33,0.1){1}{\line(1,0){0.33}}
\multiput(42.52,39.53)(0.33,0.12){1}{\line(1,0){0.33}}
\multiput(42.2,39.4)(0.32,0.13){1}{\line(1,0){0.32}}
\multiput(41.88,39.25)(0.32,0.15){1}{\line(1,0){0.32}}
\multiput(41.56,39.09)(0.32,0.16){1}{\line(1,0){0.32}}
\multiput(41.25,38.91)(0.16,0.09){2}{\line(1,0){0.16}}
\multiput(40.94,38.71)(0.15,0.1){2}{\line(1,0){0.15}}
\multiput(40.64,38.5)(0.15,0.11){2}{\line(1,0){0.15}}
\multiput(40.34,38.27)(0.15,0.11){2}{\line(1,0){0.15}}
\multiput(40.05,38.03)(0.15,0.12){2}{\line(1,0){0.15}}
\multiput(39.76,37.78)(0.14,0.13){2}{\line(1,0){0.14}}
\multiput(39.48,37.51)(0.14,0.13){2}{\line(1,0){0.14}}
\multiput(39.2,37.22)(0.14,0.14){2}{\line(0,1){0.14}}
\multiput(38.93,36.93)(0.13,0.15){2}{\line(0,1){0.15}}
\multiput(38.67,36.62)(0.13,0.16){2}{\line(0,1){0.16}}
\multiput(38.42,36.29)(0.13,0.16){2}{\line(0,1){0.16}}
\multiput(38.17,35.96)(0.12,0.17){2}{\line(0,1){0.17}}
\multiput(37.93,35.61)(0.12,0.17){2}{\line(0,1){0.17}}
\multiput(37.7,35.25)(0.12,0.18){2}{\line(0,1){0.18}}
\multiput(37.47,34.87)(0.11,0.19){2}{\line(0,1){0.19}}
\multiput(37.26,34.49)(0.11,0.19){2}{\line(0,1){0.19}}
\multiput(37.05,34.1)(0.1,0.2){2}{\line(0,1){0.2}}
\multiput(36.85,33.69)(0.1,0.2){2}{\line(0,1){0.2}}
\multiput(36.66,33.28)(0.09,0.21){2}{\line(0,1){0.21}}
\multiput(36.48,32.86)(0.18,0.42){1}{\line(0,1){0.42}}
\multiput(36.31,32.43)(0.17,0.43){1}{\line(0,1){0.43}}
\multiput(36.15,31.99)(0.16,0.44){1}{\line(0,1){0.44}}
\multiput(36,31.54)(0.15,0.45){1}{\line(0,1){0.45}}
\multiput(35.86,31.09)(0.14,0.45){1}{\line(0,1){0.45}}
\multiput(35.73,30.62)(0.13,0.46){1}{\line(0,1){0.46}}
\multiput(35.61,30.16)(0.12,0.47){1}{\line(0,1){0.47}}
\multiput(35.5,29.68)(0.11,0.47){1}{\line(0,1){0.47}}
\multiput(35.4,29.2)(0.1,0.48){1}{\line(0,1){0.48}}
\multiput(35.31,28.72)(0.09,0.48){1}{\line(0,1){0.48}}
\multiput(35.24,28.23)(0.08,0.49){1}{\line(0,1){0.49}}
\multiput(35.17,27.74)(0.07,0.49){1}{\line(0,1){0.49}}
\multiput(35.11,27.25)(0.06,0.49){1}{\line(0,1){0.49}}
\multiput(35.07,26.75)(0.04,0.5){1}{\line(0,1){0.5}}
\multiput(35.03,26.25)(0.03,0.5){1}{\line(0,1){0.5}}
\multiput(35.01,25.75)(0.02,0.5){1}{\line(0,1){0.5}}
\multiput(35,25.25)(0.01,0.5){1}{\line(0,1){0.5}}
\put(35,24.75){\line(0,1){0.5}}
\multiput(35,24.75)(0.01,-0.5){1}{\line(0,-1){0.5}}
\multiput(35.01,24.25)(0.02,-0.5){1}{\line(0,-1){0.5}}
\multiput(35.03,23.75)(0.03,-0.5){1}{\line(0,-1){0.5}}
\multiput(35.07,23.25)(0.04,-0.5){1}{\line(0,-1){0.5}}
\multiput(35.11,22.75)(0.06,-0.49){1}{\line(0,-1){0.49}}
\multiput(35.17,22.26)(0.07,-0.49){1}{\line(0,-1){0.49}}
\multiput(35.24,21.77)(0.08,-0.49){1}{\line(0,-1){0.49}}
\multiput(35.31,21.28)(0.09,-0.48){1}{\line(0,-1){0.48}}
\multiput(35.4,20.8)(0.1,-0.48){1}{\line(0,-1){0.48}}
\multiput(35.5,20.32)(0.11,-0.47){1}{\line(0,-1){0.47}}
\multiput(35.61,19.84)(0.12,-0.47){1}{\line(0,-1){0.47}}
\multiput(35.73,19.38)(0.13,-0.46){1}{\line(0,-1){0.46}}
\multiput(35.86,18.91)(0.14,-0.45){1}{\line(0,-1){0.45}}
\multiput(36,18.46)(0.15,-0.45){1}{\line(0,-1){0.45}}
\multiput(36.15,18.01)(0.16,-0.44){1}{\line(0,-1){0.44}}
\multiput(36.31,17.57)(0.17,-0.43){1}{\line(0,-1){0.43}}
\multiput(36.48,17.14)(0.18,-0.42){1}{\line(0,-1){0.42}}
\multiput(36.66,16.72)(0.09,-0.21){2}{\line(0,-1){0.21}}
\multiput(36.85,16.31)(0.1,-0.2){2}{\line(0,-1){0.2}}
\multiput(37.05,15.9)(0.1,-0.2){2}{\line(0,-1){0.2}}
\multiput(37.26,15.51)(0.11,-0.19){2}{\line(0,-1){0.19}}
\multiput(37.47,15.13)(0.11,-0.19){2}{\line(0,-1){0.19}}
\multiput(37.7,14.75)(0.12,-0.18){2}{\line(0,-1){0.18}}
\multiput(37.93,14.39)(0.12,-0.17){2}{\line(0,-1){0.17}}
\multiput(38.17,14.04)(0.12,-0.17){2}{\line(0,-1){0.17}}
\multiput(38.42,13.71)(0.13,-0.16){2}{\line(0,-1){0.16}}
\multiput(38.67,13.38)(0.13,-0.16){2}{\line(0,-1){0.16}}
\multiput(38.93,13.07)(0.13,-0.15){2}{\line(0,-1){0.15}}
\multiput(39.2,12.78)(0.14,-0.14){2}{\line(0,-1){0.14}}
\multiput(39.48,12.49)(0.14,-0.13){2}{\line(1,0){0.14}}
\multiput(39.76,12.22)(0.14,-0.13){2}{\line(1,0){0.14}}
\multiput(40.05,11.97)(0.15,-0.12){2}{\line(1,0){0.15}}
\multiput(40.34,11.73)(0.15,-0.11){2}{\line(1,0){0.15}}
\multiput(40.64,11.5)(0.15,-0.11){2}{\line(1,0){0.15}}
\multiput(40.94,11.29)(0.15,-0.1){2}{\line(1,0){0.15}}
\multiput(41.25,11.09)(0.16,-0.09){2}{\line(1,0){0.16}}
\multiput(41.56,10.91)(0.32,-0.16){1}{\line(1,0){0.32}}
\multiput(41.88,10.75)(0.32,-0.15){1}{\line(1,0){0.32}}
\multiput(42.2,10.6)(0.32,-0.13){1}{\line(1,0){0.32}}
\multiput(42.52,10.47)(0.33,-0.12){1}{\line(1,0){0.33}}
\multiput(42.84,10.35)(0.33,-0.1){1}{\line(1,0){0.33}}
\multiput(43.17,10.25)(0.33,-0.08){1}{\line(1,0){0.33}}
\multiput(43.5,10.17)(0.33,-0.07){1}{\line(1,0){0.33}}
\multiput(43.83,10.1)(0.33,-0.05){1}{\line(1,0){0.33}}
\multiput(44.17,10.05)(0.33,-0.03){1}{\line(1,0){0.33}}
\multiput(44.5,10.02)(0.33,-0.02){1}{\line(1,0){0.33}}
\put(44.83,10){\line(1,0){0.33}}
\multiput(45.17,10)(0.33,0.02){1}{\line(1,0){0.33}}
\multiput(45.5,10.02)(0.33,0.03){1}{\line(1,0){0.33}}
\multiput(45.83,10.05)(0.33,0.05){1}{\line(1,0){0.33}}
\multiput(46.17,10.1)(0.33,0.07){1}{\line(1,0){0.33}}
\multiput(46.5,10.17)(0.33,0.08){1}{\line(1,0){0.33}}
\multiput(46.83,10.25)(0.33,0.1){1}{\line(1,0){0.33}}
\multiput(47.16,10.35)(0.33,0.12){1}{\line(1,0){0.33}}
\multiput(47.48,10.47)(0.32,0.13){1}{\line(1,0){0.32}}
\multiput(47.8,10.6)(0.32,0.15){1}{\line(1,0){0.32}}
\multiput(48.12,10.75)(0.32,0.16){1}{\line(1,0){0.32}}
\multiput(48.44,10.91)(0.16,0.09){2}{\line(1,0){0.16}}
\multiput(48.75,11.09)(0.15,0.1){2}{\line(1,0){0.15}}
\multiput(49.06,11.29)(0.15,0.11){2}{\line(1,0){0.15}}
\multiput(49.36,11.5)(0.15,0.11){2}{\line(1,0){0.15}}
\multiput(49.66,11.73)(0.15,0.12){2}{\line(1,0){0.15}}
\multiput(49.95,11.97)(0.14,0.13){2}{\line(1,0){0.14}}
\multiput(50.24,12.22)(0.14,0.13){2}{\line(1,0){0.14}}
\multiput(50.52,12.49)(0.14,0.14){2}{\line(0,1){0.14}}
\multiput(50.8,12.78)(0.13,0.15){2}{\line(0,1){0.15}}
\multiput(51.07,13.07)(0.13,0.16){2}{\line(0,1){0.16}}
\multiput(51.33,13.38)(0.13,0.16){2}{\line(0,1){0.16}}
\multiput(51.58,13.71)(0.12,0.17){2}{\line(0,1){0.17}}
\multiput(51.83,14.04)(0.12,0.17){2}{\line(0,1){0.17}}
\multiput(52.07,14.39)(0.12,0.18){2}{\line(0,1){0.18}}
\multiput(52.3,14.75)(0.11,0.19){2}{\line(0,1){0.19}}
\multiput(52.53,15.13)(0.11,0.19){2}{\line(0,1){0.19}}
\multiput(52.74,15.51)(0.1,0.2){2}{\line(0,1){0.2}}
\multiput(52.95,15.9)(0.1,0.2){2}{\line(0,1){0.2}}
\multiput(53.15,16.31)(0.09,0.21){2}{\line(0,1){0.21}}
\multiput(53.34,16.72)(0.18,0.42){1}{\line(0,1){0.42}}
\multiput(53.52,17.14)(0.17,0.43){1}{\line(0,1){0.43}}
\multiput(53.69,17.57)(0.16,0.44){1}{\line(0,1){0.44}}
\multiput(53.85,18.01)(0.15,0.45){1}{\line(0,1){0.45}}
\multiput(54,18.46)(0.14,0.45){1}{\line(0,1){0.45}}
\multiput(54.14,18.91)(0.13,0.46){1}{\line(0,1){0.46}}
\multiput(54.27,19.38)(0.12,0.47){1}{\line(0,1){0.47}}
\multiput(54.39,19.84)(0.11,0.47){1}{\line(0,1){0.47}}
\multiput(54.5,20.32)(0.1,0.48){1}{\line(0,1){0.48}}
\multiput(54.6,20.8)(0.09,0.48){1}{\line(0,1){0.48}}
\multiput(54.69,21.28)(0.08,0.49){1}{\line(0,1){0.49}}
\multiput(54.76,21.77)(0.07,0.49){1}{\line(0,1){0.49}}
\multiput(54.83,22.26)(0.06,0.49){1}{\line(0,1){0.49}}
\multiput(54.89,22.75)(0.04,0.5){1}{\line(0,1){0.5}}
\multiput(54.93,23.25)(0.03,0.5){1}{\line(0,1){0.5}}
\multiput(54.97,23.75)(0.02,0.5){1}{\line(0,1){0.5}}
\multiput(54.99,24.25)(0.01,0.5){1}{\line(0,1){0.5}}

\linethickness{0.1mm}
\put(30,30){\line(1,0){30}}
\put(60,30){\vector(1,0){0.12}}
\linethickness{0.1mm}
\put(30,20){\line(1,0){30}}
\put(60,20){\vector(1,0){0.12}}

\linethickness{0.1mm}
\put(45,40){\circle{7.07}}
\put(45,10){\circle{7.07}}
\put(48,41){\line(0,1){1}}
\put(48,42){\line(1,0){1}}
\put(48,12){\line(0,1){1}}
\put(47,12){\line(1,0){1}}
}

\linethickness{0.1mm}
\multiput(10,10)(0,1.82){17}{\line(0,1){0.91}}
\put(10,40){\vector(0,1){0.12}}
\put(10,10){\vector(0,-1){0.12}}

\linethickness{0.1mm}
\put(42,28){\makebox(0,0)[cc]{$0$}}
\put(88,22){\makebox(0,0)[cc]{$x_1$}}
\put(42,48){\makebox(0,0)[cc]{$x_2$}}
\put(18,2){\makebox(0,0)[cc]{$x_3$}}
\put(20,28){\makebox(0,0)[cc]{$\sim 2 \frac{|\ln(c)|}{c}$}}
\end{picture}
\end{center}

As a matter of fact, the mountain-pass lemma can also be used in dimension three to assert the existence of a branch parametrized by the speed.

\begin{theorem}[\cite{Chiron1}]
\label{mpdim3} 
Assume $N = 3$. There exists $c_0 > 0$ such that there exists a non-constant finite energy solution $\u_c$ to \eqref{TWc} for any $0 < c < c_0$. Moreover, $v$ presents a vortex ring whose diameter is $\sim \frac{2 |\ln(c)|}{c}$, as $c \to 0$, and
$$E(\u_c) \sim \frac{2 \pi^2 |\ln(c)|^2}{c}, \ {\rm and} \ p(\u_c) \sim \frac{\pi^2 \ln(c)^2}{c^2}, \ {\rm as} \ c \to 0.$$
\end{theorem}

\begin{remark}
One may conjecture that both approaches exactly yield the same solutions.
\end{remark}

In \cite{BetGrSa1}, we revisit the minimization under constraint method, and are able to construct the full branch of minimizers parametrized by $\p$ in both dimensions.

\begin{theorem}[\cite{BetGrSa1}]
\label{dim2}
Assume $N = 2$ and $\p > 0$. There exists a non-constant finite energy solution $u_\p \in W(\R^2)$ to equation \eqref{TWc}, with $c = c(u_\p)$, and $p(u_\p) = \p$, such that $u_\p$ is solution to the minimization problem
$$E(u_\p) = E_{\min}(\p) = \inf \{ E(v), v \in W(\R^2), p(v) = \p \}.$$
\end{theorem}

\begin{remark}
In particular, Theorem \ref{dim2} shows that there exist travelling wave solutions of arbitrary small energy. This suggests that scattering in the energy space is not likely to hold.
\end{remark}

In dimension three, the existence result is somewhat different.

\begin{theorem}[\cite{BetGrSa1}]
\label{dim3}
Assume $N = 3$. There exists some constant $\p_0 > 0$ such that\\
i) For any $0 < \p < \p_0$,
$$E_{\min}(\p) = \inf \{ E(v), v \in W(\R^3), p(v) = \p \} = \sqrt{2} \p,$$
and the infimum is not achieved in $W(\R^3)$.\\
ii) For any $\p \geq \p_0$, there exists a non-constant finite energy solution $u_\p \in W(\R^3)$ to equation \eqref{TWc}, with $c = c(u_\p)$ and $p(u_\p) = \p$. Moreover, $E(u_{\p_0}) = E_{\min}(\p_0) = \sqrt{2} \p_0$, and for any $\p > \p_0$,
$$E(u_\p) = E_{\min}(\p) = \inf \{ E(v), v \in W(\R^3), p(v) = \p \} < \sqrt{2} \p.$$
\end{theorem}

Besides the existence of minimizers, our analysis yields also properties of the curve $E_{\min}$ as well as of the speed $c(u_\p)$. More precisely, we have

\begin{theorem}[\cite{BetGrSa1}]
\label {proemin}
i) Let $N = 2$ or $N = 3$. For any $\p, \q \geq 0$, we have the inequality
$$|E_{\min}(\p) - E_{\min}(\q)| \leq \sqrt{2} |\p - \q|,$$
i.e. the real-valued function $\p \mapsto E_{\min}(\p)$ is Lipschitz, with Lipschitz's constant $\sqrt{2}$. Moreover, it is positive and non-decreasing on $\R_+$.\\
ii) Assume $N = 2$. Then, the function $\p \mapsto E_{\min}(\p)$ is strictly concave.\\
iii) Assume $N = 3$. Then,
$$E_{\min}(\p) = \sqrt{2} \p, \ \forall 0 \leq \p \leq \p_0,$$
whereas $E_{\min}$ is strictly concave on $(\p_0, + \infty)$.
\end{theorem}

In particular, the function $E_{\min}$ is differentiable except possibly for a countable set of values. Its derivative at the points of differentiability is given by the speed $c(u_\p)$, which satisfies for any $\p > 0$,
$$0 < c(u_\p) < \sqrt{2}.$$
It remains an open problem to determine whether the curve $E_{\min}$ is differentiable or not. This question is related to the problem of uniqueness, up to symmetries, for the minimizer $u_\p$, which is completely open as well. As a matter of fact, uniqueness for any $\p > 0$ of the minimizer would lead to the differentiability of the full curve. We believe that, if at some point $E_{\min}$ were not differentiable, then there are at least two different minimizers with different speeds. In that case, the function $\p \mapsto c(u_\p)$ is not single valued. However, we can prove that it is a decreasing (possibly multivalued) function.

In dimension two, the function $E_{\min}$ has the following graph:

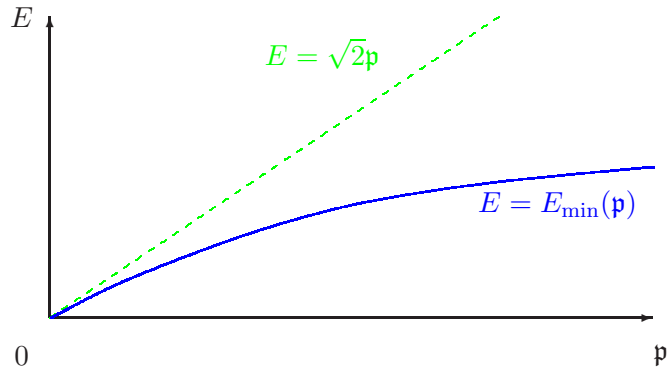
\begin{center}
\begin{picture}(90,50)(0,0)
\linethickness{0.2mm}
\put(10,10){\line(0,1){40}}
\put(10,50){\vector(0,1){0.12}}
\linethickness{0.2mm}
\put(10,10){\line(1,0){80}}
\put(90,10){\vector(1,0){0.12}}
{\color{green}
\linethickness{0.1mm}
\multiput(10,10)(1.64,1.1){37}{\multiput(0,0)(0.16,0.11){5}{\line(1,0){0.16}}}
}
{\color{blue}
\linethickness{0.2mm}
\qbezier(9,10)(9.29,10.17)(12.66,11.95)
\qbezier(12.66,11.95)(16.04,13.72)(19,15)
\qbezier(19,15)(26.49,18.14)(33.79,20.66)
\qbezier(33.79,20.66)(41.09,23.18)(49,25)
\qbezier(49,25)(60.4,27.33)(74.08,28.65)
\qbezier(74.08,28.65)(87.77,29.98)(89,30)
}
\put(5,5){\makebox(0,0)[cc]{$0$}}
\put(5,50){\makebox(0,0)[cc]{$E$}}
\put(90,5){\makebox(0,0)[cc]{$\p$}}
{\color{green}
\put(45,45){\makebox(0,0)[cc]{$E = \sqrt{2} \p$}}
}
{\color{blue}
\put(75,25){\makebox(0,0)[cc]{$E = E_{\min}(\p)$}}
}
\end{picture}
\end{center}

In dimension three, the graph of $E_{\min}$ has the following form:

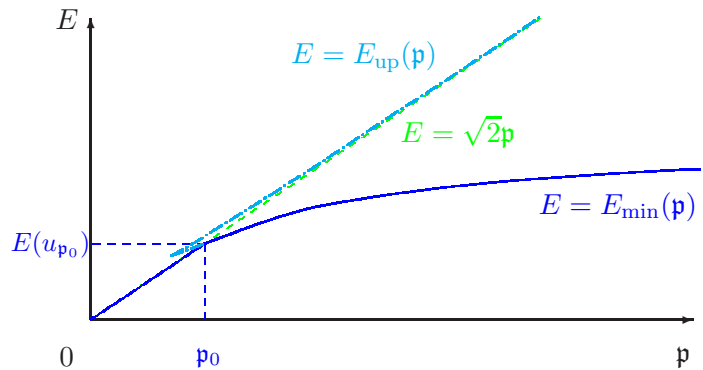
\begin{center}
\begin{picture}(90,50)(0,0)
\linethickness{0.2mm}
\put(10,10){\line(0,1){40}}
\put(10,50){\vector(0,1){0.12}}
\linethickness{0.2mm}
\put(10,10){\line(1,0){80}}
\put(90,10){\vector(1,0){0.12}}
{\color{green}
\linethickness{0.1mm}
\multiput(10,10)(1.64,1.1){37}{\multiput(0,0)(0.16,0.11){5}{\line(1,0){0.16}}}
}
{\color{blue}
\linethickness{0.1mm}
\multiput(24,10)(0,1.82){6}{\line(0,1){0.91}}
\linethickness{0.1mm}
\multiput(9,20)(2,0){8}{\line(1,0){1}}
}
{\color{blue}
\linethickness{0.2mm}
\multiput(7.4,10)(0.18,0.12){83}{\line(1,0){0.18}}
\linethickness{0.2mm}
\qbezier(22.4,20)(22.84,20.21)(27.9,22.07)
\qbezier(27.9,22.07)(32.96,23.94)(37.4,25)
\qbezier(37.4,25)(51.65,27.69)(68.75,28.86)
\qbezier(68.75,28.86)(85.86,30.03)(88.4,30)
}
{\color{cyan}
\linethickness{0.2mm}
\qbezier(21,20)(19.5,19.4)(19.5,19.4)
\qbezier(19,19.25)(18.9,19.2)(18.9,19.2)
\qbezier(18.4,19)(16.9,18.5)(16.9,18.5)
\qbezier(16.9,18.5)(18.4,19.3)(18.4,19.3)
\qbezier(18.9,19.6)(19,19.65)(19,19.65)
\qbezier(19.5,19.95)(21,20.85)(21,20.85)
\qbezier(21.5,21.2)(21.6,21.25)(21.6,21.25)
\qbezier(22.1,21.6)(23.6,22.55)(23.6,22.55)
\qbezier(24.1,22.85)(24.2,22.9)(24.2,22.9)
\qbezier(24.7,23.2)(26.2,24.15)(26.2,24.15)
\qbezier(26.7,24.45)(26.8,24.55)(26.8,24.55)
\qbezier(27.3,24.85)(28.8,25.75)(28.8,25.75)
\qbezier(29.3,26.1)(29.4,26.2)(29.4,26.2)
\qbezier(29.9,26.45)(31.4,27.4)(31.4,27.4)
\qbezier(31.9,27.75)(32,27.85)(32,27.85)
\qbezier(32.5,28.15)(34,29.1)(34,29.1)
\qbezier(34.5,29.4)(34.6,29.5)(34.6,29.5)
\qbezier(35.1,29.8)(36.6,30.75)(36.6,30.75)
\qbezier(37.1,31.1)(37.2,31.2)(37.2,31.2)
\qbezier(37.7,31.5)(39.2,32.5)(39.2,32.5)
\qbezier(39.7,32.8)(39.8,32.9)(39.8,32.9)
\qbezier(40.3,33.25)(41.8,34.25)(41.8,34.25)
\qbezier(42.3,34.55)(42.4,34.65)(42.4,34.65)
\qbezier(42.9,35)(44.4,35.95)(44.4,35.95)
\qbezier(44.9,36.25)(45,36.35)(45,36.35)
\qbezier(45.5,36.7)(47,37.7)(47,37.7)
\qbezier(47.5,38)(47.6,38.05)(47.6,38.05)
\qbezier(48.1,38.4)(49.6,39.35)(49.6,39.35)
\qbezier(50.1,39.75)(50.2,39.8)(50.2,39.8)
\qbezier(50.7,40.05)(52.2,41.05)(52.2,41.05)
\qbezier(52.7,41.4)(52.8,41.45)(52.8,41.45)
\qbezier(53.3,41.75)(54.8,42.75)(54.8,42.75)
\qbezier(55.3,43.1)(55.4,43.15)(55.4,43.15)
\qbezier(55.9,43.5)(57.4,44.5)(57.4,44.5)
\qbezier(57.9,44.85)(58,44.95)(58,44.95)
\qbezier(58.5,45.2)(60,46.2)(60,46.2)
\qbezier(60.5,46.6)(60.6,46.65)(60.6,46.65)
\qbezier(61.1,46.95)(62.6,47.95)(62.6,47.95)
\qbezier(63.1,48.35)(63.2,48.4)(63.2,48.4)
\qbezier(63.7,48.75)(65.2,49.75)(65.2,49.75)
\qbezier(65.7,50.1)(65.8,50.2)(65.8,50.2)
}
\put(3,5){\makebox(0,0)[cc]{$0$}}
\put(3,50){\makebox(0,0)[cc]{$E$}}
\put(85,5){\makebox(0,0)[cc]{$\p$}}
{\color{green}
\put(55,35){\makebox(0,0)[cc]{$E = \sqrt{2} \p$}}
}
{\color{blue}
\put(75,25){\makebox(0,0)[cc]{$E = E_{\min}(\p)$}}
}
{\color{cyan}
\put(40,45){\makebox(0,0)[cc]{$E = E_{\rm up}(\p)$}}
}
{\color{blue}
\put(18,5){\makebox(0,0)[cc]{{\small $\p_0$}}}
\put(-3,20){\makebox(0,0)[cc]{{\small $E(u_{\p_0})$}}}
}
\end{picture}
\end{center}

Our results are in full agreement with the corresponding figure given in \cite{JoneRob1}. In dimension three, the numerical value found in \cite{JoneRob1} for $\p_0$ is close to $80$. Jones and Roberts have also shown in \cite{JoneRob1}, mainly by numerical means, that in dimension three, the branch of solutions can be extended past the curve $E = \sqrt{2} \p$. Its representation in the $E$-$\p$ diagram bifurcates at some point, then is asymptotic to the curve $E = \sqrt{2} \p$ (see the curve $E_{\rm up}$ on the diagram above). At this stage, there is no mathematical proof of the existence of the upper branch of solutions. However, it is proved in \cite {BetGrSa1} that the slope of the curve at the point $(\p_0, \sqrt{2} \p_0)$ is strictly less than $\sqrt{2}$. This leaves some hope to use an implicit function theorem to construct the curve $E_{\rm up}$, at least near $(\p_0, \sqrt{2} \p_0)$.

\begin{remark}
Jones, Putterman and Roberts \cite{JoneRob1,JonPuRo1} conjectured the existence of some momentum $\p_1$ such that the minimizer $u_\p$ has vortices for $\p \geq \p_1$, and has no vortex otherwise. The numerical value found in \cite{JoneRob1} for $\p_1$ is close to $75$ in dimension three. There is no evidence of the existence of such a number $\p_1$. In particular, we do not know if $\p_1 > \p_0$ in dimension three.
\end{remark}

Taking advantage of the analyticity property of finite energy solutions to \eqref{TWc}, we may derive the following additional property of the minimizing solutions we obtained in Theorems \ref{dim2} and \ref{dim3}.

\begin{theorem}[\cite{BetGrSa1}]
\label{symetrie}
Let $N = 2$ or $N = 3$, $\p > 0$ and assume that $E_{\min}(\p)$ is achieved by $u_\p$. Then $u_\p$ is, up to a translation, axisymmetric around axis $x_1$. More precisely, there exists a function $\u_\p : \R \times \R_+$ such that
$$u_\p(x) = \u_\p(x_1, |x_\perp|), \ \forall x = (x_1, x_\perp) \in \R^N.$$
\end{theorem}

\subsection{The (KP I) transonic limit in dimension two}

In \cite{JoneRob1,JonPuRo1}, it is formally shown that, if $u_c$ is a two-dimensional solution to \eqref{TWc}, then, after a suitable rescaling, the function $1 - |u_c|^2$ converges, as the speed $c$ converges to $\sqrt{2}$, to a solitary wave solution to the Kadomtsev-Petviashvili equation \eqref{KP}, which writes
\renewcommand{\theequation}{KP I}
\begin{equation}
\label{KP}
\partial_t u + u \partial_1 u + \partial_1^3 u - \partial_1^{- 1} (\partial^2_2 u) = 0.
\end{equation}
Notice that the first terms correspond to the Korteweg-de Vries equation, whereas the last term is a transverse perturbation. As \eqref{GP}, equation \eqref{KP} is hamiltonian, with Hamiltonian given by
$$E_{KP}(u) = \frac{1}{2} \int_{\R^2} (\partial_1 u)^2 + \frac{1}{2} \int_{\R^2} (\partial_1^{-1}(\partial_2 u))^2 - \frac{1}{6} \int_{\R^2} u^3,$$
and the $L^2$-norm of $u$ is conserved as well. Solitary wave solutions $u(x, t) = w(x_1 - \sigma t, x_2)$ may be obtained in dimension two minimizing the Hamiltonian, keeping the $L^2$-norm fixed (see \cite{deBoSau3,deBoSau1}). The equation for the profile $w$ of a solitary wave of speed $\sigma = 1$ is given by
\renewcommand{\theequation}{\arabic{equation}}
\numberwithin{equation}{section}
\setcounter{equation}{0}
\begin{equation}
\label{SW}
\partial_1 w - w \partial_1 w - \partial_1^3 w + \partial_1^{- 1} (\partial_2^2 w) = 0.
\end{equation}
In contrast with \eqref{TWc}, the range of speeds is the positive axis. Indeed, for any given $\sigma > 0$, a solitary wave $w_\sigma$ of speed $\sigma$ is deduced from a solution $w$ to \eqref{SW} by the scaling
$$w_\sigma(x_1, x_2) = \sigma w(\sqrt{\sigma} x_1, \sigma x_2).$$
We term ground state, a solitary wave $w$ which minimizes the action $S$ given by
$$S(v) = E_{KP}(v) + \frac{1}{2} \int_{\R^2} v^2,$$
among all the solutions to \eqref{SW} (see \cite{deBoSau2} for more details). In dimension two, it is shown in \cite{deBoSau3} that $w$ is a ground state if and only if it minimizes the Hamiltonian $E_{KP}$ keeping the $L^2$-norm fixed. We denote $\boS_{KP}$, the action $S(w)$ of the ground states $w$.

The correspondence between \eqref{TWc} and \eqref{SW} is given as follows. Setting $\varepsilon \equiv \sqrt{2 - c^2}$ and $\eta_c \equiv 1 - |u_c|^2$, and performing the change of variables
$$N_\varepsilon(x) = \frac{6}{\varepsilon^2} \eta_c \Big( \frac{x_1}{\varepsilon}, \frac{\sqrt{2} x_2}{\varepsilon^2} \Big),$$
it is shown that $N_\varepsilon$ approximatively solves \eqref{SW} as $c$ converges to $\sqrt{2}$. The minimizing branch constructed in Theorem \ref{dim2} contains transonic solutions in the limit $\p \to 0$, as shown in the following proposition.

\begin{prop}
\label{T2bis}
Assume $N = 2$. There exist some constants $\p_1 > 0$, $K_0$ and $K_1$ such that we have the asymptotic behaviours
\begin{equation}
\label{estimE2}
\frac{48 \sqrt{2}}{\boS_{KP}^2} \p^3 - K_0 \p^4 \leq \sqrt{2} \p - E_{\min}(\p) \leq K_1 \p^3, \ \forall 0 \leq \p \leq \p_1.
\end{equation}
Let $u_\p$ be as in Theorem \ref{dim2}. Then, there exist some constants $\p_2 > 0$, $K_2 > 0$ and $K_3$ such that
\begin{equation}
\label{c-estim}
K_2 \p^2 \leq \sqrt{2} - c(u_\p) \leq K_3 \p^2, \ \forall 0 \leq \p < \p_2.
\end{equation}
\end{prop}

We next consider the map
$$N_\p = N_{\varepsilon_p},$$
where $\varepsilon_\p = \varepsilon(u_\p) = \sqrt{2 - c(u_\p)^2}$. In \cite{BetGrSa2}, we prove

\begin{theorem}[\cite{BetGrSa2}]
\label{convGPKP}
There exists a subsequence $(\p_n)_{n \in \N}$ tending to $0$, as $n \to + \infty$, a ground state $w$ of \eqref{KP}, and a universal constant $\gamma_0 > 0$ such that, for any $0 \leq \gamma < \gamma_0$, we have
$$N_{\p_n} \to w \ {\rm in} \ C^{0, \gamma}(\R^2), \ {\rm as} \ n \to + \infty.$$
\end{theorem}

\begin{remark}
There is an explicit solitary wave solution to \eqref{KP} of speed $1$, namely the so-called "lump" solution, which is written as
$$w_\ell(x_1, x_2) = 24 \frac{3 - x_1^2 + x_2^2}{(3 + x_1^2 + x_2^2)^2}.$$
It is usually conjectured that the "lump" solution is a ground state. It is also conjectured that the ground state solution is unique, up to the invariances of the problem. If this is the case, then the full family $(N_{\p})_{\p > 0} $ converges to $w_\ell$, as $\p \to 0$.
\end{remark}

\begin{remark}
If $u_c$ is a solution to \eqref{TWc} in dimension three, then it is also formally shown in \cite{IordSmi1,JoneRob1,JonPuRo1}, that the function $w_c$ defined by
$$w_c(x) = \frac{6}{\eps^2} \bigg( 1 - \Big|v_c \Big( \frac{x_1}{\varepsilon}, \frac{\sqrt{2} x_2}{\eps^2}, \frac{\sqrt{2} x_3}{\eps^2} \Big) \Big|^2 \bigg),$$
converges, as the parameter $\eps=\sqrt{2-c^2}$ converges to $0$, to a solitary wave solution $w$ to the three-dimensional Kadomtsev-Petviashvili equation \eqref{KP}, which writes
$$\partial_t u + u \partial_1 u + \partial_1^3 u - \partial_1^{- 1} (\partial_2^2 u + \partial_3^2 u) = 0.$$
In particular, the equation for the solitary wave $w$ is now written as
$$\partial_1 w - w \partial_1 w - \partial_1^3 w + \partial_1^{- 1} (\partial_2^2 w + \partial_3^2 w) = 0.$$
So far, there is no rigorous proof of the existence of a branch of solutions in the transonic limit. This branch of solution is conjectured in \cite{JoneRob1}, and represented on our graph above as the upper branch $E_{\rm up}$.
\end{remark}

\section{Related problems}
\label{potaufeu}

\subsection{Infinite energy solutions}

So far, we have restricted ourselves to finite energy solutions. Solutions with infinite energy are also of interest. In dimension two, a typical example is the stationary symmetric vortices for the Ginzburg-Landau equation (see e.g \cite{BergChe1}). In dimension three, for small speeds, Chiron \cite{Chiron2,Chiron3} established the existence of solutions having vorticity concentrated on helices. Such solutions are reminiscent to some flows for the incompressible Euler equation.

\subsection{Nonlinear Schr\"odinger flow past an obstacle}

\indent As mentioned in the introduction, the Gross-Pitaevskii equation provides a model describing superfluidity in the Hartree approximation. When impurities or an obstacle are present in the superfluid, one may model these new elements adding an external potential $V$ to the equation. As done by Hakim \cite{Hakim1} in the one-dimensional case, the flow of a superfluid past an obstacle moving with constant positive speed $c$ in direction $x_1$, may be described by the Gross-Pitaevskii equation with a coupling with an additional potential $V$, namely
\begin{equation}
\label{GP-V}
i \partial_t \Psi (x, t) = \Delta \Psi (x, t) + \Psi(x, t) \Big( 1 - |\Psi (x, t)|^2 - V(x_1 - c t, x_\perp) \Big), \ \big( x = (x_1, x_\perp), t \big) \in \R^N \times \R.
\end{equation}
As before, one adds the condition $|\Psi(x)| \to 1$, as $|x| \to + \infty$, to describe a superfluid which is at rest at infinity. In the frame of the moving obstacle, equation \eqref{GP-V} may be recast as
\begin{equation}
\label{GP-pot}
i \partial_t \Phi = i c \partial_1 \Phi + \Delta \Phi + \Phi \Big( 1 - |\Phi|^2 - V \Big),
\end{equation}
where we denote $\Phi(x, t) = \Psi(x_1 + c t, x_\perp, t)$. Stationary solutions $v$ to \eqref{GP-pot} satisfy the elliptic equation
\begin{equation}
\label{TW-pot}
i c \partial_1 v + \Delta v + v \Big( 1 - |v|^2 - V \Big) = 0.
\end{equation}
Equation \eqref{GP-pot} formally remains hamiltonian. For a superfluid at rest at infinity, the conserved Hamiltonian may be written as
\begin{equation}
\label{Fc}
F_c^V(\Phi) = \frac{1}{2} \int_{\R^N} |\nabla \Phi|^2 + \frac{1}{4} \int_{\R^N} (1 - |\Phi|^2)^2 - \frac{1}{2} \int_{\R^N} V (1 - |\Phi|^2) - \frac{c}{2} \int_{\R^N} \langle i \partial_1 \Phi, \Phi - 1 \rangle.
\end{equation}
The previous variational formulation leads us to restrict ourselves to solutions $v$ for which the following modified Ginzburg-Landau energy is finite, that is
$$E^V(v) \equiv \frac{1}{2} \int_{\R^N} |\nabla v|^2 + \frac{1}{4} \int_{\R^N} (1 - |v|^2)^2 - \frac{1}{2} \int_{\R^N} V (1 - |v|^2) < + \infty.$$

In \cite{Hakim1}, a formal and numerical analysis of stationary solutions to \eqref{GP-pot} was initiated in the one-dimensional case, whereas in \cite{Maris3}, still in dimension one, Maris established the existence of solutions to \eqref{TW-pot} minimizing the Hamiltonian $F_c^V$, assuming the potential $V$ is a bounded measure with small total variation compared to the speed $c$. Assuming additionally that $V$ is compactly supported, Maris \cite{Maris3} also exhibited excited states.

In higher dimensions, the flow of a superfluid past an obstacle has been widely investigated both formally and numerically (see, for instance, \cite{FriPoRi1,PomeRic1,PomeRic2,SrivSto1,AdaJaMc1,AdJaMcW1,BracHue1}). In \cite{Tarquin2}, Tarquini proved that, under suitable assumptions on the localized potential $V$, there are no finite energy supersonic solutions to \eqref{TW-pot}, that is with speed $c > \sqrt{2}$.
He also computed an explicit bound for the solutions to \eqref{TW-pot}, and described their regularity and convergence at infinity.

To our knowledge, the existence problem of finite energy solutions to \eqref{TW-pot} has not been addressed yet for a general class of potentials, although some special cases have been considered in \cite{AftaBla1,Aftalio1}. Therefore, we would like to provide in this paper an existence result in dimension two.

\begin{theorem}
\label{tripot}
Let $0 \leq c < \sqrt{2}$ and $V \in L^2(\R^2)$. There exists a constant $K(c)>0$ depending only on $c$, such that, if
\begin{equation}
\label{pas-de-pot}
\int_{\R^2} |V|^2 \leq K(c)^2,
\end{equation}
then equation \eqref{TW-pot} has at least one finite energy solution $u_c$.
\end{theorem}

The solutions $u_c$ provided by Theorem \ref{tripot} are constructed minimizing the Hamiltonian $F_c^V$ given by \eqref{Fc} locally near the constant solutions to \eqref{TWc}. As a matter of fact, in the case the potential is $0$, the solutions we construct correspond to the constant solutions to \eqref{TWc}. So far, we have not tried to construct excited states (as done in the one-dimensional case by Maris \cite{Maris3}). For a zero potential, these solutions would correspond to the non-constant two-dimensional solutions to \eqref{TWc} we have discussed before.

\begin{remark}
Though we did not address this question here, we believe that the solutions $u_c$ provided by Theorem \ref{tripot} are local minimizers of $F_c^V$ on suitable variational spaces (see \cite{BetGrSa1}, where a similar question is addressed for the equation without potential, and \cite{Tarquin2} for the choice of the suitable variational space, in case $V$ has compact support). In particular, the argument used to prove Theorem \ref{tripot} does not establish that either $F_c^V(u_c)$ or $p(u_c)$ are finite.
\end{remark}

As mentioned above, we construct the solutions $u_c$ by locally minimizing the Hamiltonian $F_c^V$ given by \eqref{Fc}. When the potential $V$ belongs to $L^2(\R^2)$, $F_c^V$ is well-defined and of class $C^1$ on the space $\{ 1 \} + H^1(\R^2)$. However, the construction of minimizers for $F_c^V$ presents some difficulty in view of the lack of compactness. To circumvent this difficulty, we introduce, as in \cite{BetOrSm1,BetGrSa1}, minimizing problems on expanding tori, for which the existence of minimizers presents no major difficulty, and then pass to the limit when the size of the torus tends to infinity. More precisely, we introduce the flat torus defined by
$$\T_n \simeq \Omega_n \equiv [- \pi n, \pi n]^2,$$
for any $n \in \N^*$ (with opposite faces identified), and the space
$$X_n = H^1(\T_n, \C) \simeq H^1_{\rm per}(\Omega_n, \C)$$
of $2 \pi n$-periodic $H^1$-functions. We define the energy $E_n^V$ and the momentum $p_n$ on $X_n$ by
\begin{equation}
\label{En-pot}
E_n^V(v) = \frac{1}{2} \int_{\Omega_n} |\nabla v|^2 + \frac{1}{4} \int_{\Omega_n} (1 - |v|^2)^2 - \frac{1}{2} \int_{\Omega_n} V (1 - |v|^2),
\end{equation}
and
$$p_n(v) = \frac{1}{2} \int_{\Omega_n} \langle i \partial_1 v , v - 1 \rangle,$$
so that the Hamiltonian $F_{c, n}^V$ on $X_n$ is given by
\begin{equation}
\label{Fcn-pot}
F_{c, n}^V(v) = E_n^V(v) - c p_n(v).
\end{equation}
Notice that the function $V$ in \eqref{En-pot} and \eqref{Fcn-pot} denotes the $2 \pi n$ -periodic restriction of $V$ to $\Omega_n$. For given $\Lambda > 0$, we introduce the open set $\boE_n(\Lambda)$ defined by
$$\boE_n(\Lambda) = \bigg\{ v \in X_n, \ {\rm s.t.} \ E_n(v) \equiv \frac{1}{2} \int_{\Omega_n} |\nabla v|^2 + \frac{1}{4} \int_{\Omega_n} (1 - |v|^2)^2 < \Lambda \bigg\},$$
which is a neighbourhood of the set of constant functions, and consider the minimization problem
\begin{equation}
\label{Pn-pot}
\boF_{c, n}^V(\Lambda) \equiv \underset{v \in \boE_n(\Lambda)}{\inf} \Big( F_{c,n}^V(v) \Big).
\end{equation}
Our purpose is to prove, that for a suitable choice of $\Lambda = \Lambda(c)$ with respect to the speed $c$, independent of $n$, the previous minimization is achieved. The main ingredient in the proof of Theorem \ref{tripot}, which allows in particular to prove that $F_c^V$ is bounded from below on any set of functions with sufficiently small energy, is the following estimate.

\begin{lemma}
\label{levelcurve}
Let $\delta_0 > 0$ and $\delta_1 > 0$ such that
$$\delta_c = \frac{c}{\sqrt{2}} < \delta_0 < \delta_1 < 1,$$
and consider a function $v \in \boE_n(\Lambda)$. Then, there exists some integer $n_0$, and some positive constant $\Lambda_0$, depending only on $\delta_0$ and $\delta_1$, such that
\begin{equation}
\label{colorado}
c |p_n(v)| \leq \frac{\delta_c}{\delta_0} E_n(v) + \frac{c \delta_1}{\sqrt{\pi} (1 - \delta_1^2)(\delta_1 - \delta_0)} E_n(v)^\frac{3}{2} + \frac{2 \sqrt{2} c}{(1 - \delta_1^2)(\delta_1 - \delta_0)^2 \delta_0} E_n(v)^2,
\end{equation}
for any $n > n_0$, and any $\Lambda < \Lambda_0$.
\end{lemma}

\begin{proof}
The proof of Lemma \ref{levelcurve} is reminiscent of the proof of Lemma \ref{colisee}. When $E_n(v)$ is small enough, the modulus of $v$ is close to $1$. Assuming that we may construct a lifting $v = \varrho \exp i \varphi $ of $v$, we compute
\begin{equation}
\label{toutfaux}
c |p_n(v)| \approx \frac{c}{2} \bigg| \int_{\Omega_n} \big( 1 - \varrho^2 \big) \partial_1 \varphi \bigg| \leq \frac{c}{\sqrt{2}} \int_{\Omega_n} \frac{e(v)}{\varrho} \approx \delta_c \int_{\Omega_n} e(v) = \delta_c E_n(v).
\end{equation}
Lemma \ref{levelcurve} reduces to estimate any error terms in \eqref{toutfaux} to obtain \eqref{colorado}.

Therefore, we consider some function $v \in \boE_n(\Lambda)$, assuming, up to a density argument, that $v$ is smooth on $\Omega_n$. In order to estimate how the modulus of $v$ is close to $1$, and to construct a lifting of $v$, we study the level sets
$$\Omega(\xi) = \big\{ x \in \Omega_n, \ {\rm s.t.} \ |v(x)| = \xi \big\},$$
for any given $\xi \in \R_+$. By Sard's lemma, the sets $\Omega(\xi)$ are one-dimensional, smooth submanifolds, for almost any $\xi \in \R_+$, so that their Hausdorff's length
$$L(\xi) = \boH^1 \big( \Omega(\xi) \big),$$
is finite for almost any $\xi \in \R_+$. Moreover, the coarea formula gives
$$\int_{\delta_0}^{\delta_1} \bigg( \int_{\Omega(\xi)} |\nabla v| \bigg) d\xi \leq \int_{\delta_0 < |x| < \delta_1} |\nabla v(x)|^2 dx \leq 2 E_n(v),$$
and
$$\int_{\delta_0}^{\delta_1} L(\xi) d\xi \leq \int_{\delta_0 < |x| < \delta_1} |\nabla v(x)| dx \leq \sqrt{2} E_n(v)^\frac{1}{2} |U(\delta_1)|^\frac{1}{2},$$
where we denote $U(\xi) = \{ x \in \Omega_n, \ {\rm s.t.} \ |v(x)| < \xi \}$, for any $\xi \in \R_+$. Since
$$|U(\delta_1)| \leq \frac{4}{(1 - \delta_1^2)^2} E_n(v),$$
we are led to
$$\int_{\delta_0}^{\delta_1} L(\xi) d\xi \leq \frac{2 \sqrt{2}}{1 -
 \delta_1^2} E_n(v).$$
The mean value inequality then provides the existence of $\delta_0 < \xi_0 < \delta_1$ such that the length of the one-dimensional submanifold $\Omega(\xi_0)$ verifies
\begin{equation}
\label{alaska}
L(\xi_0) \leq \frac{2 \sqrt{2}}{(1 - \delta_1^2)(\delta_1 - \delta_0)} E_n(v) \leq \frac{2 \sqrt{2} \Lambda}{(1 - \delta_1^2)(\delta_1 - \delta_0)},
\end{equation}
and
\begin{equation}
\label{illinois}
\int_{\Omega(\xi_0)} |\nabla v| \leq \frac{2}{\delta_0 - \delta_1} E_n(v) \leq \frac{2 \Lambda}{\delta_0 - \delta_1}.
\end{equation}
We first deduce from \eqref{alaska} that, up to some new unfolding of the torus $\T_n$, $\Omega(\xi_0)$ is, for $n \geq \frac{2 \sqrt{2} \Lambda}{\pi (1 - \delta_1^2)(\delta_1 - \delta_0)}$, a collection of smooth, connected, closed curves $(\gamma_i)_{i \in I}$ of $\Omega_n$. Moreover, since
$$|U(\xi_0)| \leq |U(\delta_1)| \leq \frac{4}{(1 - \delta_1^2)^2} E_n(v) \leq \frac{4 \Lambda}{(1 - \delta_1^2)^2},$$
we may assume that the complementary of $U(\xi_0)$ contains the boundary of $\Omega_n$ for any $n \geq \frac{2 \sqrt{\Lambda}}{\pi (1 - \delta_1^2)}$.

It now follows from \eqref{illinois} that
$$\underset{x \in \gamma_i}{\sup} \bigg| \frac{v(x)}{|v(x)|} - \frac{v(y)}{|v(y)|} \bigg| \leq \frac{1}{\xi_0} \int_{\Omega(\xi_0)} |\nabla v| \leq \frac{2 \Lambda}{\delta_0 (\delta_1 - \delta_0)}, \ \forall y \in \gamma_i,$$
so that, assuming that $\Lambda < \delta_0 (\delta_1 - \delta_0)$, the topological degree of $\frac{v}{|v|}$ on each curve $\gamma_i$ is equal to $0$. This provides the existence of a smooth lifting $v = \varrho \exp i \varphi$ of $v$ on the complementary of $U(\xi_0)$, so that, using the periodicity of $v$, we can compute
\begin{equation}
\label{idaho}
p_n(v) = \frac{1}{2} \int_{\Omega_n} \langle i \partial_1 v, v \rangle = \frac{1}{2} \int_{U(\xi_0)} \langle i \partial_1 v, v \rangle - \frac{1}{2} \int_{\Omega_n \setminus U(\xi_0)} \varrho^2 \partial_1 \varphi.
\end{equation}
The first term in the right-hand side of \eqref{idaho} can be estimated by
$$\bigg| \frac{1}{2} \int_{U(\xi_0)} \langle i \partial_1 v, v \rangle \bigg| \leq \frac{|\xi_0|}{2} \int_{U(\xi_0)} |\partial_1 v| \leq \frac{\delta_1}{2} |U(\xi_0)|^\frac{1}{2} \bigg( \int_{U(\xi_0)} |\partial_1 v|^2 \bigg)^\frac{1}{2},$$
so that, using the isoperimetric inequality
$$4 \pi |U(\xi_0)| \leq L(\xi_0)^2,$$
and \eqref{illinois}, we are led to
\begin{equation}
\label{alabama}
\bigg| \frac{1}{2} \int_{U(\xi_0)} \langle i \partial_1 v, v \rangle \bigg| \leq \frac{\delta_1}{\sqrt{\pi} (1 - \delta_1^2)(\delta_1 - \delta_0)} E_n(v)^\frac{3}{2}.
\end{equation}
On the other hand, the second term in the right-hand side of \eqref{idaho} is written as
\begin{equation}
\label{dakota}
- \frac{1}{2} \int_{\Omega_n \setminus U(\xi_0)} \varrho^2 \partial_1 \varphi = \frac{1}{2} \int_{\Omega_n \setminus U(\xi_0)} (1 - \varrho^2) \partial_1 \varphi - \frac{1}{2} \sum_{i \in I} \int_{\gamma_i} (\varphi - \varphi_i) \nu_1,
\end{equation}
where $\nu_1$ denotes the first component of the outward normal to $\gamma_i$, and
$$\varphi_i = \frac{1}{|\gamma_i|} \int_{\gamma_i} \varphi.$$
Using the result of Lemma \ref{colisee}, we obtain
$$\bigg| \int_{\Omega_n \setminus U(\xi_0)} (1 - \varrho^2) \partial_1 \varphi \bigg| \leq \sqrt{2} \int_{\Omega_n \setminus U(\xi_0)} \frac{e(v)}{\varrho} \leq \frac{\sqrt{2}}{\xi_0} E_n(v),$$
whereas
$$\bigg| \sum_{i \in I} \int_{\gamma_i} (\varphi - \varphi_i) \nu_1 \bigg| \leq \sum_{i \in I} \int_{\gamma_i} \bigg( \int_{\gamma_i} |\nabla \varphi| \bigg) \leq \frac{L(\xi_0)}{\xi_0} \int_{\Omega(\xi_0)} |\nabla v|,$$
so that using \eqref{alaska}, \eqref{illinois} and \eqref{dakota}, we are led to
\begin{equation}
\label{connecticut}
\bigg| \frac{1}{2} \int_{\Omega_n \setminus U(\xi_0)} \varrho^2 \partial_1 \varphi \bigg| \leq \frac{1}{\sqrt{2} \delta_0} E_n(v) + \frac{2 \sqrt{2}}{\delta_0 (1 - \delta_1^2) (\delta_1 - \delta_0)^2} E_n(v)^2.
\end{equation}
Inequality \eqref{colorado} finally follows from \eqref{idaho}, \eqref{alabama} and \eqref{connecticut}, provided that
\begin{equation}
\label{portorico}
\Lambda < \Lambda_0 \equiv \delta_0 (\delta_1 - \delta_0), \ {\rm and} \ n \geq n_0 \equiv \max \bigg\{ \frac{2 \sqrt{2} \delta_0}{\pi ( 1 - \delta_1^2)}, \frac{2 \sqrt{\delta_0 (\delta_1 - \delta_0)}}{\pi ( 1 - \delta_1^2)} \bigg\}.
\end{equation}
\end{proof}

As a consequence of Lemma \ref{levelcurve}, we have

\begin{lemma}
\label{missouri}
Let $0 < c < \sqrt{2}$ and $n \geq n_0$, where $n_0$ is the integer of Lemma \ref{levelcurve}. There exists positive constants $K(c)$, $\Lambda(c)$ and $\kappa(c) < \Lambda(c)$ such that, if $\| V \|_{L^2(\R^2)} \leq K(c)$, then
\begin{equation}
\label{aubassin}
F_{c, n}^V(v) \geq \boF_{c, n}^V(\kappa(c)) ,
\end{equation}
for any function $v \in \boE_n(\Lambda(c))$. Moreover,
\begin{equation}
\label{negative}
- \Lambda(c) \leq \boF_{c, n}^V(\Lambda(c)) \leq 0.
\end{equation}
\end{lemma}

 \begin{proof}
The proof reduces to compute an upper and a lower bound for the function $F_{c, n}^V$ given by
\begin{equation}
\label{michigan}
F_{c, n}^V(v) = E_n(v) - \frac{1}{2} \int_{\T_n} V (1 - |v|^2) - c p_n(v),
\end{equation}
on the set $\boE_n(\Lambda)$, which amounts to estimate the second and the third term of the right-hand side of \eqref{michigan}. We use Lemma \ref{levelcurve} to provide an upper bound for the scalar momentum $p_n(v)$. Therefore, we consider numbers $\delta_0$, $\delta_1$, $\delta_2$ and $\mu$ such that
\begin{equation}
\label{guam}
\delta_c = \frac{c}{\sqrt{2}} < \delta_2 < \delta_0 < \delta_1 < 1, \ {\rm and} \ \frac{\delta_c}{\delta_2} < \mu < 1,
\end{equation}
and assume that $n \geq n_0$ and $\Lambda < \Lambda_0$, where $n_0$ and $\Lambda_0$ are provided by Lemma \ref{levelcurve}. By Lemma \ref{levelcurve}, we have
$$|c p_n(v)| \leq \frac{\delta_c}{\delta_2} E_n(v),$$
provided that $E_n(v) < \Lambda_1(c)$, where $\Lambda_1(c)$ denotes some positive constant depending only on $c$. On the other hand, we compute
$$\bigg| \frac{1}{2} \int_{\Omega_n} V (1 - |v|^2) \bigg| \leq \frac{1}{4 (1 - \mu)} \int_{\Omega_n} V^2 + (1 - \mu) E_n(v),$$
so that
\begin{equation}
\label{oregon}
\Big( \mu - \frac{\delta_c}{\delta_2} \Big) E_n(v) - \frac{1}{4 (1 - \mu)} \int_{\R^2} V^2 \leq F_{c, n}^V(v) \leq \Big( 1 + \mu + \frac{\delta_c}{\delta_2} \Big) E_n(v) + \frac{1}{4 (1 - \mu)} \int_{\R^2} V^2.
\end{equation}
We now fix
$$\Lambda(c) = \min \{ \Lambda_0, \Lambda_1(c) \} > 0,$$
and choose some positive number $K(c)$ and $\kappa(c)$ such that
\begin{equation}
\label{wisconsin}
\big( 1 + 2 \mu \big) \kappa(c) \leq \Lambda(c), \ {\rm and} \ K(c)^2 < 4 \big( 1 - \mu \big) \Big( \mu - \frac{\delta_c}{\delta_2} \Big) \kappa(c),
\end{equation}
so that, assuming that $\| V \|_{L^2(\R^2)} \leq K(c)$, and $n \geq n_0$, \eqref{oregon} provides
$$F_{c, n}^V(v) \geq \Big( \mu - \frac{\delta_c}{\delta_2} \Big) \kappa(c) - \frac{1}{4 (1 - \mu)} K(c)^2 > 0 = F_{c, n}^V(1) \geq \boF_{c, n}^V(\kappa(c)),$$
for any function $v \in X_n$ such that $\kappa(c) \leq E_n(v) \leq \Lambda(c)$. Hence, we are led to \eqref{aubassin}, that is
\begin{equation}
\label{auxernests}
\boF_{c, n}^V(\kappa(c)) = \boF_{c, n}^V(\Lambda(c)) \leq 0.
\end{equation}
On the other hand, \eqref{oregon} also yields, using \eqref{wisconsin},
$$|F_{c, n}^V(v)| \leq \Big( 1 + \mu + \frac{\delta_c}{\delta_2} \Big) \kappa(c) + \frac{1}{4 (1 - \mu)} K(c)^2 \leq \big( 1 + 2 \mu \big) \kappa(c) \leq \Lambda(c),$$
for any function $v \in \boE_n(\kappa(c))$, so that \eqref{negative} follows from \eqref{auxernests}.
\end{proof}

We next have

\begin{lemma}
\label{Petit-pot}
Let $0 < c < \sqrt{2}$ and $n \geq n_0$, where $n_0$ is the integer of Lemma \ref{levelcurve}, and let $K(c)$ and $\Lambda(c)$ be as in Lemma \ref{missouri}. Assume moreover that condition \eqref{pas-de-pot} holds for $K(c)$ and $V$. Then, the minimization problem \eqref{Pn-pot} for $\Lambda = \Lambda(c)$ has a solution $u_c^n$ in $\boE_n(\Lambda(c))$. In particular, $u_c^n$ is solution to \eqref{TW-pot} on $\Omega_n$, and verifies
\begin{equation}
\label{petit-pot}
E_n(u_c^n) \leq \kappa(c)<\Lambda(c), \ {\rm and} \ - \Lambda(c) \leq F_{c, n}^V(u_c^n) = \boF_{c, n}^V(\Lambda(c)) \leq 0.
\end{equation}
\end{lemma}

\begin{proof}
Lemma \ref{Petit-pot} follows from standard variational arguments. For given $n \in \N$ and $\Lambda > 0$, $F_{c, n}^V$ is bounded from below on $\boE_n(\Lambda)$. Indeed, there exists some constant $K(n, \Lambda)$ such that
$$\int_{\Omega_n} |\nabla v|^2 + \int_{\Omega_n} |v|^4 \leq K(n, \Lambda),$$
for any function $v \in \boE_n(\Lambda)$, so that
$$F_{c, n}^V(v) \geq - \sqrt{\Lambda} \Big( \| V \|_{L^2(\T_n)} + c \| v - 1\|_{L^2(\T_n)} \Big) \geq - K(n, \Lambda).$$
Next, consider a minimizing sequence $(w_m)_{m \in \N}$ for $F_{c, n}^V$ on $\Lambda(c)$, and assume $n \geq n_0$. In view of \eqref{aubassin}, we have
$$F_{c, n}^V(w_m) \geq \boF_{c, n}^V(\kappa(c)),$$
so that we may assume, since $0 < \kappa(c) < \Lambda(c)$, and hence $\boE_n(\kappa(c)) \subset \boE_n(\Lambda(c))$, that $w_m \in \boE_n(\kappa(c))$ for any $m \in \N$, or equivalently,
$$E_n(w_n) \leq \kappa(c) < \Lambda(c).$$
Using Rellich's compactness theorem, there exists, up to some subsequence, a function $u_c^n \in X_n$ such that
$$\nabla w_m \rightharpoonup \nabla u_c^n \ {\rm in} \ L^2(\Omega_n), \ {\rm and} \ w_m \to u_c^n \ {\rm in} \ L^4(\Omega_n), \ {\rm as} \ m \to + \infty,$$
so that
$$F_{c, n}^V(u_c^n) \leq \underset{m \to + \infty}{\liminf} \Big( F_{c, n}^V(w_m) \Big) = \boF_{c, n}^V(\Lambda).$$
The same argument also provides
$$E_n(u_c^n) \leq \underset{m \to + \infty}{\liminf} \Big( E_n(w_m) \Big) \leq \kappa(c)<\lambda(c).$$
when $\| V \|_{L^2(\R^2)} \leq K(c)$. Hence, the minimization problem \eqref{Pn-pot} for $\Lambda = \Lambda(c)$ is attained by $u_c^n$ in the open set $\boE_n(\Lambda(c))$. Since $F_{c, n}^V$ is of class $C^1$ on $X_n$, $u_c^n$ is a critical point of $F_{c, n}^V$, that is a solution to \eqref{TW-pot} on $\Omega_n$.
\end{proof}

\begin{proof}[Proof of Theorem \ref{tripot}]
Passing to the limit $n \to + \infty$ requires to prove some compactness for the sequence $(u_c^n)_{n \geq n_0}$. We first provide a uniform estimate for $(u_c^n)_{n \geq n_0}$ using \eqref{petit-pot}, and standard elliptic estimates. By \eqref{petit-pot},
$$\int_{B(x_0, 2)} \Big( |\nabla u_c^n|^2 + |u_c^n|^4 \Big) \leq \Lambda(c),$$
where we denote $B(x_0, 2) = \{ x \in \Omega_n, |x - x_0| < 2 \}$ for any $x_0 \in \Omega_n$. Using H\"older's inequality, we are led to
$$\int_{B(x_0, 2)} \Big( \big| u_c^n (1 - |u_c^n|^2) \big|^\frac{4}{3} + \big| V u_c^n \big|^\frac{4}{3} + \big| i c \partial_1 u_c^n \big|^\frac{4}{3} \Big) \leq \Lambda(c),$$
so that, by equation \eqref{TW-pot}, and standard elliptic estimates, the sequence $(u_c^n)_{n \geq n_0}$ is bounded in $W^{2, \frac{4}{3}}(B(x_0, 1))$. Therefore, by Sobolev's embedding theorem, there exists some positive constant $\Lambda(c)$ such that
\begin{equation}
\label{texas}
\| u_c^n \|_{L^\infty(\Omega_n)} \leq \Lambda(c), \ \forall n \geq 1.
\end{equation}
Moreover, \eqref{petit-pot} also gives
\begin{equation}
\label{arizona}
\| \nabla u_c^n \|_{L^2(\Omega_n)} \leq \Lambda(c), \ \forall n \geq 1.
\end{equation}
Since any compact set $K$ of $\R^2$ is included in $\Omega_n$ for $n$ sufficiently large, we then construct, using \eqref{texas}, \eqref{arizona}, Rellich's compactness theorem, and a diagonal argument, a function $u_c \in H^1_{\rm loc}(\R^2)$, such that, up to some subsequence,
\begin{equation}
\label{nevada}
\nabla u_c^n \rightharpoonup \nabla u_c \ {\rm in} \ L^2(K), \ {\rm and} \ u_c^n \to u_c \ {\rm in} \ L^q(K), \ {\rm as} \ n \to + \infty,
\end{equation}
for any compact set $K$ of $\R^2$, and any $1 \leq q < + \infty$. Since $u_c^n$ is solution to \eqref{TW-pot} on $\Omega_n$, convergences \eqref{nevada} first yield that $u_c$ is solution to \eqref{TW-pot} on $\R^2$. Moreover, they also provide
$$\frac{1}{2} \int_K |\nabla u_c|^2 + \frac{1}{4} \int_K (1 - |u_c|^2)^2 \leq \underset{n \to +\infty}{\liminf} \bigg( \frac{1}{2} \int_K |\nabla u_c^n|^2 + \frac{1}{4} \int_K (1 - |u_c^n|^2)^2 \bigg) \leq \Lambda(c),$$
so that by monotone convergence,
$$E(u_c) \leq \Lambda(c).$$
Hence,
$$E^V(u_c) \leq E(u_c) + \| V \|_{L^2(\R^2)} E(u_c)^\frac{1}{2} < + \infty,$$
so that $u_c$ is a finite energy solution to \eqref{TW-pot} on $\R^2$.
\end{proof}

\begin{remark}
Notice that it follows from \eqref{portorico}, \eqref{guam} and \eqref{wisconsin} that the constant $K(c)$ in Theorem \ref{tripot} tends to $0$, as $c \to \sqrt{2}$, so that the considered potential $V$ must be smaller and smaller as $c \to \sqrt{2}$.
\end{remark}

\subsection{Nucleation by impurities}

An interesting generalization of \eqref{GP} (see \cite{GranRob1}) concerns the modelling of vortex nucleation by an impurity, e.g. an electron. In the Hartree approximation, the equations governing the one-particle wave function of the condensate $\psi$, and the wave function of the impurity $\phi$, are a pair of coupled equations
\begin{equation}
\label{impur1}
\begin{split}
i\hbar \partial_t \Psi & = -\frac{\hbar^2}{2 M} \Delta \Psi + (U_0 |\Phi|^2 + V_0 |\Psi|^2 - E) \Psi,\\
i\hbar \partial_t \Phi & = - \frac{\hbar^2}{2 \mu} \Delta \Phi +(U_0 |\Psi|^2 - E_e) \Phi,
\end{split}
\end{equation}
where $M$ and $E$ are the mass and single-particle energy for the bosons, and $\mu$ and $E_e$ are the mass and energy of the impurity, while $U_0$ (resp. $V_0$) denotes the mass of the Dirac interaction potentials between bosons and impurity (resp. bosons).\\
In dimensionless variables, \eqref{impur1} is written as (see \cite{GranRob1})
\begin{equation}
\label{impur2}
\begin{split}
i \partial_t \Psi & = -\Delta \Psi + \frac{1}{\varepsilon^2} \Big(\frac{1}{\varepsilon^2} |\Phi|^2 + |\Psi|^2 - 1 \Big) \Psi,\\
i \delta \partial_t \Phi & = -\Delta \Phi + \frac{1}{\varepsilon^2} \Big( q^2 |\Psi|^2 - \varepsilon^2 k^2 \Big) \Phi.
\end{split}
\end{equation}
Here, $\delta = \frac{\mu}{M}$ is the ratio of the mass of the impurity over the boson mass ($\delta \ll 1$), $q^2 = \delta \frac{U_0}{V_0}$, $k^{10} = \frac{\mu^5 E_e^5 U_0^2}{2 \pi^2 M^2 E^4 \hbar^6}$ is a dimensionless measure of the single-particle impurity energy, and $\varepsilon^{10} = \frac{2 \pi^2 \hbar^6}{E M^3 U_0^2}$ is a dimensionless constant, which in applications is about $0.2$.\\
Assuming that we are in a frame in which the condensate is at rest at infinity yields the formal boundary conditions
$$|\psi(x)| \to 1, \ {\rm and} \ \phi(x) \to 0, \ {\rm as} \ |x| \to + \infty.$$
Similarly to the Gross-Pitaevskii equation, \eqref{impur2} has a formally conserved Hamiltonian,
$$E(\Psi, \Phi) = \frac{1}{2} \int_{\R^N} \bigg( \varepsilon^4 |\nabla \Psi|^2 + \frac{\varepsilon^2}{q^2} |\nabla\Phi|^2 + \frac{\varepsilon^2}{2} \Big( 1 - |\Psi|^2 \Big)^2 + |\Phi|^2 |\Psi|^2 - \frac{\varepsilon^2 k^2}{q^2} |\Phi|^2 \bigg) \equiv \int_{\R^N} e(\Psi, \Phi).$$
It moreover conserves
$${\rm m}(\Phi) = \int_{\R^N} |\Phi|^2.$$
The finite energy travelling waves of \eqref{impur2} have been formally investigated in \cite{GranRob1}. No rigorous existence results seem to be known so far in dimensions two and three, but Bouchel \cite{Bouchel0}, extending the analysis of \cite{Graveja2,Graveja3} for \eqref{GP}, has proved decay estimates for finite energy travelling waves of \eqref{impur2}, together with the non-existence of supersonic travelling waves. On the other hand, Maris \cite{Maris4} has proved that one-dimensional travelling waves to \eqref{impur2} exist if and only if their velocity is less than the sound velocity at infinity, and that in this case, the set of travelling waves contains global subcontinua in appropriate Sobolev spaces.

\section{Conclusion}

Although travelling waves have been widely discussed in the physical literature, the field is still widely open for mathematical investigation. In the course of our discussion, we already mentioned a number of open problems and research directions: differentiability of the curve $E_{\min}$, existence of solutions to \eqref{TWc} for all values of speeds $0 < c < \sqrt{2}$, uniqueness, up to the invariances, of constrained minimizers as well as of ground states of \eqref{KP}, existence of the upper branch $E_{\rm up}$ in dimension three, and its possible transonic limit, non-existence of sonic travelling waves in dimension three, similar questions for the related problems above... To complete this survey, we wish to emphasize the question of stability: even the orbital stability of constrained minimizers presents serious yet unsolved difficulties in dimensions two and three.

\bibliographystyle{plain}
\bibliography{Bibliogr}

\end{document}